\documentclass[preprint,
3p,
times
]{elsarticle}
\usepackage{ifdraft} 

\makeatletter
\renewcommand\subsection{\@startsection{subsection}{2}{\z@}%
           {12\p@ \@plus 6\p@ \@minus 3\p@}%
           {3\p@ \@plus 6\p@ \@minus 3\p@}%
           {\normalsize\bfseries}}
\renewcommand\subsubsection{\@startsection{subsubsection}{3}{\z@}%
           {12\p@ \@plus 6\p@ \@minus 3\p@}%
           {\p@}%
           {\normalsize\bfseries}}
\makeatother

\makeatletter
\def\ps@pprintTitle{%
 \let\@oddhead\@empty
 \let\@evenhead\@empty
 \def\@oddfoot{}%
 \let\@evenfoot\@oddfoot}
\makeatother

\usepackage[american]{babel} 
\usepackage[babel]{csquotes} 
\usepackage{mathtools}
\usepackage[OMLmathsfit]{isomath} 
\usepackage{booktabs}        
\usepackage{tabularx}             
  \newcolumntype{R}{>{\raggedleft\arraybackslash}X}  \newcolumntype{L}{>{\raggedright\arraybackslash}X}  \newcolumntype{C}{>{\centering\arraybackslash}X}
\usepackage{siunitx}                 
\usepackage{float} 

\usepackage{hyperref}
\usepackage{cleveref} 
  \crefformat{algorithm}{Alg.~#2#1#3}        \crefname{algorithm}{Alg.}{Algs.}
  \Crefformat{algorithm}{Algorithm~#2#1#3}   \Crefname{algorithm}{Algorithm}{Algorithms}
  \crefformat{definition}{Def.~#2#1#3}
  \Crefformat{definition}{Definition~#2#1#3}
  \crefformat{equation}{Eq.~(#2#1#3)}
  \Crefformat{equation}{Equation~(#2#1#3)}
  \crefformat{example}{Expl.~#2#1#3}
  \Crefformat{example}{Example~#2#1#3}
  \crefformat{figure}{Fig.~#2#1#3}           \crefname{figure}{Fig.}{Figs.}
  \Crefformat{figure}{Figure~#2#1#3}         \Crefname{figure}{Figure}{Figures}
  \crefformat{lemma}{Lem.~#2#1#3}
  \Crefformat{lemma}{Lemma~#2#1#3}
  \crefformat{page}{p.~#2#1#3}
  \Crefformat{page}{Page~#2#1#3}
  \crefformat{problem}{Prob.~#2#1#3}
  \Crefformat{problem}{Problem~#2#1#3}
  \crefformat{proposition}{Prop.~#2#1#3}
  \Crefformat{proposition}{Proposition~#2#1#3}
  \crefformat{remark}{Rem.~#2#1#3}
  \Crefformat{remark}{Remark~#2#1#3}
  \crefformat{section}{Sec.~#2#1#3}         \crefname{section}{Sec.}{Secs.}
  \Crefformat{section}{Section~#2#1#3}      \Crefname{section}{Section}{Sections}
  \crefformat{table}{Tab.~#2#1#3}           \crefname{table}{Tab.}{Tables}
  \Crefformat{table}{Table~#2#1#3}          \Crefname{table}{Table}{Tables}
  \crefformat{theorem}{Thm.~#2#1#3}         \crefname{theorem}{Thm.}{Thms.}
  \Crefformat{theorem}{Theorem~#2#1#3}      \Crefname{theorem}{Theorem}{Theorems}

\newcommand*{\E}[1]{\ensuremath{\mathrm{E}{#1}}}

\providecommand*{\gray}{\leavevmode\color{gray}}
\newcommand*{\setE}{\ensuremath{\mathcal{E}}}
\usepackage{upgreek}                                             
\newcommand*{\muCT}{\mbox{$\upmu$CT}}                            
\usepackage{dsfont}                                              
  \newcommand{\IN}{\ensuremath\mathds{N}}
  \newcommand{\IP}{\ensuremath\mathds{P}}
  \newcommand{\IR}{\ensuremath\mathds{R}}

  \renewcommand*{\vec}[1]{{\boldsymbol{#1}}}                     
\DeclareMathAlphabet{\mathbfsf}{\encodingdefault}{\sfdefault}{bx}{n}
  \newcommand*{\vecc}[1]{\mathbfsf{#1}}                          
\newcommand*{\laplace}{\upDelta}                                 
\newcommand*{\grad}{\vec{\nabla}}                                
\renewcommand*{\div}{\vec{\nabla}\cdot}                          
\newcommand*{\llbrace}{\lbrace\hspace*{-0.18em}\vert}
\newcommand*{\rrbrace}{\vert\hspace*{-0.18em}\rbrace}
\newcommand*{\Pe}{\mathrm{Pe}}                                   
\newcommand*{\avg}[1]{\llbrace{#1}\rrbrace}                      
\newcommand*{\avgHarmonic}[1]{\langle{#1}\rangle}                
\newcommand*{\jump}[1]{\left\llbracket{#1}\right\rrbracket}      
\newcommand*{\dd}{\mathrm{d}}                                    
\newcommand*{\abs}[1]{\ensuremath{|#1|}}                         
\newcommand*{\Nel}{N_\mathrm{el}}                                
\newcommand*{\Nst}{N_\mathrm{st}}                                
\newcommand*{\Nloc}{N_\mathrm{loc}}                              
\DeclareMathOperator*{\sig}{sig}                                 
\newcommand*{\brokennorm}[2]{|||#1|||_{#2}}                      
\newcommand*{\norm}[2]{\|#1\|_{#2}}                              
\newcommand*{\normal}{\vec{n}}                                   
\newcommand*{\transpose}[1]{{#1}^\mathrm{T}}                     
\DeclareMathOperator*{\round}{round}                             

\journal{Computer Methods in Applied Mechanics and Engineering}
\begin{document}
\begin{frontmatter}

\title{A~finite volume\,/\,discontinuous Galerkin method for the\\advective Cahn--Hilliard~equation with degenerate mobility\\on porous domains stemming from micro-CT~imaging}
\author[label1]{Florian Frank}
\author[label1]{Chen Liu}
\author[label2]{Faruk O.~Alpak}
\author[label1]{Beatrice Riviere}
\address[label1]{Rice University, CAAM Department, 6100~Main~Street, Houston, TX 77005, USA}
\address[label2]{Shell International Exploration and Production Inc., Shell Technology Center, 3333 Highway 6 South, Houston, TX 77082, USA}


\begin{abstract}
A~numerical method is formulated for the solution of the advective Cahn--Hilliard~(CH) equation with constant and degenerate mobility in three-dimensional porous media with non-vanishing velocity on the exterior boundary.  
The CH~equation describes phase separation of an~immiscible binary mixture at constant temperature in the presence of a~mass constraint and dissipation of free energy.
Porous media\,/\,pore-scale problems specifically entail high-resolution images of rocks in which the solid matrix and pore spaces are fully resolved.
The interior penalty discontinuous Galerkin method is used for the spatial discretization of the CH~equation in mixed form, while a~semi-implicit convex--concave splitting is utilized for temporal discretization. 
The spatial approximation order is arbitrary, while it reduces to a~finite volume scheme for the choice of elementwise constants.
The resulting nonlinear systems of equations are reduced using the Schur complement and solved via Newton's method.  
The numerical scheme is first validated using numerical convergence tests and then applied to a~number of fundamental problems for validation and numerical experimentation purposes including the case of degenerate mobility.
First-order physical applicability and robustness of the numerical method are shown in a~breakthrough scenario on a~voxel set obtained from a~micro-CT scan of a~real sandstone rock sample.
\end{abstract}

\begin{keyword}
Cahn--Hilliard equation  \sep phase field \sep advection \sep fourth order parabolic equation \sep discontinuous Galerkin \sep finite volumes \sep micro-CT scan \sep voxel sets \sep porous media

\end{keyword}
\end{frontmatter}

\section{Introduction}
Phase-field models are becoming increasingly popular for modeling multi-phase fluid flow below Darcy scale. The fundamental idea underlying phase-field models is to introduce a~possibly conserved order parameter, e.\,g., a~mass fraction that varies continuously over thin interfacial layers. The order parameter is uniform in pure (\enquote{bulk}) phases, and instead of sharp phase interfaces there are diffuse finite-thickness transition regions with smoothly varying interfacial forces. The development of diffuse-interface models in hydrodynamics and their application to a~wide variety of interfacial phenomena is presented in a~comprehensive review paper by~\cite{AndersonFaddenWheeler1998}. A~more recent survey of phase-field models~\cite{Kim2012review} addresses improvements as modeling ternary fluids~\cite{MorralCahn1971,BoyerLapuerta2006,Kim2007,BoyerEtAl2009CHNS,KimLowengrub2005}, drop coalescence and retraction in viscoelastic fluids~\cite{YueEtAl2005}, and contact angle boundary conditions~\cite{DingSpelt2007,He2011Contact,Jacqmin2000,KhatavkarEtAl2007,LeeKim2011,YueFeng2011,ARF2016}.
\par
The Cahn--Hilliard~(CH) equation~\cite{CahnHilliard1958,Cahn1961} is a~stiff, fourth-order, nonlinear parabolic partial differential equation, which is conservative for its order parameter.
This equation has the same equilibria as the Allen--Cahn equation, which is of second order.  Both equations can be obtained as limits of the so-called viscous CH~equation~\cite{BaiEtAl1995CH}.
The CH~equation was initially proposed as a~phenomenological model for phase separation that occurs in binary alloys. Since then it has undergone considerable development both in terms of physics of phase separation phenomena as well as of computational methods in order to make applications feasible. Phase separation phenomena can be observed in a~variety of materials in addition to alloys including two-phase fluid systems~\cite{Emmerich2011Book}, polymers~\cite{SaxenaCaneba2002,Aristotelous2013}, ceramics~\cite{Cogswell2010Thesis}, and tumor tissues~\cite{WiseEtAl2008Tumor,WuZwietenZee2014Tumor,Colli2015Tumor}. The CH~equation has also an~application in image processing for enhancing the sharpness of images or for image in-painting~\cite{BertozziEtAl2007}.  The standard variant of the CH~equation characterizes phase segregation of a~binary mixture, which is the alignment of a~system into spatial domains predominated by one of the two components, in the presence of a~mass constraint and dissipation of free energy.  This comprises a~rapid phase separation process, which gives rise to interface creation, and a~slow aggregation process, which in turn leads to the development of bulk phases with diffuse interfaces separating them.  Different temporal and spatial scales characterize these two phenomena rendering the accurate and efficient solution of the CH~equation a~computational challenge.  
\par
The main unknown of the CH~equation is the~\emph{order parameter}~$c$, which can either be a~volume or mass fraction of one of the two components~$c_1, c_2$, i.\,e.~$c = c_1$, $c_2 = 1 - c_1$, $c\in[0,1]$ 
or the difference between fractions, i.\,e.~$c= c_1 - c_2 = 2c_1 - 1$, $c\in[-1,1]$.  We use the latter choice throughout this paper.
Cahn and Hilliard~\cite{CahnHilliard1958} proposed the following approximation of \emph{total Helmholtz free energy}~$F(c)$ for an incompressible binary isotropic solution in a~domain~$\Omega\subset\IR^3$ :
\begin{equation}\label{eq:HelmholtzFreeEnergy}
F(c) ~\coloneqq~ \int_\Omega \Big(\Phi(c)        + \frac{\kappa}{2}\, \abs{\grad c}^2\Big)\,\dd\vec{x}~,
\end{equation}
where $\Phi$ is called \emph{chemical energy density} or \emph{homogeneous free energy density} corresponding to the Helmholtz free energy of a~unit volume of homogeneous material
of composition~$c$.  The function~$\Phi$ has the shape of a~double-well potential with minima at the bulk compositions, i.\,e.~$c\in\{-1,1\}$, provided that the (constant) system temperature is below the critical temperature of the mixture (i.\,e.~where a~homogeneous mixture becomes energetically unfavorable).  A~frequently used expression for~$\Phi$ is 
\begin{equation}\label{eq:PhiOfc}
\Phi(c) ~\coloneqq~\frac{1}{4} (c^2 - 1)^2~.
\end{equation}
The second term in~\eqref{eq:HelmholtzFreeEnergy} is sometimes called \enquote{\emph{gradient energy density}} or \enquote{\emph{interfacial energy density}}, where the integrand~$\frac{\kappa}{2} \abs{\grad c}^2$ with $0<\kappa \ll 1$ is the first term of an~expansion
representing the increase in free energy due to introducing a gradient of composition.
The \emph{chemical potential} is the variational/functional derivative of the total Helmholtz free energy with respect to~$c$:
\begin{equation}\label{eq:definitionChemicalPotential}
\mu ~\coloneqq~ \delta_c F(c) ~=~ \Phi'(c) - \kappa\, \laplace c~.
\end{equation}  
A~chemical potential gradient implies a~mass flux~$\vec{j}^\mathrm{diff}\coloneqq -M(c)\grad\mu$, where the phenomenological coefficient~$M(c)\geq 0$ is called \emph{mobility}.  A~consequence of mass conservation is the CH~equation
(without advection):
\begin{equation*}
\partial_t c + \div \vec{j}^\mathrm{diff} ~=~ 0~.
\end{equation*}
\par
The wide spectrum of the CH~equation applications makes the development of its efficient numerical solution highly desirable. However, the nonlinear character and the fourth order of the equation render the task demanding, especially for 3D~problems. There have been significant developments in numerical methods both in terms of spatial and temporal approximation techniques for solving the CH~equation.
In general, the spatial discretization is either based on the mixed CH~formulation and solves for both the order parameter~$c$ and the chemical potential~$\mu$, or on a~primal (also primitive-variable) formulation by substitution~$\mu$ by its expression in~$c$.  The first case yields a~coupled system of two second order elliptic equations;  one of which contains no time-derivative.  This coupled formulation has the advantage of allowing lower order of approximation, however it leads to a~saddle-point problem. 
\par
Various numerical schemes have been developed based on finite difference~\cite{SaylorEtAl2007,Furihata2001,Kim2007}, finite volume~(FV)~\cite{CuetoFelguerosoPeraire2008,KimKang2009}, finite element~\cite{ElliottFrench1989,ElliottGarcke1996,ZhangWang2010,DuNicolaides1991,WodoGanapathysubramanian2011,CopettiElliott1992,BarrettBloweyGarcke1998},  mixed finite element~\cite{ElliottFrankMilner1989,FengProhl2004}, and spectral methods~\cite{He2007} with specific advantages and drawbacks. 
It has to be noted that the majority of the reported applications were developed for two dimensional problems. Developments for simulations of the evolution of 3D~morphology of separating phases are more recent~\cite{Gomez2008Iso,WodoGanapathysubramanian2011}. Another recent development is the use of adaptive time stepping~\cite{WuZwietenEtAl2013Adaptive,TierraGonzales2015}.
\par
The discontinuous Galerkin~(DG) method has been applied to the primal (i.\,e.~forth order) CH~equation in~\cite{WellsKG2006,FengKarakashian2007,AristotelousKarakashianWise2014} and to the mixed version in~\cite{Aristotelous2013,KaySS2009}; the work, \cite{KaySS2009} also includes the advection operator, however with vanishing velocity across the domain boundary.
By introducing auxiliary flux variables, the mixed CH~equation can also be formulated as a~system of coupled first order equations, where some of the equations are vector-valued.  DG~methods based on this ansatz are known as~local discontinuous Galerkin~(LDG) schemes~\cite{CockburnShu1998,AizingerDCC2000}.  LDG~methods for the CH~equation are found in~\cite{XiaEtAl2007LDG,XiaEtAl2009,GuoXu2013LDG}.  Except of~\cite{WellsKG2006,XiaEtAl2007LDG,XiaEtAl2009,GuoXu2013LDG}, all cited DG~papers use constant mobility.  To the best knowledge of the authors, there are no publications that show DG~simulations in three space dimensions except of~\cite{GuoXu2013LDG}, which studies a multigrid solver for LDG and third-order implicit Runge--Kutta method. 

\par
In this paper, we formulate and implement a~DG~method for the effective numerical solution of the mixed CH~equation in three space dimensions. 
The model equations take a~degenerate mobility coefficient into account and include an~advection term and inhomogeneous inflow boundary conditions, which become important when simulating flow and transport through a~porous domain (note that there are models that also respect the impact of the interface on the flow field~\cite{LowengrubTruskinovsky1998}).
The first order time discretization is implicit--explicit based on a~convex--concave decomposition of the nonlinear chemical energy density~$\Phi$. 
We follow the approach of discretizing the CH~equation in mixed form, in particular, since DG~schemes based on the primal formulation require an~additional penalty term, penalizing the normal fluxes across element boundaries.  The scaling of this term has a~significant unfavorable impact on the condition number of the system matrix, and also on the magnitude of the discretization error~\cite[cf.][]{WellsKG2006}, and according to our own experience adjusting the scaling of this term is not trivial. A~Schur-complement reduction halves the degrees of freedom by resubstituting\,/\,eliminating the chemical potential on the fully discrete level.
Thereby, we avoid solving saddle point problems, which are unsuitable for standard iterative linear solvers. 
Our implementation aims at low order schemes and owing to the structure of the grid (voxel set), it is even possible to reduce the polynomial order to zero, yielding a~cell-centered~FV approach. Numerical validation tests are conducted both for constant and nonlinear mobility cases. These tests confirm the theoretical orders of convergence for all investigated polynomial orders of approximation. Intuitive numerical experiments are performed to investigate discrete energy dissipation and to compare the effect of using constant vs.~nonlinear mobility.  A~flow simulation through porous medium stemming from microtomography demonstrates the applicability of the numerical scheme toward realistic scenarios.
\par
The paper is organized in the following manner: We introduce the~CH~model problem in~\cref{sec:ModelProblem}.
The triangulation is discussed and our numerical scheme defined in~\cref{sec:TriangulationAndScheme}.
\Cref{sec:SimulationResults} is dedicated to the numerical validation tests and numerical simulations of physical experiments
including various phase segregation phenomena including spinodal decomposition, droplet in a~fluid, and pore-scale flow. Some concluding remarks follow in~\cref{sec:Conclusion}.

\section{Model problem}\label{sec:ModelProblem}
Let~$J\coloneqq (0,T)$ denote the considered time interval with end time $T\in\IR^+$.
For an~unknown $c:\overline{J}\times\overline{\Omega}\rightarrow [-1,1]$, we consider the \emph{advective Cahn--Hilliard equation} with possibly degenerate mobility~$M$ with $M(c)\geq 0$:
\begin{subequations}\label{eq:model}
\begin{align}
\partial_t c  - \frac{1}{\Pe} \,\div\big(M(c)\,\grad \mu \big)   + \div(\vec{v}\, c)   &~=~ 0   && \text{in}~J\times\Omega~, \label{eq:model:a}\\
c     &~=~ c^\mathrm{in}                                                                        && \text{on}~J\times\partial\Omega^\mathrm{in}(t)~,\label{eq:model:in}\\
\grad c\cdot\normal               &~=~ 0                                                        && \text{on}~J\times\partial\Omega^\mathrm{out}(t)~,\label{eq:model:out}\\
 M(c)\,\grad\mu  \cdot\normal     &~=~ 0                                                        && \text{on}~J\times\partial\Omega~,\label{eq:model:normalFluxMu} \\
c                                 &~=~ c^0                                                      && \text{on}~\{0\}\times\Omega \label{eq:model:initial}
\end{align}
\end{subequations}
with (time-dependent) solenoidal mixture velocity~$\vec{v}:J\times\Omega\rightarrow\IR^3$, P\'eclet number~$\Pe>0$, inflow boundary~$\partial\Omega^\mathrm{in}(t) \coloneqq \{ \vec{x} \in \partial\Omega \,;\, \vec{v}(t, \vec{x}) \cdot \normal < 0 \}$, $\normal$~denoting the unit normal on~$\partial \Omega$ exterior to~$\Omega$, outflow boundary~$\partial\Omega^\mathrm{out}(t) \coloneqq \partial\Omega \setminus \partial\Omega^\mathrm{in}(t)$, $c^0:\Omega\rightarrow [-1,1]$ and~$c^\mathrm{in}:J\times \partial\Omega^\mathrm{in}(t)\rightarrow[-1,1]$ the given initial and inflow boundary data, respectively. 
The partition of~$\partial\Omega$ into inflow boundary and outflow boundary may change in time due to the time-depending velocity field~$\vec{v}$.
The boundary conditions~\mbox{\eqref{eq:model:in}--\eqref{eq:model:normalFluxMu}} are found, e.\,g., in~\cite{BoterMinjeaud2011}.  
\par
The chemical potential~$\mu: J\times\Omega\rightarrow\IR$ was defined in~\eqref{eq:definitionChemicalPotential}.  
Typical choices for~$M$ are constant mobility~$M(c)\equiv M_0>0$, where~$M_0$ can be set to one after rescaling in time, and~$M(c) = 1-c^2$, one version of the so-called \enquote{degenerate mobility} \cite{jingxue1992existence,ElliottGarcke1996}.
For convenience, we set the expression
\begin{equation}\label{eq:mobility}
M(c) ~\coloneqq~ 1 - \beta\,c^2
\end{equation}
in this paper with $\beta\in\{0,1\}$, such that we obtain the constant mobility case for~$\beta = 0$ and the degenerate case for~$\beta = 1$.
The parameter~$0<\kappa\ll 1$ in \eqref{eq:HelmholtzFreeEnergy} is typically chosen as small as possible to generate a~small interface width~$\delta$.
Equilibrium profiles of~$c$ are obtained by solving~$\mu(c) = 0$ for~$c$, which yields an interface width of approximately~$\delta\approx 4\,\sqrt{\kappa}$ in one space dimension for~$\Phi$ as given in~\eqref{eq:PhiOfc} (see~\cite{WodoGanapathysubramanian2011} and references therein for possible definitions of~$\delta$).  Simulation studies have shown that this estimate also holds in three dimensions, and that 
 in numerical simulations~$\kappa$ has to be limited from below by the mesh size~$h$.  We choose $\kappa$ such that $\kappa \geq h^2$, which means there are at least four elements across the interface.
\par
The CH~equation with constant mobility is mathematically studied, e.\,g., in~\cite{BaiEtAl1995CH,Elliott1989CH,ElliottSongmu1986}.  Existence and regularity results for the degenerate mobility case with vanishing velocity are presented in~\cite{jingxue1992existence,ElliottGarcke1996}.

\subsection{Continuous solution properties}\label{sec:SolutionProperties}
Solutions of~\eqref{eq:model} enjoy the following properties:
\begin{enumerate}
\item 
\emph{Mass conservation}:
The change in time of the total amount is
\begin{equation}\label{eq:massConservation}
\dd_t \int_\Omega c 
~=~  \int_\Omega  \partial_t c
\overset{\eqref{eq:model:a}}{~=~}\frac{1}{\Pe}\int_\Omega \div\Big(M(c)\,\grad\mu \Big) - \int_\Omega\div(\vec{v}\,c)
\overset{\eqref{eq:model:normalFluxMu}}{~=~}  - \int_{\partial\Omega} \vec{v}\cdot\normal\,c ~.
\end{equation}
In particular, if~$\int_{\partial\Omega}\vec{v}\cdot\normal\,c = 0$, which is equivalent to $\int_{\partial\Omega^\mathrm{in}(t)}\vec{v}\cdot\normal\,c^\mathrm{in} = -\int_{\partial\Omega^\mathrm{out}(t)}\vec{v}\cdot\normal\,c$ (the mass entering the domain is equal to the mass transported out of the domain) then \eqref{eq:massConservation} implies that the total amount of the order parameter~$c$ is preserved, i.\,e.
\begin{equation*}
\forall t\in J~, \quad \frac{1}{\abs{\Omega}}  \int_\Omega c(t,\vec{x})\,\dd\vec{x} ~=~\bar{c}~,
\end{equation*}
where $\bar{c} \coloneqq  \frac{1}{\abs{\Omega}}  \int_\Omega c^0$.
\item 
\emph{Energy dissipation}:
For a~closed system (i.\,e.~$\vec{v}\cdot\vec{n}$ on~$\partial\Omega$), the total Helmholtz free energy is non-increasing in time, i.\,e.,
\begin{equation}\label{eq:discreteEnergyDissipation}
\forall t\in J~, \quad \dd_t F(c)(t) ~\leq~0~.
\end{equation}
Using the chain rule for functionals and the homogeneous boundary condition for the normal gradient of the chemical potential, we obtain (pointwise in~$t$), 
\begin{align*}
\dd_t F(c) &~=~ \int_\Omega \delta_c F(c)\,\partial_t c
\overset{\eqref{eq:definitionChemicalPotential},\eqref{eq:model:a}}{~=~}\int_\Omega \mu\,\div\Big( M(c)\,\grad\mu\Big) ~=~ - \int_\Omega \grad\mu\cdot \Big( M(c)\,\grad\mu\Big) + \int_{\partial\Omega} \mu\,M(c)\,\grad\mu\cdot\normal
\\
&~=~ - \norm{M(c)^{1/2}\grad\mu}{L^2(\Omega)}^2
~\leq~ 0~.
\end{align*}
\item
\emph{Maximum principle} \cite{ElliottGarcke1996}:
For the case of degenerate mobility~$\beta = 1$, under additional smoothness assumptions on~$\partial\Omega$ and the initial data~$c^0(\vec{x})\in[-1,1]$ for~a.\,e.~$\vec{x}\in\Omega$, the order parameter~$c$ is in the physically meaningful interval of~$[-1,1]$ for a.\,e.~$(t,\vec{x})$ in $J\times\Omega$.
\end{enumerate}
\par
A~maximum principle does also hold when using a~logarithmic type of energy density~$\Phi$ instead of~\eqref{eq:PhiOfc} and either choice of mobility, \cite[cf.][(1.3)]{TierraGonzales2015}.
The computational domain~$\Omega_h$ introduced in~\cref{sec:TriangulationAndScheme} clearly does not meet the smoothness requirements for item~3.  However, it is desirable that numerical schemes 
at least approximate similar discrete versions of the properties above.  In which extend this is the case for our DG~scheme will be discussed in~\cref{sec:spaceDiscretization}.

\section{Triangulation and numerical scheme}\label{sec:TriangulationAndScheme}

\subsection{Domain and triangulation}

\begin{figure}[ht!]
\begin{tabularx}{\linewidth}{@{}LCR@{}}
\includegraphics[width=\linewidth]{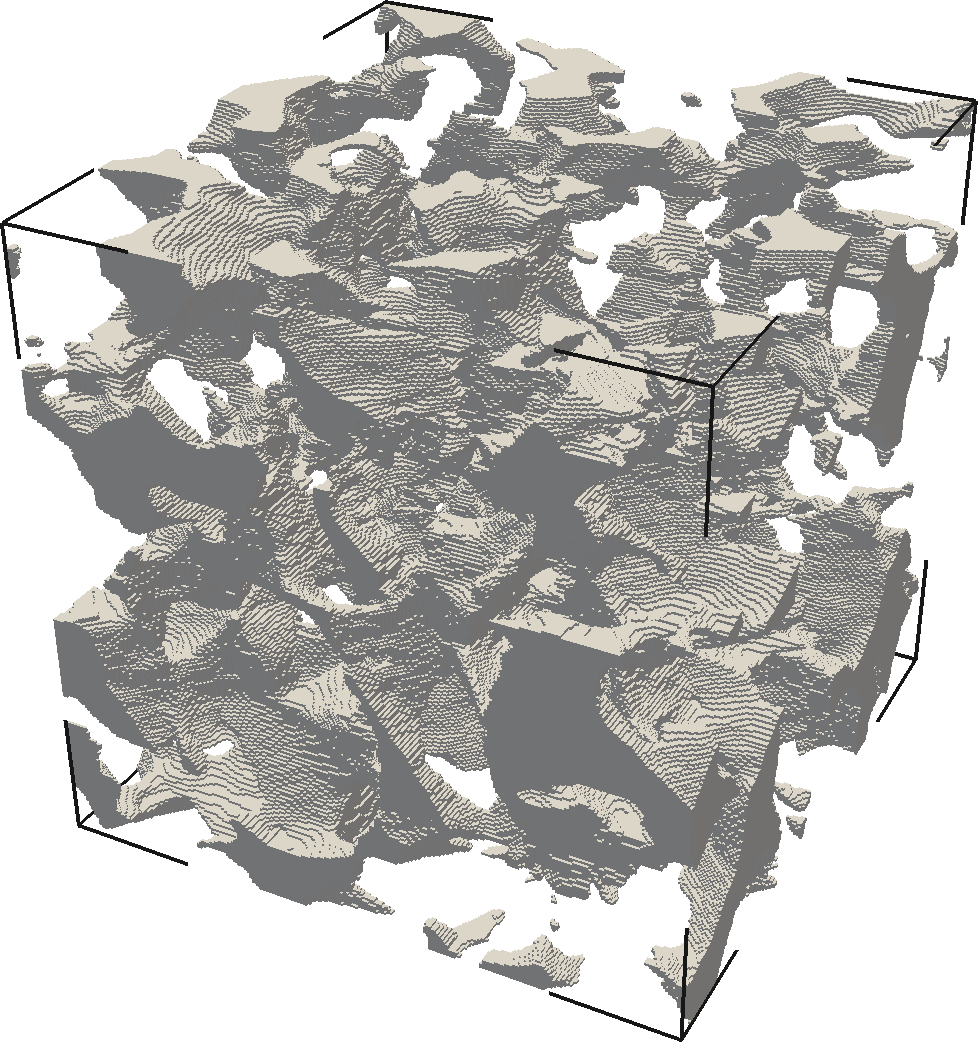} &
\includegraphics[width=\linewidth]{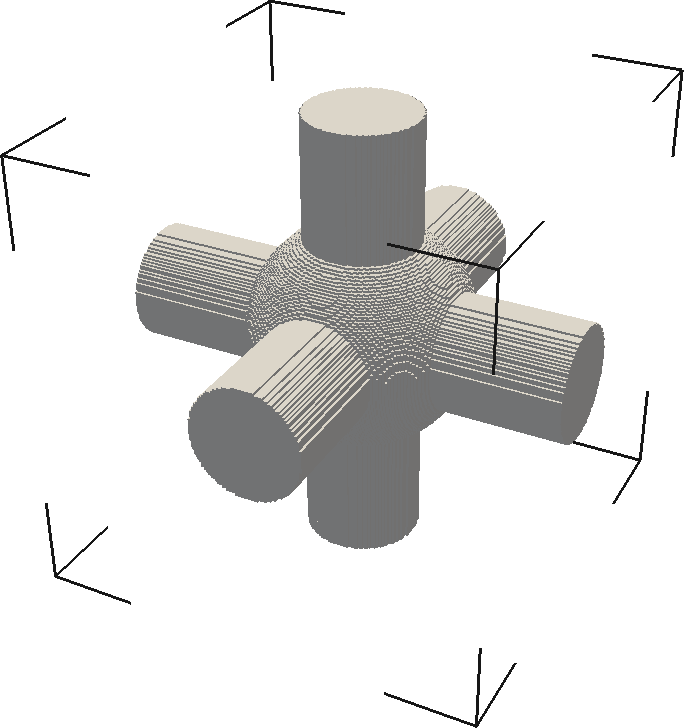} &
\includegraphics[width=\linewidth]{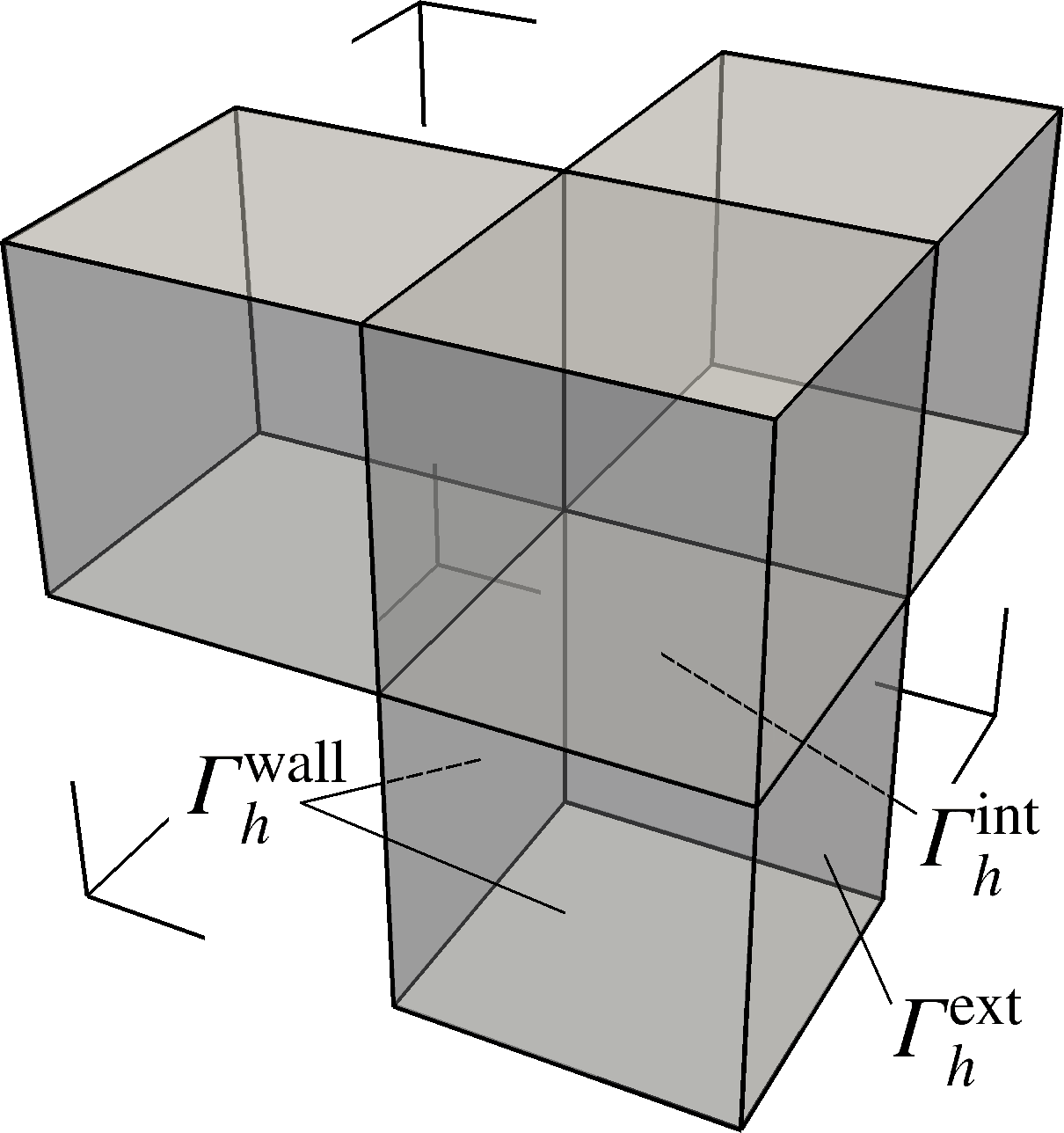} 
\end{tabularx}
\caption{Examples of computational domains~$\Omega_h\subset[0,1]^3$ illustrated in gray: \mbox{$256^3$-voxel} \muCT~image~\emph{(left)}, \mbox{$200^3$-voxel} pore network~\emph{(middle)}, \mbox{$2^3$-voxel} example containing 4~elements illustrating the face sets~$\Gamma_h^\mathrm{int}$, $\Gamma_h^\mathrm{ext}$, $\Gamma_h^\mathrm{wall}$~\emph{(right)}.  In the latter case, faces on~$\partial[0,L]^3\cap\partial\Omega_h$ can either be declared to be part of~$\Gamma_h^\mathrm{wall}$ or $\Gamma_h^\mathrm{ext}$.}
\label{fig:compDomainsAndIndex}
\end{figure}

\paragraph{Micro-CT images}
Modern microtomography (micro-CT, \muCT) can create cross-sections of small rock samples on $\mathrm{nm}$ to~$\upmu\mathrm{m}$~scale by means of x-rays, which subsequently can be used to construct a~model by 3D~imaging software.  In our case, this model consists of a~set of binary voxels, each of which is either associated with the pore space of the rock sample or with its solid matrix. We consider the matrix itself rigid and impermeable for fluids.  Our goal is to perform simulations directly on the given data, more precisely on the cubic-shaped voxels that represent the \emph{pore space}, which we assume saturated by wetting fluids, cf.~\cref{fig:compDomainsAndIndex}. As~every voxel is a~regular cube of size~$h$ we can heavily exploit the properties of a~Cartesian grid. The advantage here is that except of the voxel set itself, there is no need to store grid-topological data as, e.\,g., the indices of adjacent elements or geometric properties of every element.  Instead, these data can be computed \enquote{on the fly} and thus no additional memory is consumed here.  The voxel set is efficiently stored as a~three-dimensional boolean array.

\paragraph{Computational domain, indexing, and boundary partitions}
The unit cube~$[0,1]^3$ is divided into~$N$ voxels in each direction, i.\,e.~the mesh size is~$h\coloneqq 1/N$.  We assign two types of indices to voxels, between which one-to-one mappings exist: \emph{tuple indices} and \emph{flat indices}.  The first one bears the voxel coordinate from which adjacencies can be determined immediately, while the latter one realizes a~consecutive scalar index that is required by the numerical scheme.
\par
The voxel containing the origin is associated with the tuple index~$(0,0,0)$ and the one in the opposite corner with~$(N-1,N-1,N-1)$.   
For a~tuple index~$(i,j,k)\in\{0,\ldots,N-1\}^3$, we define the flat index~$n\in\{0,\ldots,N^3-1\}$ by $n \coloneqq i + N\,j + N^2\,k$ and vice versa, $i= n\%N$, $j=  (n/N)\%N$, $k = n/N^2$, where the operators~\enquote*{$/$} and~\enquote*{$\%$} yield the integer quotient and the remainder of an~integer division, respectively.
\par
We call the (topologically closed) voxels that are associated with the pore space \enquote{elements}~$E_k$, $k\in\{0,\ldots,\allowbreak\mbox{$\Nel-1$}\}$ with edge length~$h$, i.\,e.,~$\setE_h\coloneqq\{E_k\}_k$ is a non-overlapping partition~$\overline{\Omega}_h\coloneqq\cup_k E_k$.  The computational domain~$\Omega_h$ is not necessarily connected, i.\,e.,~it may consist of several, non-connected subdomains.   The union of all boundary parts is denoted by~$\partial\Omega_h$.  Faces that are subset of~$\partial\Omega_h$ divide into two sets: $\Gamma_h^\mathrm{wall}$ containing faces that lie between an~element and a~solid voxel and faces on~$\partial[0,L]^3\cap\partial\Omega_h$ that are associated with an~impermeable (solid wall) boundary, and  $\Gamma_h^\mathrm{ext}$ containing faces on~$\partial[0,L]^3\cap\partial\Omega_h$ that are associated with an~exterior fluid domain, cf.~\cref{fig:compDomainsAndIndex}.
Let~$\Gamma_h^\mathrm{int}$ denote the set of \emph{interior faces}, i.\,e.,~faces connecting two elements. 
The resulting triangulation~$\setE_h$ has the structure of a~pseudo-Cartesian grid, which is exploited in the implementation.

\subsection{Time discretization}
An~explicit time stepping of system~\eqref{eq:definitionChemicalPotential}--\eqref{eq:model} is only possible when based on the primal (non-mixed) formulation since~\eqref{eq:definitionChemicalPotential} does not contain an~evolution term; however this entails a~severe time step size restriction of order~$O(h^4)$.  For fully implicit time steps on the other hand, instabilities have been observed~\cite{DuNicolaides1991,QiaoSun2014}, possibly due to a~loss of uniqueness for not sufficiently small time steps.  We adopt a~popular method that is unconditionally stable by using a~convex--concave splitting of~$\Phi$ of the form
\begin{equation} \label{eq:convexConcaveSplitting}
\Phi(c)~=~\Phi_{+}(c) + \Phi_{-}(c)~,
\end{equation}
where the convex part~$\Phi_{+}$ is treated implicit and the concave part~$\Phi_{-}$ explicit in time~\cite{ElliottStuart1993,Eyre1997ConvexConcave}.  Certainly, this approach yields a~nonlinear system in each time step, if~$\Phi_{+}$ is nonlinear.
\par
Let $0\eqqcolon<t_1<\ldots< t_{\Nst}\coloneqq T$
be a~decomposition of~$J$ into~$\Nst$ subintervals, not necessary equidistant. 
The \mbox{$n$th}~time step size is denoted by~$\tau_n \coloneqq t_n - t_{n-1}$.
Using~\eqref{eq:definitionChemicalPotential} to write~\eqref{eq:model} in mixed form the time-discrete problem reads:
\par
For $n\in\{1,\ldots,\Nst\}$, find $\mu^n$, $c^n$ such that
\begin{subequations}\label{eq:model:timeDiscrete}
\begin{align}
 \label{eq:model:timeDiscrete:a}
 - \kappa\, \laplace c^n - \mu^n +\Phi_+'(c^n)                                                &~=~ -  \Phi_-'(c^{n-1})   && \text{in}~\Omega_h~,
 \\
 \label{eq:model:timeDiscrete:b}
c^n  - \frac{\tau_n}{\Pe} \,\div\big(M(c^{n-1})\,\grad \mu^n   \big) + \tau_n \, \div\big(\vec{v}(t_n)\, c^n\big) &~=~  c^{n-1}              && \text{in}~\Omega_h~,
\end{align}
\end{subequations}
subjected to respective boundary and initial conditions, cf.~\eqref{eq:model:in}--\eqref{eq:model:initial}.
For a closed system (i.\,e.~$\vec{v}\cdot\vec{n}$ on~$\partial\Omega_h$), the resulting scheme is uniquely solvable and unconditionally energy stable, i.\,e., $\forall n\in\{1,\ldots,\Nst\},~ F\big(c^{n}\big) \leq F\big(c^{n-1}\big)$.
\par
For a~special form of chemical energy density, namely~\eqref{eq:PhiOfc},
the corresponding nonlinear term~$\Phi'(c)$ in~\eqref{eq:model:a} can be treated with \enquote{Eyre's scheme}~\cite{Eyre1998EyreScheme,deMello2005}, which leads to a~\emph{linear} (however conditionally stable~\cite{TierraGonzales2015}) approximation of the convex--concave splitting~\eqref{eq:convexConcaveSplitting}.  
Treating the mobility explicit in time is a~common practice \cite{GrunKlingbeil2014,BarrettBloweyGarcke1999,BarrettBloweyGarcke1998} to overcome delicate complexities in the nonlinear solver.  In particular, our numerical simulations have shown that the overshoots beyond the meaningful interval of~$[-1,1]$ are larger when the mobility is taken implicitly into account. 
The advection term in~\eqref{eq:model:timeDiscrete:b} was taken into account time-implicitly, as the expected overshoot will be less than in the explicit case for high P\'eclet numbers.
A~review of first and second order of approaches and their comparison regarding time accuracy error, energy stability, unique solvability is found in~\cite{TierraGonzales2015}.

\subsection{Approximation spaces}\label{sec:ApproximationSpaces}
We use a~mixture of algebraic and numerical indexing: for instance, 
$e_{km}\in\partial E_k\cap\Gamma_h^\mathrm{int}$ is associated with \emph{all possible} element indices~$k\in\{0,\allowbreak\ldots,\allowbreak\Nel-1\}$ 
and local face indices $m\in\{0,\ldots,5\}$ such that~$e_{km}$ lies in~$\partial E_k\cap\Gamma_h^\mathrm{int}$. 
This implicitly fixes the numerical indices which accordingly can be used in the numerical schemes to index matrices or arrays.
\par
Let~$\IP_p(\hat{E})$ denote the space of polynomials of degree at most~$p$ on the \emph{reference element}~$\hat{E}\coloneqq[-1,1]^3$ and let $\Nloc\coloneqq\dim \IP_p(\hat{E}) = (p+1)(p+2)(p+3)/6$ denote the \emph{number of local degrees of freedom}, i.\,e.,~$\Nloc$~linearly independent \emph{basis functions}~$\hat{\varphi}_i:\hat{E}\rightarrow \IR$ span the space~$\IP_p(\hat{E})$. 
Let~$\vec{F}_k:\hat{E}\rightarrow E_k$ denote the mapping from~$\hat{E}$ to the elements~$E_k\in\setE_h$, which is in our setting the composition of a~scaling and a~shift.
The broken polynomial space on the triangulation~$\setE_h$ is given as
\begin{equation*}
\IP_p(\setE_h)~\coloneqq~\prod\nolimits_{E_k\in\setE_h}\IP_p(E_k)~,
\end{equation*}
with $\IP_p(E_k) \coloneqq \mathrm{span}\{\varphi_{ki}\}_{i\in\{0,\ldots,\Nloc-1\}}$,
where the basis functions are~$\varphi_{ki}: E_k\ni\vec{x}\mapsto \hat{\varphi}_i\circ\vec{F}_k^{-1}(\vec{x})\in\IR$, i.\,e., $\varphi_{ki}\coloneqq\hat{\varphi}_i\circ\vec{F}_k^{-1}$. Here we assumed tacitly that the polynomial degree\,/\,local degrees of freedom~$\Nloc$ is the same on each element, which means that we have~$\Nel \Nloc$ global degrees of freedom.
Note that elements of~$\IP_p(\setE_h)$ are in general \emph{discontinuous} on interior faces~$e\in\Gamma_h^\mathrm{int}$.
\par
We use a~hierarchical modal orthonormal basis on the reference element~$\hat{\varphi}_0,\ldots,\hat{\varphi}_{\Nloc-1}$ spanning~$\IP_p(\hat{E})$, constructed using tensor products of one~dimensional Legendre polynomials, cf.~\cite{EG2004}. 
The basis functions~$\hat{\varphi}_i$ are \emph{orthonormal} with respect to the~$L^2(\hat{E})$ inner product, i.\,e., $\int_{\hat{E}}\hat{\varphi}_i\,\hat{\varphi}_j = \delta_{ij}$ holds, where~$\delta_{ij}$ denotes the Kronecker symbol.  Thus, the entries of the mass matrix~$\vecc{M}$ become
\begin{equation}\label{eq:MassMatrix}
[\vecc{M}]_{k\Nloc + i, l \Nloc + j}  ~\coloneqq~ \int_{\Omega_h}\varphi_{lj} \, \varphi_{ki} ~=~ \delta_{lk} \int_{E_k} \varphi_{lj} \, \varphi_{ki} ~=~ \frac{h^3}{8}\,\delta_{lk} \int_{\hat{E}} \hat{\varphi}_i\,\hat{\varphi}_j ~=~ \frac{h^3}{8}\,\delta_{ij}\,\delta_{lk}
\end{equation}
for $k,l\in\{0,\ldots,\Nel-1\}$ and $i,j\in\{0,\ldots,\Nloc-1\}$.
The \emph{hierarchy} property, i.\,e., $\IP_0(\hat{E}) = \mathrm{span}\{\hat{\varphi}_0 \}$, $\IP_1(\hat{E}) = \mathrm{span}\{\hat{\varphi}_0, \ldots,  \hat{\varphi}_3 \}$, $\IP_2(\hat{E}) = \mathrm{span}\{\hat{\varphi}_0, \ldots,  \hat{\varphi}_9\}$, etc., is convenient for an~arbitrary order implementation.   The functions~$\hat{\varphi}_i$ are so-called \emph{modal} basis functions that do \emph{not} possess interpolation properties at nodes unlike Lagrangian\,/\,nodal basis functions, which are often used by the continuous finite element or nodal DG~methods. 

\subsection{Space discretization}\label{sec:spaceDiscretization}
For an~interior face $e\in\Gamma_h^\mathrm{int}$ shared by elements $E^-$ and $E^+$, we define a~unit normal vector $\normal_e$ that is exterior to $E^-$ pointing to~$E^+$.  Let the \emph{average}~$\avg{\,\cdot\,}$ and the \emph{jump}~$\jump{\,\cdot\,}$ of a~scalar quantity~$w$ on a~face~$e$ with normal~$\normal_e=\normal_{E^-}$ exterior to~$E^-\in\setE_h$ be defined by
\begin{equation}\label{eq:AverageAndJump}
\avg{w} ~\coloneqq~
\begin{Bmatrix*}[l]
\frac{1}{2} w|_{E^-} + \frac{1}{2}w|_{E^+}& \text{if}~e = \partial E^-\cap \partial E^+\\
w|_{E^-}                                  & \text{if}~e = \partial E^-\cap \partial\Omega_h
\end{Bmatrix*}
,\quad
\jump{w} ~\coloneqq~
\begin{Bmatrix*}[l]
w|_{E^-} - w|_{E^+}                   & \text{if}~e=\partial E^-\cap \partial E^+\\
w|_{E^-}                                  & \text{if}~e=\partial E^-\cap \partial\Omega_h
\end{Bmatrix*}.
\end{equation}
The definition of the jump is invariant with respect to the orientation of the face normals~$\normal_e$ since changing~$\normal_e$ to $-\normal_e$ means swapping~$E^-$ with $E^+$.
Furthermore, for $\vec{x}\in e$, we define a~first order \emph{upwinding} by
\begin{equation*}
c^\uparrow(t,\vec{x})~\coloneqq~\sig(z)\,c|_{E^-}(t,\vec{x}) + \big(1 - \sig(z)\big)\,c|_{E^+}(t,\vec{x}) \qquad\text{where}\qquad z~\coloneqq~\avg{\vec{v}(t,\vec{x})}\cdot\normal_e~,
\end{equation*}
and $\sig(z)\coloneqq 1/\big(1 + \exp(- k\,z)\big)$ is a~sigmoid function (we use $k=100$ in our simulations).  In particular, we have~$c^\uparrow = \avg{c}$ at points where the normal velocity vanishes, and approximately full upwinding for high velocities.  This pointwise upwind rule will be applied to quadrature points in the implementation. 

\subsubsection{Fully discrete problem}
In the numerical scheme, we have to limit the mobility from below and use $M_{+}(c)\coloneqq\max\{M(c),0\}$ instead, as we have to expect overshoots of~$c$ beyond the physical meaningful interval of~$[-1,1]$, cf.~\cref{sec:bulkShift}.  
\par
Let $p\in\IN_0$. The \emph{fully discrete problem} is defined as follows:
For $n\in\{1,\ldots,\Nst\}$, find $c_h^n, \mu_h^n\in\IP_p(\setE_h)$ such that $\forall w_h\in\IP_p(\setE_h)$,
\begin{subequations}\label{eq:CHmixed:fullydiscrete}
\begin{align}
\kappa\,a^\mathrm{diff}(1,c_h^n,w_h) + \kappa\,b^{\mathrm{diff},n}(c_h^n,w_h) - (\mu_h^n,w_h) +  \big(\Phi'_{+}(c_h^n),w_h\big) &~=~\kappa\,d^{\mathrm{diff},n}(w_h)  -  \big(\Phi'_{-}(c_h^{n-1}),w_h\big)~, \label{eq:CHmixed:fullydiscrete:a} 
\\
(c_h^n,w_h) +  \frac{\tau_n}{\Pe}\, a^\mathrm{diff}\Big(M_{+}(c_h^{n-1}), \mu_h^n,w_h\Big)  + \tau_n\, a^{\mathrm{adv},n}\big(c_h^n,w_h\big)
&~=~ (c_h^{n-1},w_h) + \tau_n\, b^{\mathrm{adv},n}(w_h)~,\label{eq:CHmixed:fullydiscrete:b} 
\\
\label{eq:CHmixed:fullydiscrete:init}
\big(c_h^0,w_h\big) &~=~ (c^0,w_h)~,
\end{align}
\end{subequations}
where $(\cdot\,,\,\cdot)$ denotes the $L^2(\Omega_h)$~inner product.  System~\eqref{eq:CHmixed:fullydiscrete} is linear in~$\mu_h^n$ and nonlinear in~$c_h^n$.  In the following, the symbol~$\grad$ is to be interpreted as the broken gradient, i.\,e., $\forall E\in\setE_h,~(\grad_h w)|_{E} \coloneqq \grad(w|_{E})$ while suppressing the subscript~$h$.
For positive integers~$p\in\IN$, let $a^\mathrm{diff} = a_\mathrm{DG}^\mathrm{diff}$ with
\begin{subequations}\label{eq:bilinearforms}
\begin{equation}\label{eq:bilinearform:DG}
a_\mathrm{DG}^\mathrm{diff}(z,c,w) ~\coloneqq~
 \int_{\Omega_h} z\,\grad c \cdot \grad w
-\sum_{e\in\Gamma_h^\mathrm{int}} \int_e \Big( \avg{z\,\grad c\cdot\normal_e} \jump{w}
+ \avg{z\,\grad w \cdot \normal_e} \jump{c}\Big)
+ \frac{\sigma}{h} \sum_{e\in\Gamma_h^\mathrm{int}}\int_e \jump{c}\jump{w}
\end{equation}
be the usual symmetric interior penalty Galerkin (SIPG)~bilinear form for the operator~$-\div(z\grad c)$, where $\sigma>0$ denotes the penalty parameter. With~$\sigma$ chosen sufficiently large, coercivity of the bilinear form~$a^\mathrm{diff}(z;\,\cdot\,,\,\cdot\,)$ can be ensured~\cite{Riviere2008}.
Due to the Dirichlet conditions~\eqref{eq:model:in} on the inflow boundary, and the homogeneous Neumann conditions~\eqref{eq:model:out} on the outflow boundary, the forms
\begin{align*}
b^{\mathrm{diff},n}(c,w) &~\coloneqq~  - \int_{\partial\Omega_h^\mathrm{in}(t_n)} \Big( \grad c\cdot\normal \,w
+  \grad w \cdot \normal\,c\Big)
+ \frac{\sigma}{h}  \int_{\partial\Omega_h^\mathrm{in}(t_n)} c\,w
~,
\\
d^{\mathrm{diff},n}(w) &~\coloneqq~  - \int_{\partial\Omega_h^\mathrm{in}(t_n)} \grad w \cdot \normal\,c^\mathrm{in}(t_n)
+ \frac{\sigma}{h}  \int_{\partial\Omega_h^\mathrm{in}(t_n)} c^\mathrm{in}(t_n) \,w
\end{align*}
are included in~\eqref{eq:CHmixed:fullydiscrete:a}.  Throughout this paper, we fix $\sigma\coloneqq 2^p$.  
The forms for the advection operator~$\div(c\,\vec{v})$ at~$t_n\in J$ subjected to the boundary condition~\eqref{eq:model:in}, \eqref{eq:model:out} are
\begin{align*}
a^{\mathrm{adv},n}(c,w) &~\coloneqq~  
-  \int_{\Omega_h} c\,\vec{v}(t_n) \cdot \grad w
+\sum_{e\in\Gamma_h^\mathrm{int}} \int_e c^\uparrow   \avg{\vec{v}(t_n) \cdot\normal_e} \jump{w} 
+ \int_{\partial\Omega_h^{\mathrm{out}}(t_n)} c\,\vec{v}(t_n)\cdot\normal\,w~,
\\
b^{\mathrm{adv},n}(w)
&~\coloneqq~ 
- \int_{\partial\Omega_h^{\mathrm{in}}(t_n)} c^\mathrm{in}(t_n)\,\vec{v}(t_n)\cdot\normal\,w~.
\end{align*}
The decomposition $\partial\Omega = \partial\Omega^\mathrm{in}(t_n) \cup \partial\Omega^\mathrm{out}(t_n)$ is determined in every time step.
\par
For elementwise constant approximations, i.\,e., for the FV~case~$p=0$, we define~$a^\mathrm{diff} = a_\mathrm{FV}^\mathrm{diff}$ by
\begin{equation}\label{eq:bilinearform:FV}
a_\mathrm{FV}^\mathrm{diff}(z,c,w) \coloneqq \frac{1}{h} \sum_{e\in\Gamma_h^\mathrm{int}}\int_e \avg{z}\jump{c}\jump{w}\;.
\end{equation}
\end{subequations}
For the zero order case, some of the above forms simplify: 
\begin{align*}
b^{\mathrm{diff},n}(c,w) &~=~ \frac{1}{h} \int_{\partial\Omega_h^\mathrm{in}(t_n)} c\,w   ~,
\qquad
d^{\mathrm{diff},n}(w) ~=~ \frac{1}{h}  \int_{\partial\Omega_h^\mathrm{in}(t_n)} c^\mathrm{in}(t_n)\,w
~,
\\
a^{\mathrm{adv},n}(c,w) &~=~  
\sum_{e\in\Gamma_h^\mathrm{int}} \int_e c^\uparrow \avg{\vec{v}(t_n) \cdot\normal_e} \jump{w}
+ \int_{\partial\Omega_h^{\mathrm{out}}(t_n)} c\,\vec{v}(t_n)\cdot\normal\,w~.
\end{align*}
The numerical scheme~\eqref{eq:CHmixed:fullydiscrete} realized with the forms for~$p=0$  is equivalent to a~cell-centered FV~approach~\cite[cf.][]{Chidyagwai2011}.

\subsubsection{Discrete solution properties}
It is desirable that fully discrete solutions of~\eqref{eq:CHmixed:fullydiscrete} satisfy discrete versions of the properties listed in~\cref{sec:SolutionProperties}.  For a~similar discretization of the non-advective CH~equation with constant mobility, i.\,e., \eqref{eq:model} with~$\vec{v}=\vec{0}$, \mbox{$\beta = 0$}, Aristotelous~\cite{AristotelousThesis2011} showed that~\eqref{eq:CHmixed:fullydiscrete} is uniquely solvable and that a~discrete version~$F_h$ of~$F$ as given in~\eqref{eq:HelmholtzFreeEnergy} is unconditionally non-increasing in time:
\begin{equation*}
\forall n\in\{1,\ldots,\Nst\}~,\quad F_h(c_h^n) ~\leq~ F_h(c_h^{n-1})~,
\end{equation*}
where the \emph{discrete Helmholtz free energy} is defined as\footnote{The Helmholtz free energies in~\cite{AristotelousThesis2011} are scaled by~$\epsilon\coloneqq\sqrt{\kappa}$, and hence the chemical potential by~$1/\epsilon$.}
\begin{equation}\label{eq:HelmholtzFreeEnergyDiscrete}
F_h(c) ~\coloneqq~ \big(\Phi(c),1\big) + \frac{\kappa}{2} \,a^\mathrm{diff}(1,c,c)~.
\end{equation}
\par
Choosing the function $w_h=1$ constant over $\Omega$ in \eqref{eq:CHmixed:fullydiscrete:b} yields
\begin{equation*}
\forall n\in\{1,\ldots,\Nst\}~,\quad \int_\Omega \frac{c_h^n - c_h^{n-1}}{\tau_n} + \int_{\partial\Omega_h^\mathrm{in}(t_n)} \vec{v}(t_n) \cdot \normal \, c^\mathrm{in}
+\int_{\partial\Omega_h^\mathrm{out}(t_n)} \vec{v}(t_n) \cdot \normal \, c_h^n ~=~ 0~,
\end{equation*}
which is the discrete version of the mass conservation equation~\eqref{eq:massConservation}.
\par
Note that since we use DG, we also have a~\emph{local mass conservation property} up to penalty:
For simplicity, consider an element that is not adjacent to the boundary, i.\,e., $E = E^{-}\in\setE_h$ with $\partial E\cap\partial\Omega_h=\emptyset$.
From \eqref{eq:massConservation} it follows
\begin{equation*}
\dd_t \int_\Omega c 
-
\frac{1}{\Pe}\int_{\partial E} M(c)\,\grad\mu\cdot\normal_E + \int_{\partial E} c\, \vec{v}\cdot\normal_E  ~=~ 0  ~,
\end{equation*}
where $\normal_E$ denotes the outward normal to~$E$.  Choosing $w_h = 1$ on~$E$ and zero elsewhere in~\eqref{eq:CHmixed:fullydiscrete:b}
\begin{equation*}
\int_E \frac{c_h^n-c_h^{n-1}}{\tau_n} - \frac{1}{\Pe} \int_{\partial E} \avg{ M_+(c_h^{n-1}) \nabla \mu_h^n \cdot \normal_E}
 +\int_{\partial E} (c_h^n)^\uparrow  \avg{\vec{v}(t_n) \cdot\normal_E}  ~=~ - \frac{\sigma}{h\,\Pe} \sum_{e\in\partial {E^{-}}}\int_e \Big(c_h^n|_{E^{-}} - c_h^n|_{E^{+}}\Big)~.
\end{equation*}
Note that the penalty term on the right-hand side contributes to a~local mass error, which is known and computable.
\par
To the best knowledge of the authors, there exists no numerical scheme for CH~systems that satisfies a~\emph{discrete maximum principle}.

\subsubsection{Discussion on other forms}
In the continuous formulation~\eqref{eq:model}, degenerate mobility suppresses the mass flux in the bulk phases, which motivates the use of bilinear forms that suppress the fluxes across elements in this situation:
\par
For the FV case with degenerate mobility, i.\,e., $p=0$ and $\beta=1$, the initial solution must contain at least one element with a~constant value unequal to~$-1$ and~$1$ such that the mobility is positive.  If this is not the case there will be no numerical flux across any element and the solution will not evolve in time.  
For the case of~$p=0$, as opposed to definition~\eqref{eq:bilinearform:FV} that uses the geometric mean~$\avg{z}$ on faces~$e\in\Gamma_h^\mathrm{int}$, the classical FV~method often uses some form of harmonic average instead, e.\,g.,
$\avgHarmonic{z} \coloneqq {2\,z^+\,z^-}/({z^+ + z^-})$, where~$z^\pm \coloneqq z|_{E^{\pm}}$.  On the other hand, one~generalization of the DG~form~$a_\mathrm{DG}^\mathrm{diff}$ is the weighted bilinear form introduced in~\cite{ErnSZ2009}.  For the case of~$z\in\IP_0(\setE_h)$ strictly positive and assuming no advection in the system, the so-called \emph{SWIP~bilinear form} reads\footnote{In~\cite{ErnSZ2009}, half of the harmonic average is used, however, we would like this term to be consistent with~$a_\mathrm{DG}^\mathrm{diff}$ for the case of~$z\equiv 1$.} for $c, w\in\IP_p(\setE_h)$ with $p\geq 1$:
\begin{equation*}
a^\mathrm{W}(z,c,w) \coloneqq \int_{\Omega_h} z\,\grad c \cdot \grad w
-\sum_{e\in\Gamma_h^\mathrm{int}} \int_e \Big(\avg{z\,\grad c\cdot\normal_e}_\omega \jump{w}
+  \avg{z\,\grad w \cdot \normal_e}_\omega \jump{c}\Big)
+ \frac{\sigma }{h}  \sum_{e\in\Gamma_h^\mathrm{int}}\int_e \avgHarmonic{z} \jump{c}\jump{w} ,
\end{equation*}
where the \emph{weighted average} on a~face~$e\in\Gamma_h^\mathrm{int}$ is defined as
\begin{equation*}
\avg{w}_\omega ~\coloneqq~ \omega^- w|_{E^-} + \omega^+ w|_{E^+}~,\qquad\text{where}\qquad
\omega^- ~\coloneqq~ \frac{z^+}{z^- + z^+}~,
\quad \omega^+~\coloneqq 1-\omega^-~.
\end{equation*}
\par
In the case of~FV with harmonic average as well as when using the SWIP DG~form, numerical fluxes vanish across faces that are adjacent to one element with zero mobility.  This means that the interface area ($-1<c<1$) is trapped between two rigid phases making these forms practically unfeasible for CH~simulations with degenerate mobility.  Using~\eqref{eq:bilinearform:DG} with averaged~$z$ in the penalty term similar to~\eqref{eq:bilinearform:FV} turned out to yield a~numerical ill-posed scheme (with regard to both non-linear and linear solver), probably since there is too less diffusion in the system.

\subsubsection{Matrix formulation}
Let $c_h^n,\mu_h^n\in\IP_p(\setE_h)$ have the representations
\begin{equation}\label{eq:representationVectors}
c_h^n(\vec{x}) ~=~ \sum_{k=0}^{\Nel-1} \sum_{j=0}^{\Nloc-1}
c_{kj}^n\,\varphi_{kj}(\vec{x}) ~,\qquad 
\mu_h^n(\vec{x}) ~=~ \sum_{k=0}^{\Nel-1} \sum_{j=0}^{\Nloc-1}
\mu_{kj}^n\,\varphi_{kj}(\vec{x})~,
\end{equation}
where the basis~$\{\varphi_{kj}\}$ was defined in~\cref{sec:ApproximationSpaces}.
Then, \eqref{eq:CHmixed:fullydiscrete} is equivalent to the following system of equations:
\begin{equation}\label{eq:cahi:impl:mixedLinearSystem}
\begin{bmatrix}
\kappa\,(\vecc{A}^\mathrm{diff} +\vecc{B}^{\mathrm{diff},n}) & - \vecc{M}
\\
\vecc{M} + \tau_n \, \vecc{A}^{\mathrm{adv},n} & \frac{\tau_n}{\Pe}\, \vecc{A}^\mathrm{diff}_M (\vec{X}_c^{n-1})
\end{bmatrix}
\,
\begin{bmatrix}
\vec{X}_c^n\\
\vec{X}_\mu^n
\end{bmatrix}
+
\begin{bmatrix}
\vec{E}_{\Phi'_+}(\vec{X}_c^n)
\\
\vec{0}
\end{bmatrix}
 ~=~ 
\begin{bmatrix}
\vec{F}_c^n
\\
\vec{F}_\mu^n
\end{bmatrix}
~\coloneqq~
\begin{bmatrix}
\kappa\,\vec{D}^{\mathrm{diff},n}  - \vec{E}_{\Phi'_-}(\vec{X}_c^{n-1}) 
\\
\vecc{M}\,\vec{X}_c^{n-1}  + \tau_n \, \vec{B}^{\mathrm{adv},n}
\end{bmatrix}
\end{equation}
\ifdraft{
{\gray FOR DOCUMENTATION, WILL DISAPPEAR WITHOUT 'draft' OPTION\\
With source terms~$f, g_c, g_\mu$ we had
\begin{align*}
\vec{F}_c^n    &~=~  \kappa\,\vec{D}^{\mathrm{diff},n}  - \vec{E}_{\Phi'_-}(\vec{X}_c^{n-1})  +  \kappa\,\vec{G}_{g_c}^n 
\qquad \text{with} \qquad  [\vec{G}_{g_c}^n ]_{\Nloc k + i}~\coloneqq~\int_{\partial\Omega_h^\mathrm{out}(t_n)} g_c(t_n)\,\varphi_{ki}~,
\\
\vec{F}_\mu^n  &~=~  \vecc{M}\,\vec{X}_c^{n-1}  + \tau_n \, \vec{B}^{\mathrm{adv},n} +    \tau_n\,\vec{G}_f^n   + \frac{\tau_n}{\Pe}\,\vec{G}_{g_\mu}^n
~~ \text{with} ~~  [\vec{G}_{g_\mu}^n ]_{\Nloc k + i} \coloneqq \int_{\partial\Omega_h} g_\mu(t_n)\,\varphi_{ki}
, ~~[\vec{G}_f^n ]_{\Nloc k + i}                    \coloneqq \big(f(t_n),\varphi_{ki}\big) 
\end{align*}}
}{}
with~$\vec{X}_c^0$ satisfying the initial condition~\eqref{eq:CHmixed:fullydiscrete:init}. 
Here, $[\vec{X}_c^n]_{l\Nloc + j}\coloneqq c_{lj}^n$ and $[\vec{X}_\mu^n]_{l\Nloc + j}\coloneqq \mu_{lj}^n$
denote the vectors of degrees of freedom with respect to $c_h^n$ and $\mu_h^n$ and 
\begin{align*}
&[\vecc{A}^\mathrm{diff}]_{k\Nloc + i, l \Nloc + j}                ~\coloneqq~  a^\mathrm{diff}(1, \varphi_{lj},   \varphi_{ki})~,
& 
&[\vecc{A}^\mathrm{diff}_M(\vec{X}_c^n)]_{k\Nloc + i, l \Nloc + j} ~\coloneqq~  a^\mathrm{diff}\Big(M_{+}(c_h^n), \varphi_{lj},   \varphi_{ki} \Big)~,
\\
&[\vecc{B}^{\mathrm{diff},n}]_{k\Nloc + i, l \Nloc + j}                ~\coloneqq~  b^{\mathrm{diff},n}(\varphi_{lj},   \varphi_{ki})~,
&
&[\vecc{A}^{\mathrm{adv},n}]_{k\Nloc + i, l \Nloc + j}             ~\coloneqq~  a^{\mathrm{adv},n}(\varphi_{lj},   \varphi_{ki})~,
\\
&[\vec{D}^{\mathrm{diff},n}]_{k\Nloc + i}                          ~\coloneqq~  d^{\mathrm{diff},n}(\varphi_{ki})~,
&
&[\vec{B}^{\mathrm{adv},n}]_{k\Nloc + i}                           ~\coloneqq~  b^{\mathrm{adv},n}(\varphi_{ki})~,
\\
&[\vec{E}_{\Phi'_\pm}(\vec{X}_c^n)]_{k\Nloc + i}                   ~\coloneqq~ (\Phi'_\pm(c_h^n) , \varphi_{ki})
&
\end{align*}
for $k,l\in\{0,\ldots,\Nel-1\}$ and $i,j\in\{0,\ldots,\Nloc-1\}$, while the mass matrix~$\vecc{M}$ was defined in~\eqref{eq:MassMatrix}.
\par
The matrix~$\vecc{A}^\mathrm{diff}$ of size~$\Nel\Nloc \times \Nel\Nloc$ is sparse, symmetric, positive semidefinite if the penalty parameter~$\sigma$ is sufficiently large, and has a~rank deficiency of one~\emph{per disconnected subdomain} since it is associated with a~Poisson problem subjected to homogeneous Neumann boundary conditions (as a such, solutions would be defined up to a constant on each subdomain). Associated with Dirichlet condition yielded a~condition number of order~$O(1/h^2)$~\cite{Castillo2002Performance}.  The matrix~$\vecc{A}^\mathrm{diff}_M(\vec{X}_c^n)$ may have additional zero eigenvalues due to the possibly locally vanishing mobility coefficient.

\paragraph{System reduction}
The structure of~\eqref{eq:cahi:impl:mixedLinearSystem} is suitable to build the Schur complement of the upper-right block:
The first block-equation yields 
\begin{equation}\label{eq:cahi:impl:linearSystem:a}
\vec{X}_\mu^n~ =~ \vecc{M}^{-1}\Big( \kappa\,(\vecc{A}^\mathrm{diff} + \vecc{B}^{\mathrm{diff},n})\vec{X}_c^n   +  \vec{E}_{\Phi'_+}(\vec{X}_c^n) - \vec{F}_c^n  \Big)~.
\end{equation}  
Substitution into the second block-equation yields the linear system
\begin{multline}\label{eq:cahi:impl:linearSystem:b}
\Big(\vecc{M}  + \tau_n \, \vecc{A}^{\mathrm{adv},n} + \frac{\kappa\,\tau_n}{\Pe}\,\vecc{A}^\mathrm{diff}_M (\vec{X}_c^{n-1})\,\vecc{M}^{-1} \,(\vecc{A}^\mathrm{diff} + \vecc{B}^{\mathrm{diff},n} ) \Big) \, \vec{X}_c^n 
+ \frac{\tau_n}{\Pe}\, \vecc{A}^\mathrm{diff}_M (\vec{X}_c^{n-1})\,\vecc{M}^{-1} \, \vec{E}_{\Phi'_+}(\vec{X}_c^n) 
\\
~=~\vec{F}_\mu^n + \frac{\tau_n}{\Pe}\,\vecc{A}^\mathrm{diff}_M (\vec{X}_c^{n-1})\,\vecc{M}^{-1}\,\vec{F}_c^n~.
\end{multline}
Thus solving block-system~\eqref{eq:cahi:impl:mixedLinearSystem} of size~$2\Nel\Nloc \times 2\Nel\Nloc$ is reduced to system~\eqref{eq:cahi:impl:linearSystem:b} of size~$\Nel\Nloc \times \Nel\Nloc$ and an~optional direct computation of~$\mu_h^n$ by~\eqref{eq:cahi:impl:linearSystem:a}.

\subsubsection{Linearization}

\paragraph{Inexact Newton method}
We have to solve a~sequence of nonlinear systems~\eqref{eq:cahi:impl:linearSystem:b}, the residual form of which reads
\begin{equation}\label{eq:cahi:impl:linearSystem:c}
\underbrace{\Bigg(\frac{h^6\,\Pe}{64\,\tau_n}\vecc{I}  + \frac{h^3}{8}\vecc{A}^{\mathrm{adv},n} + \kappa\,\vecc{A}^\mathrm{diff}_M (\vec{X}_c^{n-1})\,(\vecc{A}^\mathrm{diff} + \vecc{B}^{\mathrm{diff},n} ) \Bigg) \, \vec{X}_c^n }_{\eqqcolon\,\mathcal{L}\, \vec{X}_c^n }
+ \underbrace{\vecc{A}^\mathrm{diff}_M (\vec{X}_c^{n-1})\,\vec{E}_{\Phi'_+}(\vec{X}_c^n)}_{\eqqcolon\,\mathcal{N}(\vec{X}_c^n)}
\underbrace{ {}- \vecc{A}^\mathrm{diff}_M (\vec{X}_c^{n-1})\,\vec{F}_c^n 
               - \frac{h^3\,\Pe}{8\,\tau_n}\,\vec{F}_\mu^n }_{\eqqcolon\,\mathcal{C}^n} = \vec{0}\,,
\end{equation}
containing a~\emph{linear part}~$\mathcal{L}\, \vec{X}_c^n$, a~\emph{nonlinear part}~$\mathcal{N}(\vec{X}_c^n)$, and a~\emph{constant part}~$\mathcal{C}^n$.
Here, we used $\vecc{M} = \frac{h^3}{8}\vecc{I}$ due to the orthonormal basis on~$\hat{E}$, cf.~\eqref{eq:MassMatrix}, $\vecc{I}$ denoting the unit matrix.
In the case of degenerate mobility, even for the case of no advection, the system matrix~$\mathcal{L}$ is not symmetric (even though $\vecc{A}^\mathrm{diff}_M(\vec{X}_c^n)$, $\vecc{A}^\mathrm{diff}$, and $\vecc{B}^{\mathrm{diff},n}$ are).
The condition number of system~\eqref{eq:cahi:impl:linearSystem:c} is of order~$O(1/h^4)$ as the condition number of~$\vecc{A}^\mathrm{diff}_M(\vec{X}_c^n)$ and $\vecc{A}^\mathrm{diff}$ scales as~$1/h^2$~\cite[cf.][]{Hanisch1993,Aristotelous2013}.
At this linear algebra level, it it easy to see that the condition number of~$\mathcal{L}$ can be decreased by decreasing~$\kappa$ or~$\tau_n$ (keeping in mind that~$\kappa$ is limited by~$h$, cf.~\cref{sec:ModelProblem}).
\par
For $\vec{F}:\IR^N\ni\vec{X}\mapsto\vec{F}(\vec{X})\in\IR^N$ let $\vecc{D}\vec{F} \coloneqq \partial\vec{F} /\partial\vec{X} = [\partial F_j /\partial X_i]_{ij}$, i.\,e.,~$\vecc{D}\vec{F}:\IR^N\rightarrow\IR^{N,N}$ denote the Jacobian of~$\vec{F}$.
For a~fixed step~$n$, the \emph{Newton method} applied to~\eqref{eq:cahi:impl:linearSystem:c} reads:
\par
Set $\vec{Y}^0\coloneqq\vec{X}_c^{n-1}$. Iteratively, for~$s=1,2,\ldots$, seek~$\vec{U}^s$ such that
\begin{subequations}\label{eq:newton}
\begin{align}\label{eq:newton:a}
\big( \mathcal{L} + \mathrm{D}\mathcal{N}(\vec{Y}^{s-1})\big)\,\vec{U}^s ~=~& \big(  \mathcal{L}\,\vec{Y}^{s-1} +  \mathcal{N}(\vec{Y}^{s-1}) + \mathcal{C}^n \big),\qquad \text{and set} 
\\\label{eq:newton:b}
\vec{Y}^{s} ~\coloneqq~& \vec{Y}^{s-1} - \vec{U}^s
\end{align}
\end{subequations}
as long as the norm of the system's residual does not fall below a~given tolerance, i.\,e., stop as soon as
\begin{equation}\label{eq:Newton:stop}
\norm{ \mathcal{L}\,\vec{Y}^{s}+  \mathcal{N}(\vec{Y}^{s}) + \mathcal{C}^n}{L_h^2}  ~\leq~\mathrm{tol_{abs}} + \mathrm{tol_{rel}}\,\norm{ \mathcal{L}\,\vec{Y}^{0} +  \mathcal{N}(\vec{Y}^{0}) + \mathcal{C}^n}{L_h^2}
\end{equation}
with given tolerances~$\mathrm{tol_{abs}}, \mathrm{tol_{rel}}>0$ and norm~$\norm{\,\cdot\,}{L_h^2} \coloneqq h^{3/2}\norm{\,\cdot\,}{2}$.  
We use the discrete $L^2$~norm to have a~similar scaling with respect to~$h$ as the~\mbox{$L^2(\Omega_h)$-norm}, i.\,e., $\norm{\vec{X}_c^n}{L_h^2}\sim \norm{c_h^n}{L^2(\Omega_h)}$.
If~\eqref{eq:Newton:stop} is met then we accept the iterative solution as solution of the current time step: $\vec{X}_c^n\coloneqq \vec{Y}^s$.
Before applying a~non-linear step, we check the initial absolute residual: if
\begin{equation}\label{eq:Newton:stopInitial}
\norm{ \mathcal{L}\,\vec{X}_c^{n-1}+  \mathcal{N}(\vec{X}_c^{n-1}) + \mathcal{C}^n}{L_h^2}  ~\leq~\mathrm{tol_{abs}}
\end{equation}
holds a~stationary state has been reached and the discrete solution will stay constant.
The Armijo rule is used to damp the Newton step length in order to guarantee a~decreasing residual and to enlarge the area of convergence~\cite{Armijo1966,Kelley2003}.
The scheme is expected to converge locally quadratic with a~local area of convergence around the fixed point.

\paragraph{Computation of the Jacobian and linear solver} We have to compute the Jacobian of the nonlinear term~$\mathcal{N}(\vec{X}_c^n)$ appearing in~\eqref{eq:cahi:impl:linearSystem:c}.  Using the matrix identity
$\dd_{\vec{x}}\big(\vecc{A}\,\vec{v}(\vec{x})\big)= \vecc{A}\,\dd_{\vec{x}} \vec{v}(\vec{x})$
for matrices~$\vecc{A}$ not depending on~$\vec{x}$, we obtain
\begin{equation*}
\mathrm{D}\mathcal{N} ~=~ \mathrm{D}\Big( \vecc{A}_M^\mathrm{diff} (\vec{X}_c^{n-1}) \, \vec{E}_{\Phi'_+}(\vec{X}_c^n) \Big) =  \vecc{A}_M^\mathrm{diff} (\vec{X}_c^{n-1})  \,  \mathrm{D}\vec{E}_{\Phi'_+}(\vec{X}_c^n) ~.
\end{equation*}
The matrix~$\mathrm{D}\vec{E}_{\Phi'_+}(\vec{X}_c^n)$ is block-diagonal with symmetric diagonal blocks: recalling~\eqref{eq:representationVectors}, we have
\begin{equation*}
\frac{\dd [\vec{E}_{\Phi'_+}(\vec{X}_c^n	)]_{k\Nloc + j}}{\dd [\vec{X}_c^n]_{k\Nloc + i} }
~=~ \int_{E_k} \Phi''_+\Big(\sum_{s = 0}^{\Nloc - 1} c_{ks}^n\,\varphi_{ks}\Big)\,\varphi_{kj}\,\varphi_{ki}
~\eqqcolon~ [\vecc{M}_{\Phi''_+}(\vec{X}_c^n)]_{k\Nloc + j, k\Nloc + i}
~=~ [\vecc{M}_{\Phi''_+}(\vec{X}_c^n)]_{k\Nloc + i, k\Nloc + j}
\end{equation*}
for $k\in\{0,\ldots,\Nel-1\}$, $i,j\in\{0,\ldots,\Nloc-1\}$, where $\vecc{M}_{\Phi''_+}(\vec{X}_c^n)$ is a~weighted mass matrix.
Thus in the Newton scheme above, we have
\begin{equation*}
 \mathrm{D}\mathcal{N}(\vec{Y}^{s-1}) ~=~  \vecc{A}_M^\mathrm{diff} (\vec{X}_c^{n-1}) \, \vecc{M}_{\Phi''_+}(\vec{Y}^{s-1})~.
\end{equation*}
Since~$\vecc{M}_{\Phi''_+}(\vec{Y}^{s-1})$ is block-diagonal and $\vecc{A}_M^\mathrm{diff} (\vec{X}_c^{n-1})$ has to be assembled for other terms too, the additional cost accompanied by Newton's method is small.
\par
The linear system~\eqref{eq:newton:a}, which has to be solved in every Newton step, contains matrix--matrix products in the system matrix.  Building the product of two sparse matrices is computationally expensive and furthermore yields a~significantly denser matrix.  In order to avoid building this resulting matrix, iterative solvers can be used that only require the acting of a~matrix on an~iteration vector, e.\,g., $\vec{z} = \vecc{A}\vecc{B}\vec{x} \Leftrightarrow\vec{y} = \vecc{B}\vec{x}, \vec{z} = \vecc{A}\vec{y}$.

\section{Simulation results}\label{sec:SimulationResults}
In~\cref{sec:convergenceStudy} we validate the numerical scheme as defined in~\cref{sec:spaceDiscretization} for the case of constant mobility and nonlinear mobility by numerically showing the expected orders of convergence for the discretization error.  Mass is conserved (up to arithmetic precision) and a~discrete version~$F_h$ of the total free energy~$F$ dissipates in time throughout all simulations, cf.~\cref{sec:CHmixed:spinodalDecomposition}.
All simulations are performed with constant and degenerate mobility and a~chemical energy density defined by
\eqref{eq:PhiOfc}.

\begin{table}[ht!]
\centering
\footnotesize
\begin{tabularx}{\linewidth}{@{}ccc|RcRc|RcRc@{}}
\toprule
$p$ & $\Nel$ & $\Nst$ 
&  \multicolumn{2}{r}{$\|c_h^{\Nst} - c(T)\|_{L^2(\Omega_h)}$} & \multicolumn{2}{r|}{$\brokennorm{\grad \big( c_h^{\Nst} - c(T)\big)}{\vec{H}^0(\setE_h)}$}
&  \multicolumn{2}{r}{$\|c_h^{\Nst} - c(T)\|_{L^2(\Omega_h)}$} & \multicolumn{2}{r}{$\brokennorm{\grad \big( c_h^{\Nst} - c(T)\big)}{\vec{H}^0(\setE_h)}$} \\
\midrule
$0$ & $1^3$  & $1$  & $3.364\E{-1}$ & ( --- )    & $9.076\E{-1}$ & ( --- )    &  n\,/\,a\textsuperscript{*}       & (n\,/\,a\textsuperscript{*})   & n\,/\,a\textsuperscript{*}       & (n\,/\,a\textsuperscript{*})\\
{}  & $2^3$  & $2$  & $7.280\E{-1}$ & ($-1.114$) & $1.197\E{-0}$ & ($-0.399$) &  n\,/\,a\textsuperscript{*}       & (n\,/\,a\textsuperscript{*})   & n\,/\,a\textsuperscript{*}       & (n\,/\,a\textsuperscript{*})\\
{}  & $4^3$  & $4$  & $1.339\E{-1}$ & ($2.443$)  & $1.415\E{-0}$ & ($-0.241$) & $1.591\E{-1}$ & ( --- )   & $1.415\E{-0}$ & ( --- )\\
{}  & $8^3$  & $8$  & $5.494\E{-2}$ & ($1.285$)  & $1.415\E{-0}$ & ($0.000$)  & $5.745\E{-2}$ & ($1.470$) & $1.415\E{-0}$ & ($0.000$) \\
{}  & $16^3$ & $16$ & $2.603\E{-2}$ & ($1.078$)  & $1.415\E{-0}$ & ($0.000$)  & $2.640\E{-2}$ & ($1.122$) & $1.415\E{-0}$ & ($0.000$)\\
{}  & $32^3$ & $32$ & $1.284\E{-2}$ & ($1.020$)  & $1.415\E{-0}$ & ($0.000$)  & $1.291\E{-2}$ & ($1.032$) & $1.415\E{-0}$ & ($0.000$)\\
{}  & $64^3$ & $64$ & $6.394\E{-3}$ & ($1.006$)  & $1.415\E{-0}$ & ($0.000$)  & $6.412\E{-3}$ & ($1.001$) & $1.415\E{-0}$ & ($0.000$)\\
\midrule
$1$ & $1^3$  & $1$    & $3.364\E{-1}$ & ( --- )   & $9.076\E{-1}$ & ( --- )     & $3.364\E{-1}$ & ( --- )   & $9.076\E{-1}$ & ( --- )\\
{}  & $2^3$  & $4$    & $1.302\E{-1}$ & ($1.369$) & $1.452\E{-0}$ & ($-0.678$)  & $1.121\E{-1}$ & ($1.585$) & $1.406\E{-0}$ & ($-0.632$)\\
{}  & $4^3$  & $16$   & $8.244\E{-2}$ & ($0.659$) & $1.083\E{-0}$ & ($0.423$)   & $9.045\E{-2}$ & ($0.310$) & $1.104\E{-0}$ & ($0.349$)\\
{}  & $8^3$  & $64$   & $3.181\E{-2}$ & ($1.374$) & $6.065\E{-1}$ & ($0.836$)   & $3.327\E{-2}$ & ($1.443$) & $6.106\E{-1}$ & ($0.854$)\\
{}  & $16^3$ & $256$  & $9.542\E{-3}$ & ($1.737$) & $2.981\E{-1}$ & ($1.025$)   & $1.003\E{-2}$ & ($1.729$) & $2.994\E{-1}$ & ($1.028$) \\
{}  & $32^3$ & $1024$ & $2.545\E{-3}$ & ($1.907$) & $1.432\E{-1}$ & ($1.058$)   & $2.684\E{-3}$ & ($1.896$) & $1.434\E{-1}$ & ($1.061$)\\
{}  & $64^3$ & $4096$ & $6.484\E{-4}$ & ($1.973$) & $7.016\E{-2}$ & ($1.029$)   & n\,/\,a\textsuperscript{**} & (n\,/\,a\textsuperscript{**})   &   n\,/\,a\textsuperscript{**} & (n\,/\,a\textsuperscript{**})  \\
\midrule
$2$ & $1^3$  & $1$        & $4.191\E{-2}$ & ( --- )   & $8.230\E{-1}$ & ( --- )   &  $4.191\E{-2}$ &  ( --- )    &  $8.230\E{-1}$ &  ( --- )\\
{}  & $2^3$  & $8$        & $5.543\E{-2}$ & ($-0.403$)   & $8.773\E{-1}$ & ($-0.092$)  &  $7.030\E{-2}$ &  ($-0.746$) &  $8.979\E{-1}$ &  ($-0.126$) \\
{}  & $4^3$  & $64$       & $2.328\E{-2}$ & ($1.252$) & $4.910\E{-1}$ & ($0.837$) &  $2.349\E{-2}$ &  ($1.581$)  &  $4.913\E{-1}$ &  ($0.870$)\\
{}  & $8^3$  & $512$      & $2.979\E{-3}$ & ($2.966$) & $1.366\E{-1}$ & ($1.846$) & $2.983\E{-3}$  &  ($2.977$)  &  $1.366\E{-1}$ &  ($1.847$)\\
{}  & $16^3$ & $4\,096$   & $3.222\E{-4}$ & ($3.209$) & $3.472\E{-2}$ & ($1.976$) &  $3.159\E{-4}$ &  ($3.239$)  &  $3.475\E{-2}$ &  ($1.975$)\\
{}  & $32^3$ & $32\,768$  & $3.682\E{-5}$ & ($3.129$) & $8.747\E{-3}$ & ($1.989$) &  $3.739\E{-5}$ &  ($3.079$)  &  $8.756\E{-3}$ &  ($1.989$)\\
{}  & $64^3$ & $262\,144$ & $4.504\E{-6}$ & ($3.031$) & $2.191\E{-3}$ & ($1.997$) & n\,/\,a\textsuperscript{**}  & (n\,/\,a\textsuperscript{**})  & n\,/\,a\textsuperscript{**}  &  (n\,/\,a\textsuperscript{**})\\
\bottomrule
\end{tabularx}
\caption{Discretization errors and estimated minimum orders of convergence \emph{(in parentheses)} for a~sequence of temporal and spatial refinements for the constant mobility case~$\beta=0$~\emph{(middle column)} and the degenerate mobility case~$\beta=1$~\emph{(right column)}.}
\label{tab:caha:convergence:1}
\end{table}

\subsection{Convergence study}\label{sec:convergenceStudy}
Consider the CH~equation~\eqref{eq:model} on the time-space cylinder $J\times\Omega_h = (0,1)\times(0,1)^3$ using~$\kappa=1$ and~$\vec{v}=\vec{0}$.
For prescribed stationary solutions that are element of~$\IP_p$, $p\in\{0,\ldots,2\}$, we obtain zero error up to arithmetic accuracy for all refinement levels.
We prescribe a~non-polygonal solution~$c(t, \vec{x}) \coloneqq  \exp(-t)\sin(2\,\pi\,x)\sin(2\,\pi\,y)\sin(2\,\pi\,z)$.
\par
We start a~refinement sequence with a~time step size of~$\tau = T = 1$ using a~triangulation with~$h=1$, while dividing~$h$ by two in each refinement level and $\tau$ by two, four, and eight for $p=0,1,2$, respectively, due to the first order discretization in time.  The discretization errors~$\|c_h^{\Nst} - c(T)\|_{L^2(\Omega_h)}$ and 
\begin{equation*}
\brokennorm{\grad \big( c_h^{\Nst} - c(T)\big)}{\vec{H}^0(\setE_h)}  \coloneqq \left(\sum_{E\in\setE_h} \norm{\grad \big( c_h^{\Nst} - c(T)\big)}{\vec{L}^2(E)}\right)^{1/2}
\end{equation*}
together with corresponding estimated convergence orders are listed in~\cref{tab:caha:convergence:1} using~$\IP_p(\setE_h)$ ($\beta$ was given in~\eqref{eq:mobility}).  The cases marked with~\enquote{\mbox{n\,/\,a\textsuperscript{*}}} are not feasible since the initial data consists only of constant elements with zero mobility and thus there is no evolution.  The linear systems that arise in the~\enquote{\mbox{n\,/\,a\textsuperscript{**}}} cases are too ill-conditioned to find a~solution. 
Optimal convergence is obtained in all cases.

\begin{figure}[hb!]
\centering
\begin{tabularx}{\linewidth}{@{}cCCCCC@{}}
\raisebox{2em}{\rotatebox{90}{$p = 0$, $\beta = 0$}} &
\includegraphics[width=\linewidth]{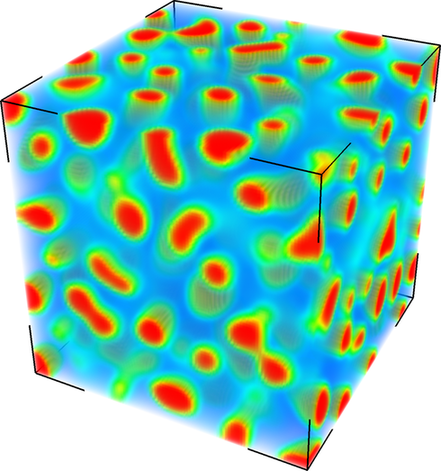} & 
\includegraphics[width=\linewidth]{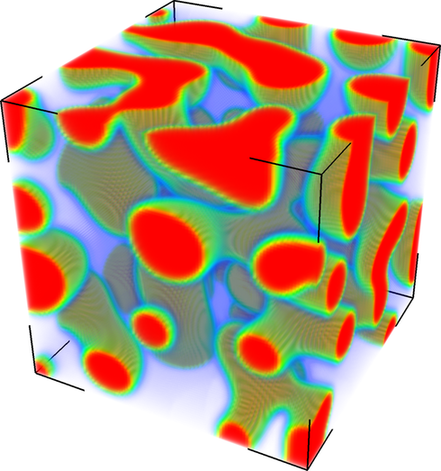} & 
\includegraphics[width=\linewidth]{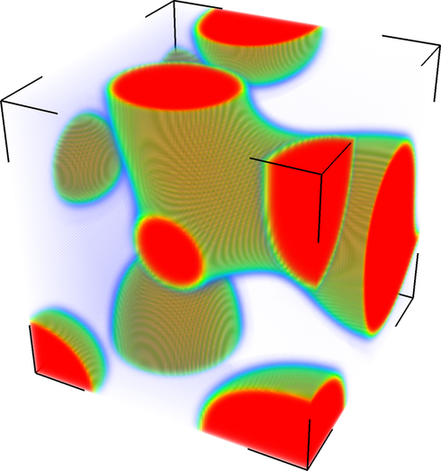} & 
\includegraphics[width=\linewidth]{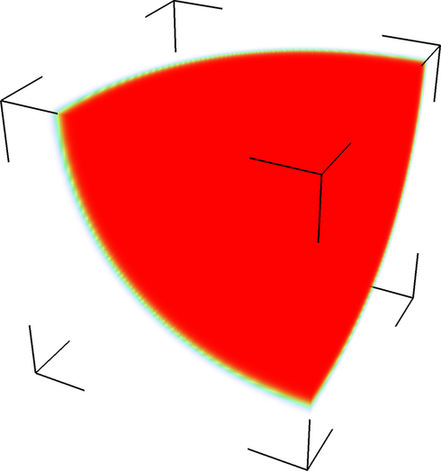} & 
\includegraphics[width=\linewidth]{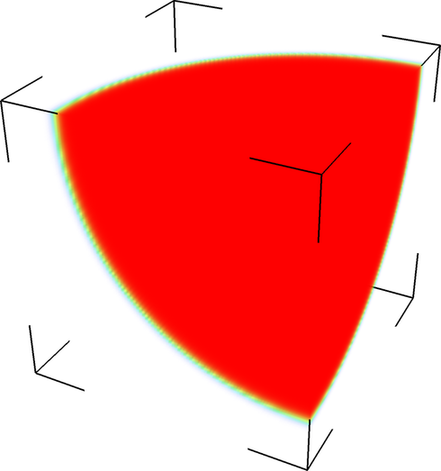}
\end{tabularx}
\\
\noindent\rule{\linewidth}{0.4pt}
\\[0.5\baselineskip]
\begin{tabularx}{\linewidth}{@{}cCCCCC@{}}
\raisebox{2em}{\rotatebox{90}{$p = 0$, $\beta = 1$}} &
\includegraphics[width=\linewidth]{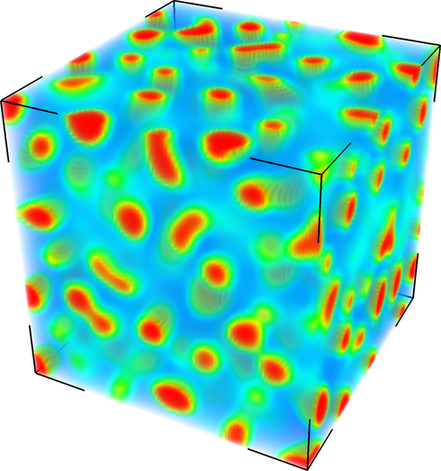} & 
\includegraphics[width=\linewidth]{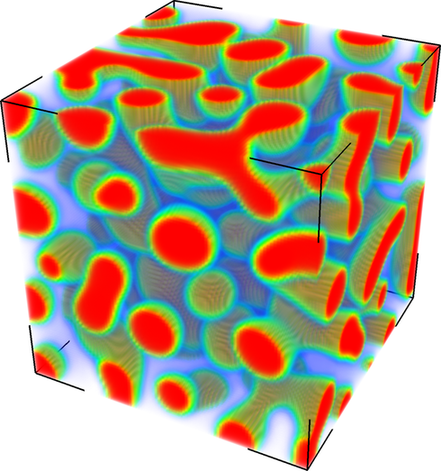} & 
\includegraphics[width=\linewidth]{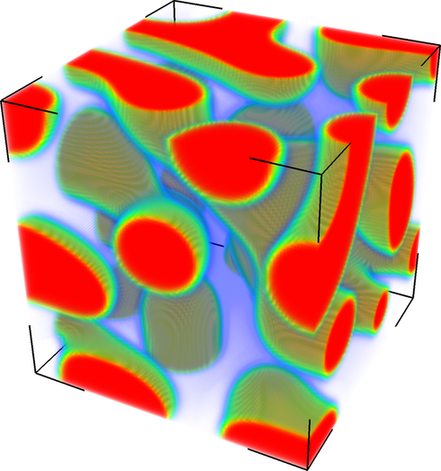} & 
\includegraphics[width=\linewidth]{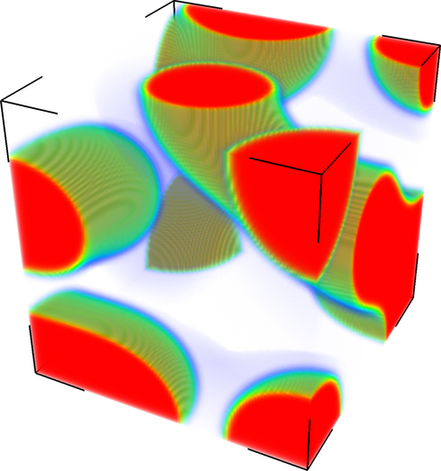} &
\includegraphics[width=\linewidth]{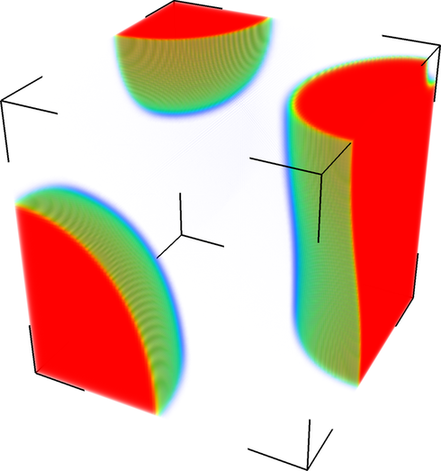}
\end{tabularx}
\\
\noindent\rule{\linewidth}{0.4pt}
\\[0.5\baselineskip]
\begin{tabularx}{\linewidth}{@{}cCCCCC@{}}
\raisebox{2em}{\rotatebox{90}{$p = 1$, $\beta = 0$}} &
\includegraphics[width=\linewidth]{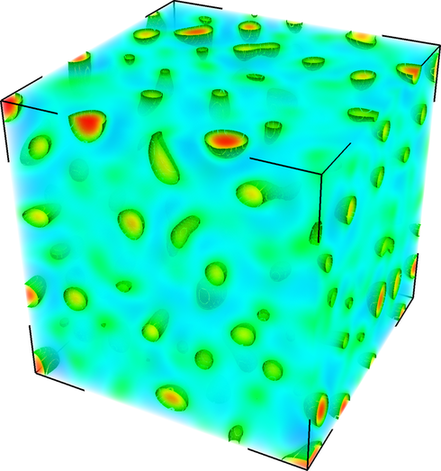} & 
\includegraphics[width=\linewidth]{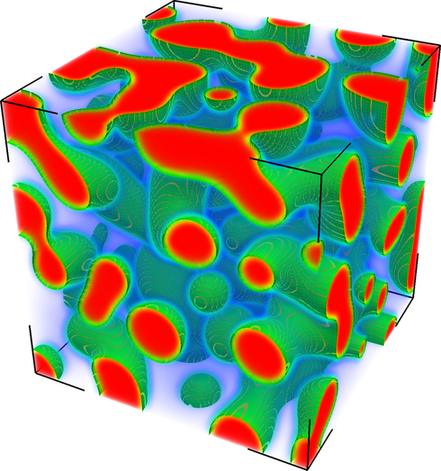} & 
\includegraphics[width=\linewidth]{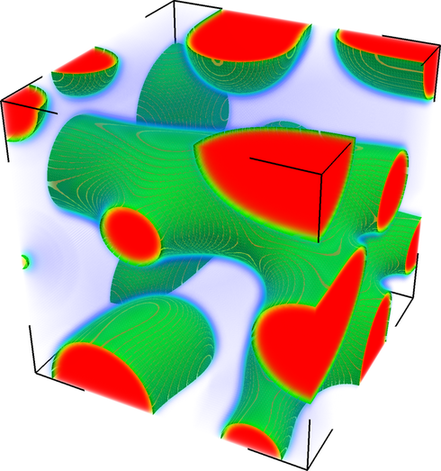} & 
\includegraphics[width=\linewidth]{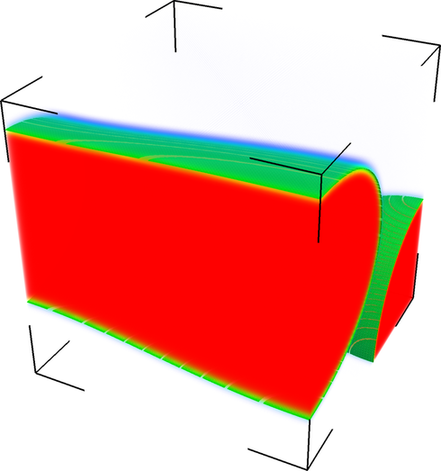} & 
\includegraphics[width=\linewidth]{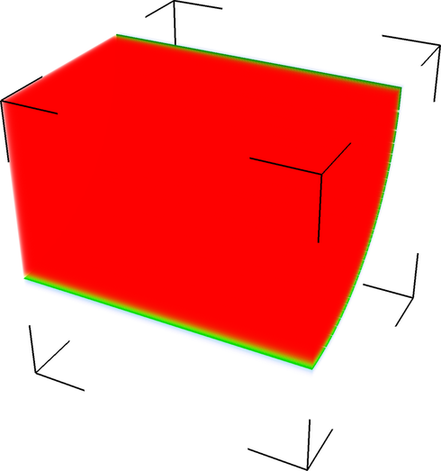} 
\end{tabularx}
\\
\noindent\rule{\linewidth}{0.4pt}
\\[0.5\baselineskip]
\begin{tabularx}{\linewidth}{@{}cCCCCC@{}}
\raisebox{2em}{\rotatebox{90}{$p = 1$, $\beta = 1$}} &
\includegraphics[width=\linewidth]{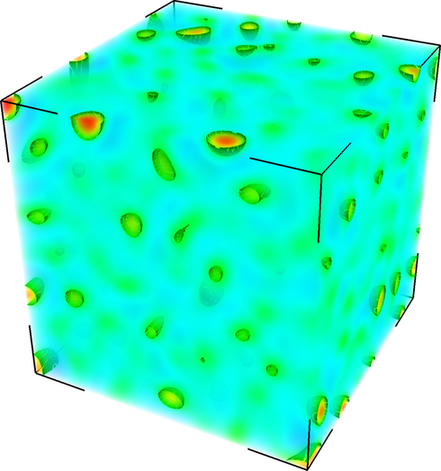} & 
\includegraphics[width=\linewidth]{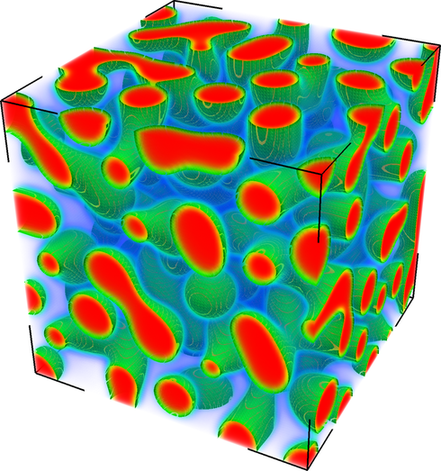} &
\includegraphics[width=\linewidth]{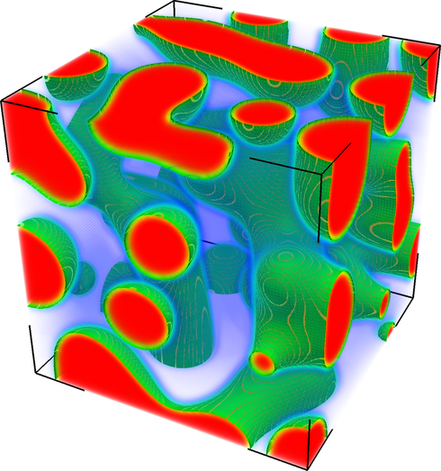} &
\includegraphics[width=\linewidth]{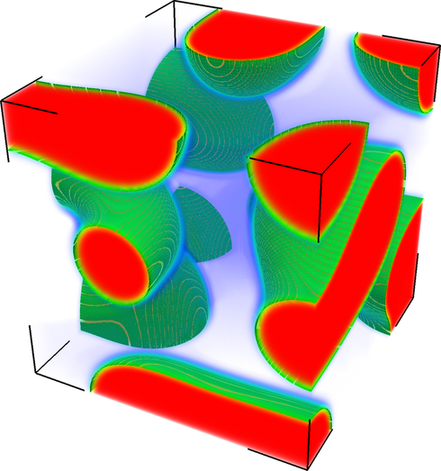} &
\includegraphics[width=\linewidth]{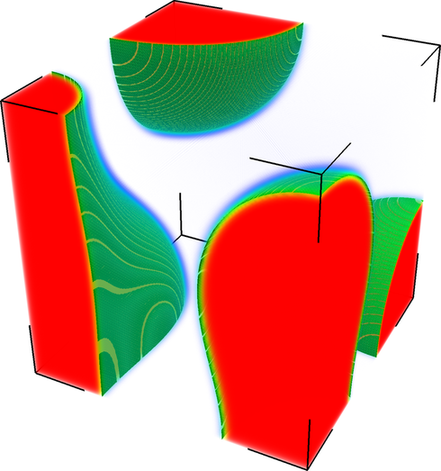}
\end{tabularx}
\caption{Simulation according to~\cref{sec:CHmixed:spinodalDecomposition} using after $2^4$, $2^7$, $2^{10}$, $2^{13}$, $2^{16}$~steps.  The colors indicate the phase~\emph{(blue\,/\,transparent \mbox{\normalfont for $c{<}0$}, red \mbox{\normalfont for $c{>}0$} and green \mbox{\normalfont for $c{=}0$})}. For~$p=1$, the displayed contour surface with value zero represents the center of the interface.}
\label{fig:sim:SD04pjbetai}
\end{figure}

\subsection{Spinodal decomposition---constant versus degenerate mobility}\label{sec:CHmixed:spinodalDecomposition}
Spinodal decomposition is a~mechanism of phase transformation through which an~initially homogeneous mixture decomposes into separate zones rich in separate components, for instance, due to a~spontaneous pressure drop of the system.
This process gives rise to the creation of interfaces and introduces additional gradient energy into the system, which is subsequently minimized together with the chemical energy (cf.~\eqref{eq:HelmholtzFreeEnergy}) during the evolution of the system from a~thermodynamically unstable state to a~stable one.
The time scale for the degenerate mobility case \mbox{($\beta=1$)} is slower than for the constant mobility case \mbox{($\beta=0$)} since the mass flux between bulk phases is suppressed.
While for~$\beta=0$, a~system at stationary state always consists of exactly two bulk phases, the stationary state for~$\beta=1$ may reveal multiple bulk phases (e.\,g.~non-coalescent droplets)---a~phenomenon that turns out to be controlled by the chosen absolute tolerance of the nonlinear numerical solver.
\paragraph{Scenario definition}
The motivation of this scenario is the quantification of the effects of constant and degenerate mobility in conjunction with an~approximation order~$p\in\{0,1\}$.  A~qualitatively similar scenario was used by~\cite{BarrettBloweyGarcke1999,Kim2007}.
We use a~grid with $100^3$~elements, i.\,e.~$h=1/100$, a~minimal interface width~$\kappa = h^2$, and $\tau=2\E{-3}$ for $p=0$, $\tau=1\E{-3}$ for $p=1$. 
The initial data, which is imposed for all investigated cases in order to render the comparisons on equal footing, is chosen element-wise constant with $c^0|_{E_k} = -0.4 \pm 0.05\round(\alpha_k)$ for every element~$E_k\in\setE_h$, where~$\alpha_k$ is random number uniformly distributed in~$[-1,1]$.
\paragraph{Results and discussion of stationary state}
Fig.~\cref{fig:sim:SD04pjbetai} shows snapshots of the order parameter at different times. We vary $p$ and $\beta$.
The paths and time-scales of evolution are clearly different between constant and degenerate mobility cases for a~given $p$.
In early stages of the simulation, the Newton solver will stop by the relative criterion, while in latter ones by the absolute criterion, cf.~\eqref{eq:Newton:stop}.
A~stationary state is reached as soon as the initial residual of a~time step satisfies the nonlinear stopping criterion.
For $\beta=0$ this is typically the case when the mixture has decomposed into two bulk phases and the surface is minimized.  
For $\beta=1$, the bulk diffusion is severely restricted by the degenerate mobility coefficient and the evolution is thus much slower, especially toward later stages of the simulation.  
Thus, it depends very much on the chosen nonlinear absolute tolerance at which state of the spinodal decomposition a~stationary state is reported.  Due to the fact that there is always numerical diffusion in the system (FV \mbox{($p=0$)} is diffusive and DG \mbox{($p>0$)} contains diffusion due to the penalty term), we claim that also for the case~$\beta=1$ a~state with two bulk phases can be reached by choosing this tolerance small enough (if permitted by the condition numbers of the linear systems per nonlinear iteration).
\paragraph{Choice of time step size}
The numerical scheme is by design sufficiently robust, which  permits the use of an arbitrarily large time step size. Having stated that, a~too large time step size gives rise to larger time~discretization errors and higher condition numbers in the system matrices.  
The pattern of the phases is sensitive to the time step size: with smaller time steps, phases tend to evolve in a~more droplet-like structure while larger time steps will create a~more inter-connected one. 
The likely reason is the damping character of the time discretization, which for larger time steps smears out the solution such that coarser phases evolve.
%
\begin{figure}[t!]%
\begin{tabular}{@{}lr@{}}
\begin{tabularx}{0.5\linewidth}{@{}LL@{}}
\includegraphics[width=\linewidth]{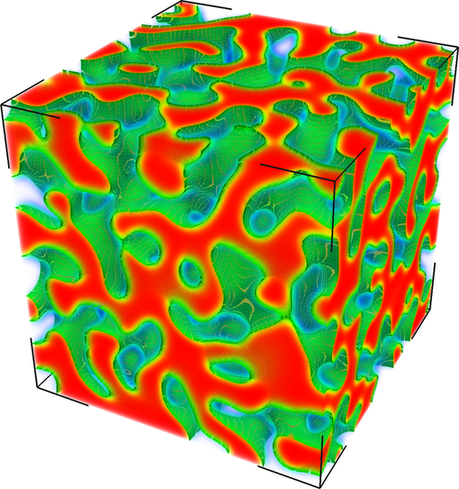} &
\includegraphics[width=\linewidth]{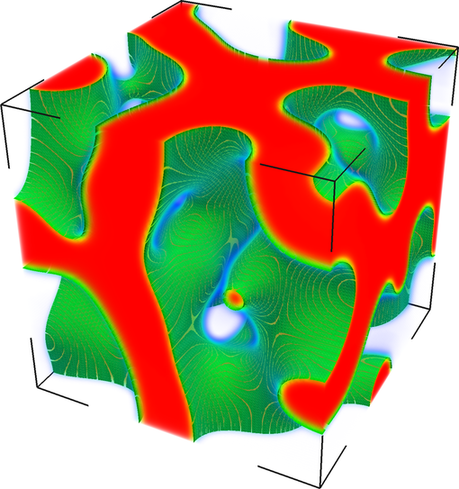} \\
\includegraphics[width=\linewidth]{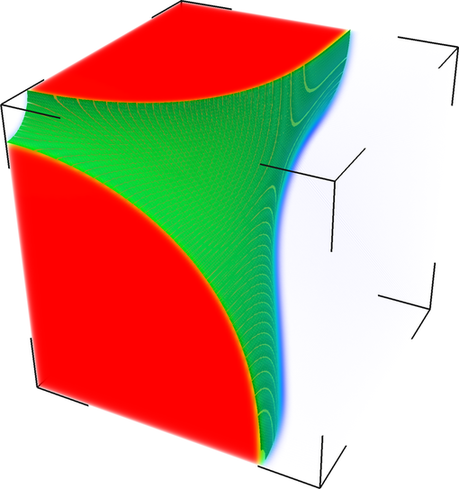} &
\includegraphics[width=\linewidth]{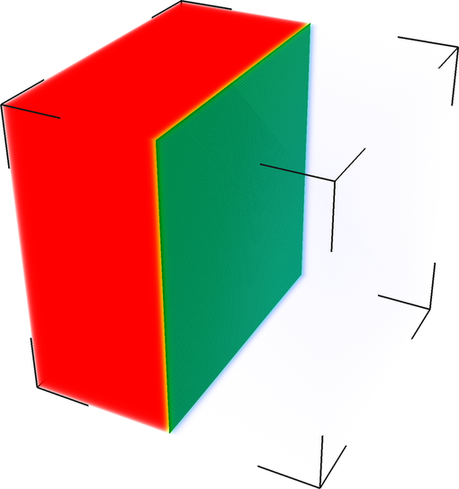}
\end{tabularx}
&
\begin{tabularx}{0.45\linewidth}{@{}R@{}}
\includegraphics[width=\linewidth]{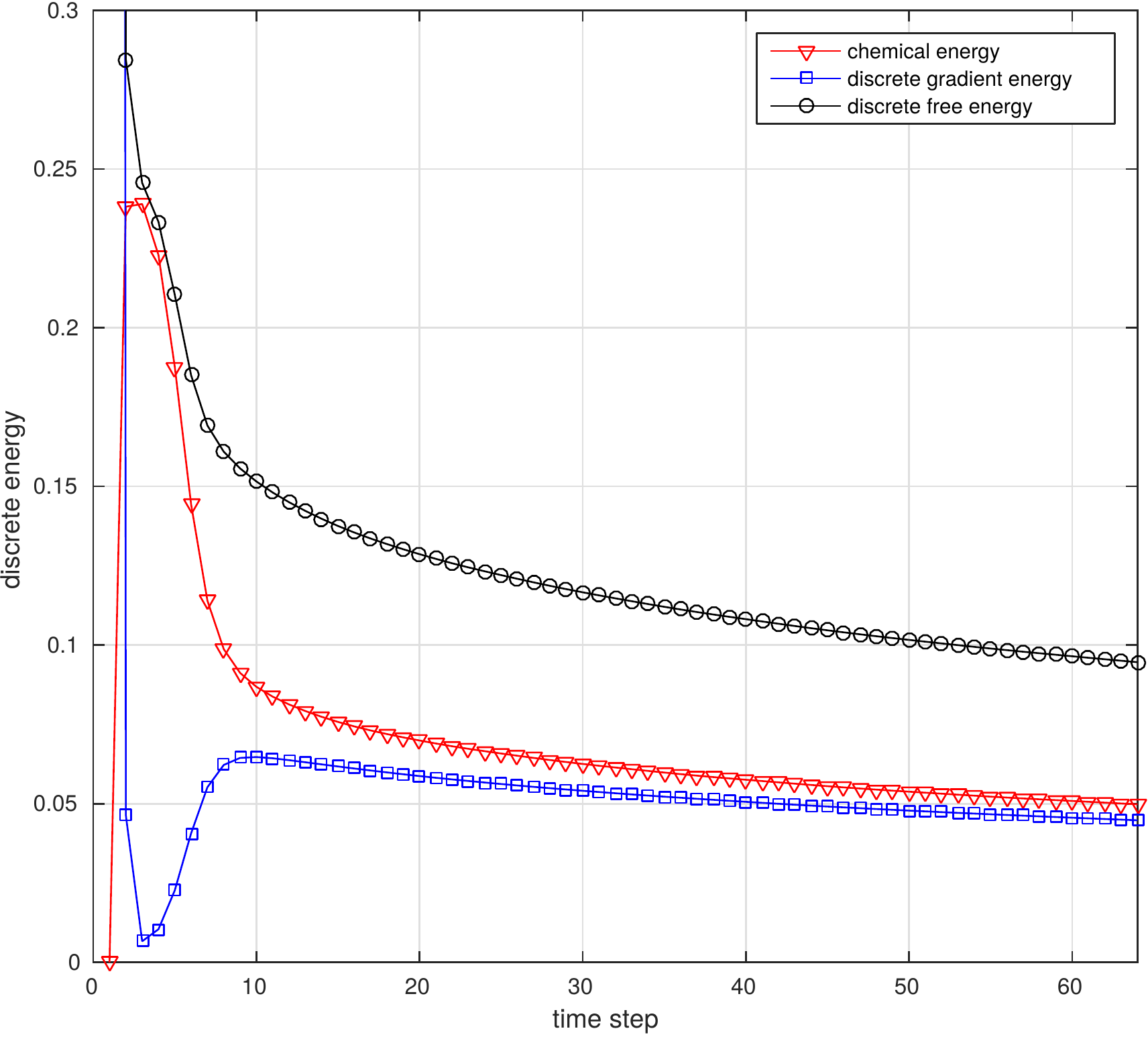}
\end{tabularx}
\end{tabular}
\caption{Scenario of~\cref{sec:CHmixed:spinodalDecomposition}: 
solution profile after $2^4$, $2^9$, $2^{14}$, $2^{16}$~steps~\emph{(left, same visualization setting as in~\cref{fig:sim:SD04pjbetai})}. Chemical energy, discrete gradient energy~(cf.~\eqref{eq:HelmholtzFreeEnergyDiscrete}), and discrete free energy for the first 128~steps~\emph{(right)}.   }
\label{fig:sim:SD00p1beta0:energyDissipation}
\end{figure}
\paragraph{Spinodal pattern and energy dissipation}
As opposed to the droplet-like pattern evolving from the initial data above, the phase structure during evolution is highly interconnected and \enquote{fat spaghetti}-like~\cite[][]{Elliott1989CH} when using random data with expected value zero instead of~$-0.4$, cf.~\cref{fig:sim:SD00p1beta0:energyDissipation}.  Here we use~$c^0|_{E_k} = \round(\alpha_k)$, i.\,e., for every element, the order parameter at initial time is either~$-1$ or~$1$.  For this setting, \cref{fig:sim:SD00p1beta0:energyDissipation}~shows the chemical energy~$\int_{\Omega_h}\Phi(c_h^n)$, the discrete gradient energy~\raisebox{0pt}[0pt][0pt]{$a^\mathrm{diff}(1,c_h^n,c_h^n)$}, and the discrete free energy~$F_h(c_h^n)$ against time: The chemical energy is zero at initial state since~$\forall \vec{x}\in\Omega_h,~c_h^0(\vec{x})\in\{-1,1\}$.  Due to the jumps of~$c_h^0$ across interior faces~$e\in\Gamma^\mathrm{int}$, the discrete gradient energy takes a~maximum at initial state (note that the continuous version of the gradient energy as given in~\eqref{eq:HelmholtzFreeEnergy} is zero here).
Even though neither the chemical nor the discrete gradient energy is monotonic, the discrete free energy~$F_h$---as a~sum of the two---is monotonically decreasing, cf.~\eqref{eq:HelmholtzFreeEnergyDiscrete}.
In this simulation, the interface has no curvature at stationary state and satisfies the maximum principle.  An~interesting finding here is that for a non-curved diffuse interface, the concentration profile does not exhibit any bulk shift that violated the maximum principle, a phenomenon which will be analyzed in detail in the following section.
\begin{table}[ht!]\centering\footnotesize
\begin{tabularx}{\linewidth}{@{}cc@{\quad}|@{\qquad}C@{}C@{}CC@{}C@{}CC@{}C@{}CC@{}C@{}C@{}}
\toprule
$p$ & $\beta$ & {} & $\tau$ & {} & {} & $h$ & {} & {} & $L$ & {} & {} & $\kappa$
\\\midrule
$0$ & $\{0,1\}$ & $1/2000$ & $\boldsymbol{1/1000}$ &  $1/500$ & $1/80$& $\boldsymbol{1/40}$& $1/20$ & $1/4$& $\boldsymbol{1/2}$& $3/4$ & $1/1600$& $\boldsymbol{1/800}$& $1/400$
\\
$1$ & $\{0,1\}$ & $1/4000$ & $\boldsymbol{1/2000}$& $1/1000$ & $1/80$& $\boldsymbol{1/40}$& $1/20$ & $1/4$& $\boldsymbol{1/2}$& $3/4$ & $1/1600$& $\boldsymbol{1/800}$& $1/400$
\\
$2$ & $\{0,1\}$ & $1/8000$ & $\boldsymbol{1/4000}$& $1/2000$ & $1/80$& $\boldsymbol{1/40}$& $1/20$ & $1/4$& $\boldsymbol{1/2}$& $3/4$ & $1/1600$& $\boldsymbol{1/800}$& $1/400$
\\
\midrule
$0$ & $0$ & 0.0283 & \textbf{0.0283} &   0.0283 &    0.0227 & \textbf{0.0283} & 0.0329 & --1.8285 & \textbf{0.0283} & 0.0175 & 0.0360 & \textbf{0.0283} & 0.0283\\
$0$ & $1$ & 0.0015 & \textbf{0.0018} &   0.0022 &  --0.0083 & \textbf{0.0018} & 0.0155 & --0.0426 & \textbf{0.0018} & 0.0006 & 0.0002 & \textbf{0.0018} & 0.0015\\
$1$ & $0$ & 0.0287 & \textbf{0.0288} &   0.0288 &    0.0214 & \textbf{0.0288} & 0.0313 & --1.8285 & \textbf{0.0288} & 0.0179 & 0.0376 & \textbf{0.0288} & 0.0287\\
$1$ & $1$ & 0.0085 & \textbf{0.0090} &   0.0099 &  --0.0081 & \textbf{0.0090} & 0.0313 & --0.0218 & \textbf{0.0090} & 0.0069 & 0.0163 & \textbf{0.0090} & 0.0085\\
$2$ & $0$ & 0.0285 & \textbf{0.0285} &   0.0285 &    0.0204 & \textbf{0.0285} & 0.0287 & --1.8285 & \textbf{0.0285} & 0.0176 & 0.0368 & \textbf{0.0285} & 0.0285\\
$2$ & $1$ & 0.0003 & \textbf{0.0004} &   0.0005 &  --0.0109 & \textbf{0.0004} & 0.0093 & --0.0265 & \textbf{0.0004} & 0.0000 & 0.0007 & \textbf{0.0004} & 0.0003\\
\midrule 
$0$ & $0$ & 0.0321 & \textbf{0.0318} &   0.0316 &    0.0404 & \textbf{0.0318} & 0.0372 &   0.1662 & \textbf{0.0318} & 0.0188 & 0.0452 & \textbf{0.0318} & 0.0321\\
$0$ & $1$ & 0.0400 & \textbf{0.0392} &   0.0379 &    0.0472 & \textbf{0.0392} & 0.0372 &   0.0765 & \textbf{0.0392} & 0.0201 & 0.0476 & \textbf{0.0392} & 0.0400\\
$1$ & $0$ & 0.0327 & \textbf{0.0326} &   0.0323 &    0.0419 & \textbf{0.0326} & 0.0315 &   0.1662 & \textbf{0.0326} & 0.0192 & 0.0453 & \textbf{0.0326} & 0.0327\\
$1$ & $1$ & 0.0404 & \textbf{0.0400} &   0.0394 &    0.0464 & \textbf{0.0400} & 0.0316 &   0.0784 & \textbf{0.0400} & 0.0208 & 0.0482 & \textbf{0.0400} & 0.0404\\
$2$ & $0$ & 0.0326 & \textbf{0.0325} &   0.0323 &    0.0425 & \textbf{0.0325} & 0.0317 &   0.1662 & \textbf{0.0325} & 0.0190 & 0.0451 & \textbf{0.0325} & 0.0326\\
$2$ & $1$ & 0.0407 & \textbf{0.0406} &   0.0402 &    0.0467 & \textbf{0.0406} & 0.0319 &   0.0746 & \textbf{0.0406} & 0.0207 & 0.0488 & \textbf{0.0406} & 0.0407\\
\bottomrule
\end{tabularx}
\caption{Parameters of the scenario of~\cref{sec:bulkShift}~\emph{(top)}, the deviation to~$1$~\emph{(middle)}, and~$-1$~\emph{(bottom)} of the order parameter~$c_h$ at stationary state.  The default cases are bold, variations are performed in one parameter at a~time only.}
\label{tab:bulkShiftData}
\end{table}
%
\subsection{Investigation of the bulk shift---a~sensitivity study}\label{sec:bulkShift}
\begin{figure}[ht!]\centering\footnotesize
\begin{tabularx}{\linewidth}{@{}LL|LL@{}}
\includegraphics[width=\linewidth]{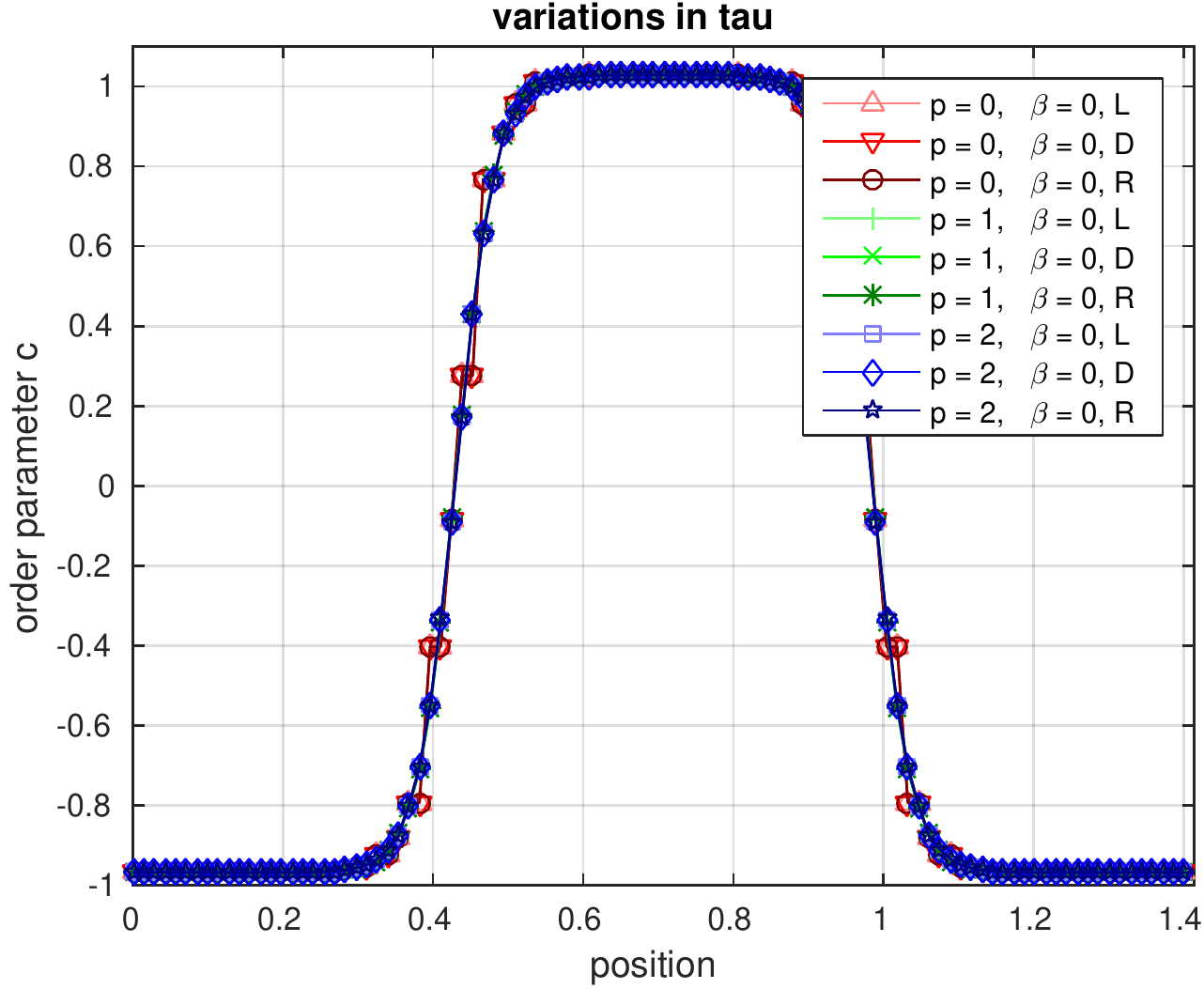} & 
\includegraphics[width=\linewidth]{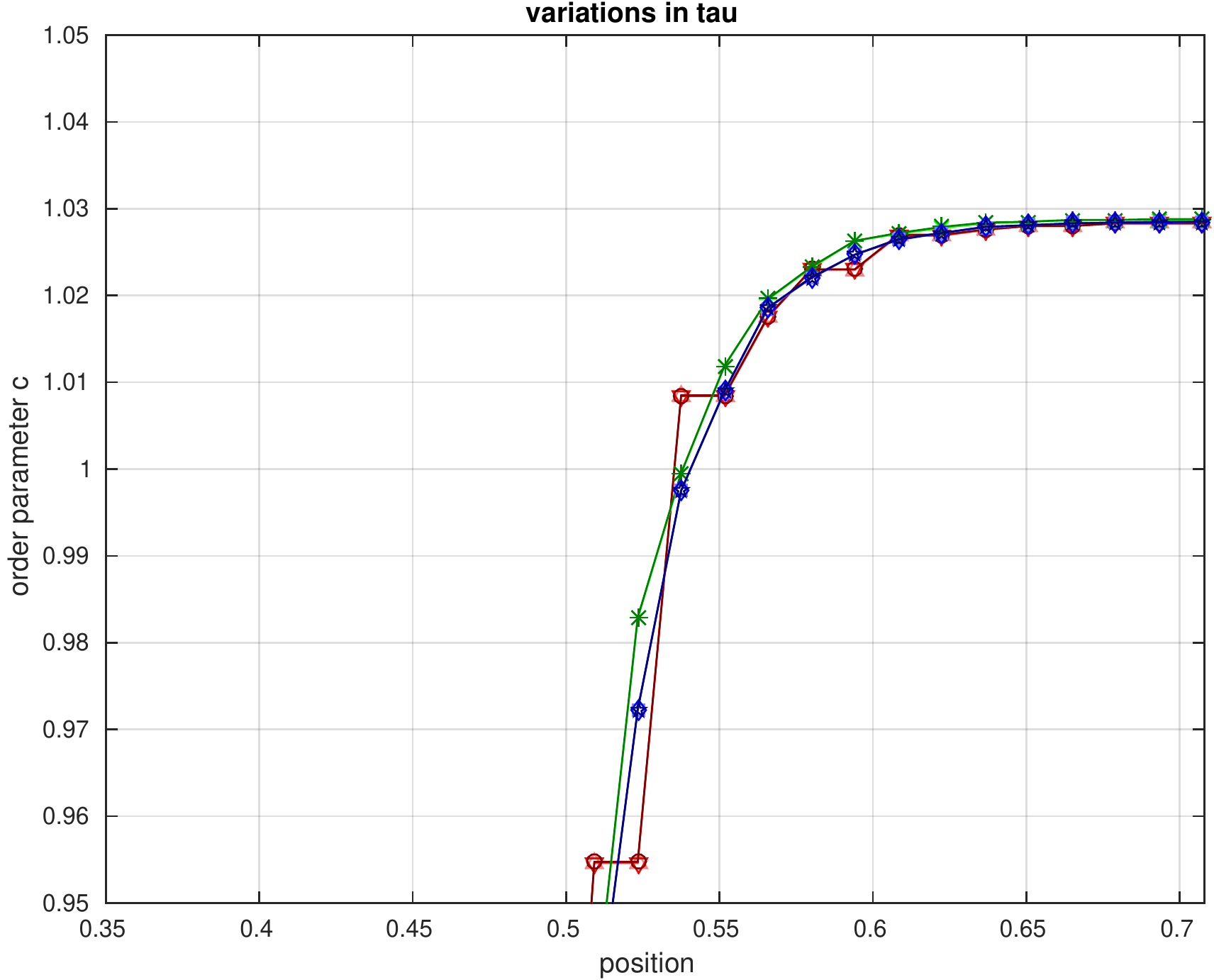} & 
\includegraphics[width=\linewidth]{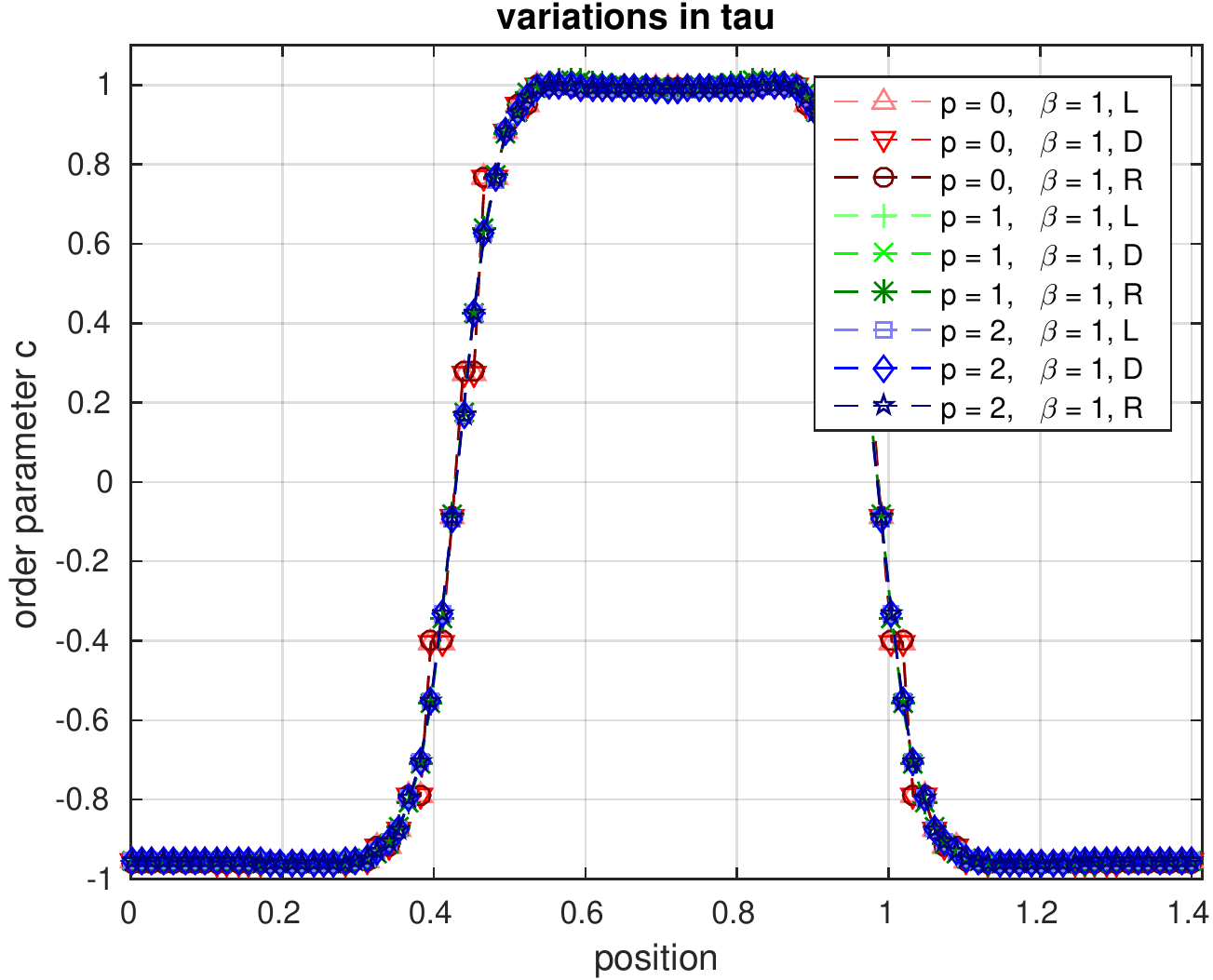} & 
\includegraphics[width=\linewidth]{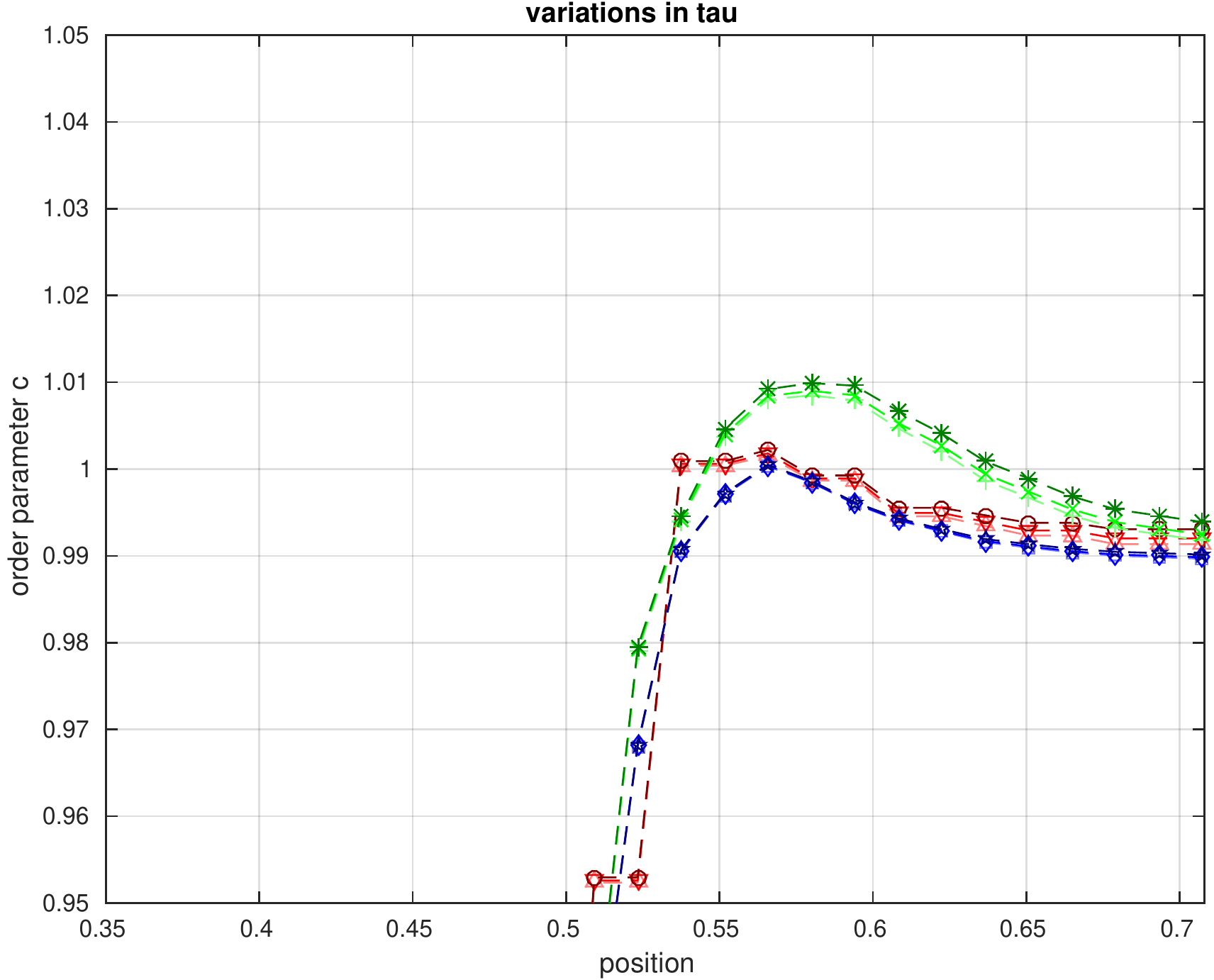} 
\\
\includegraphics[width=\linewidth]{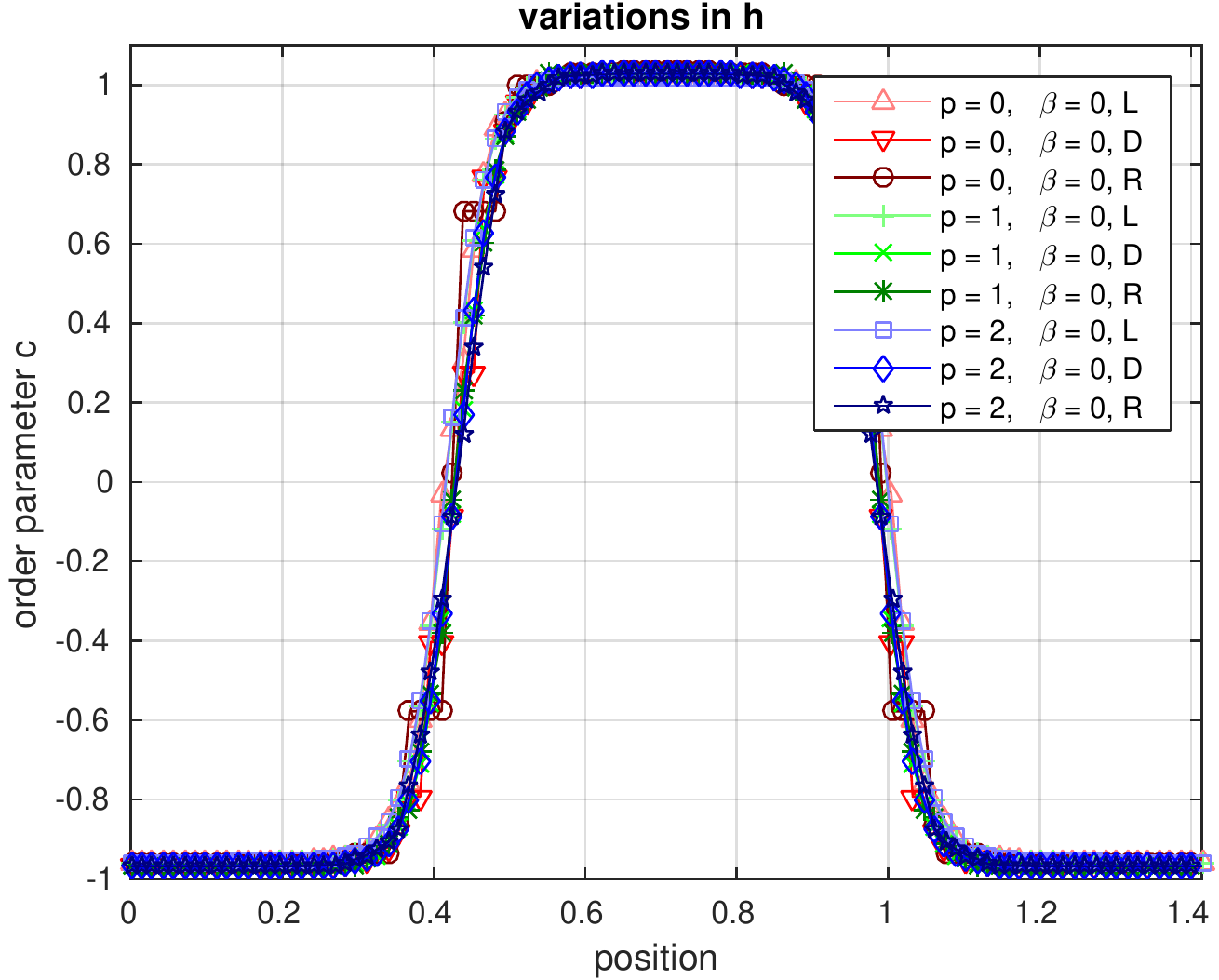} & 
\includegraphics[width=\linewidth]{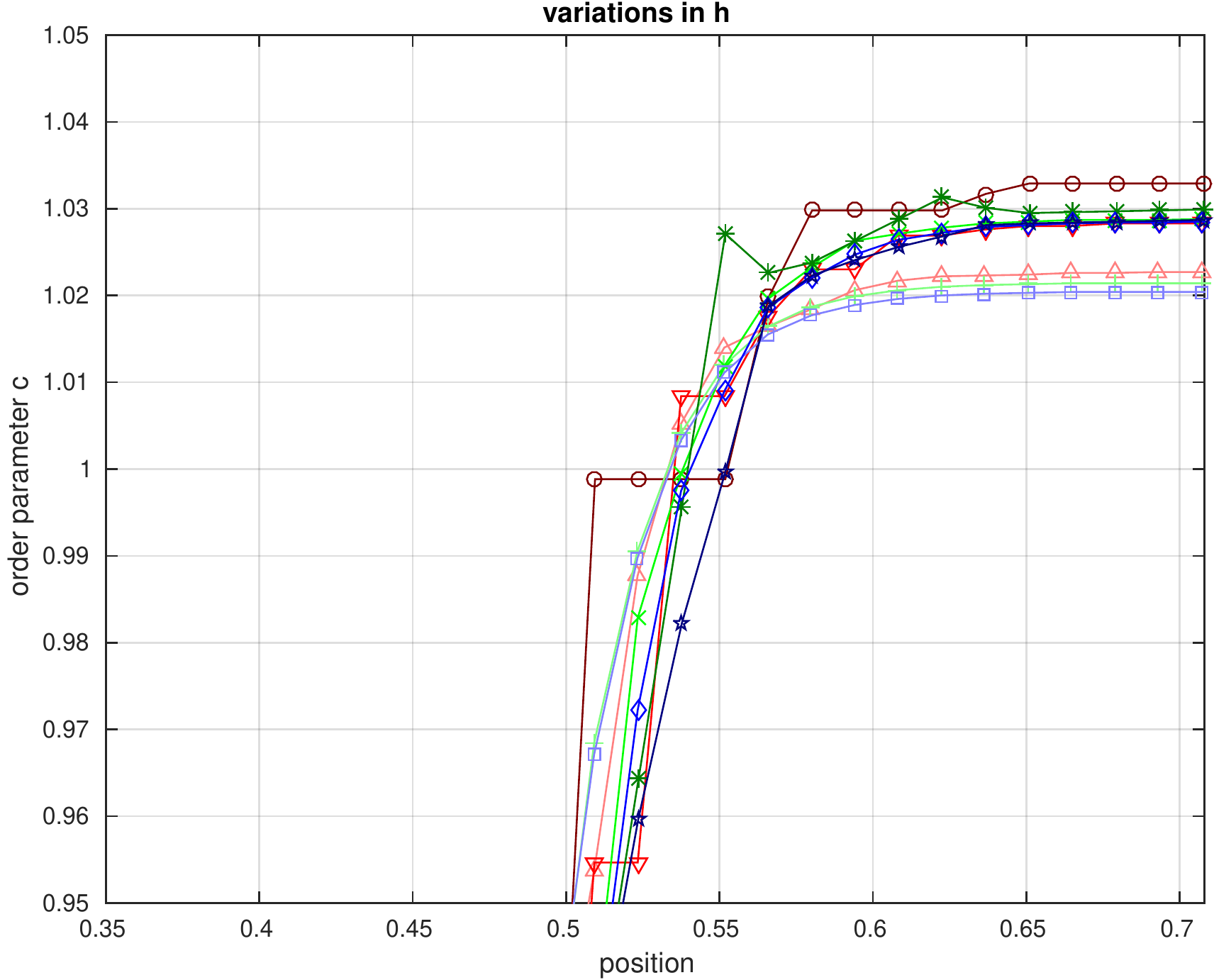} & 
\includegraphics[width=\linewidth]{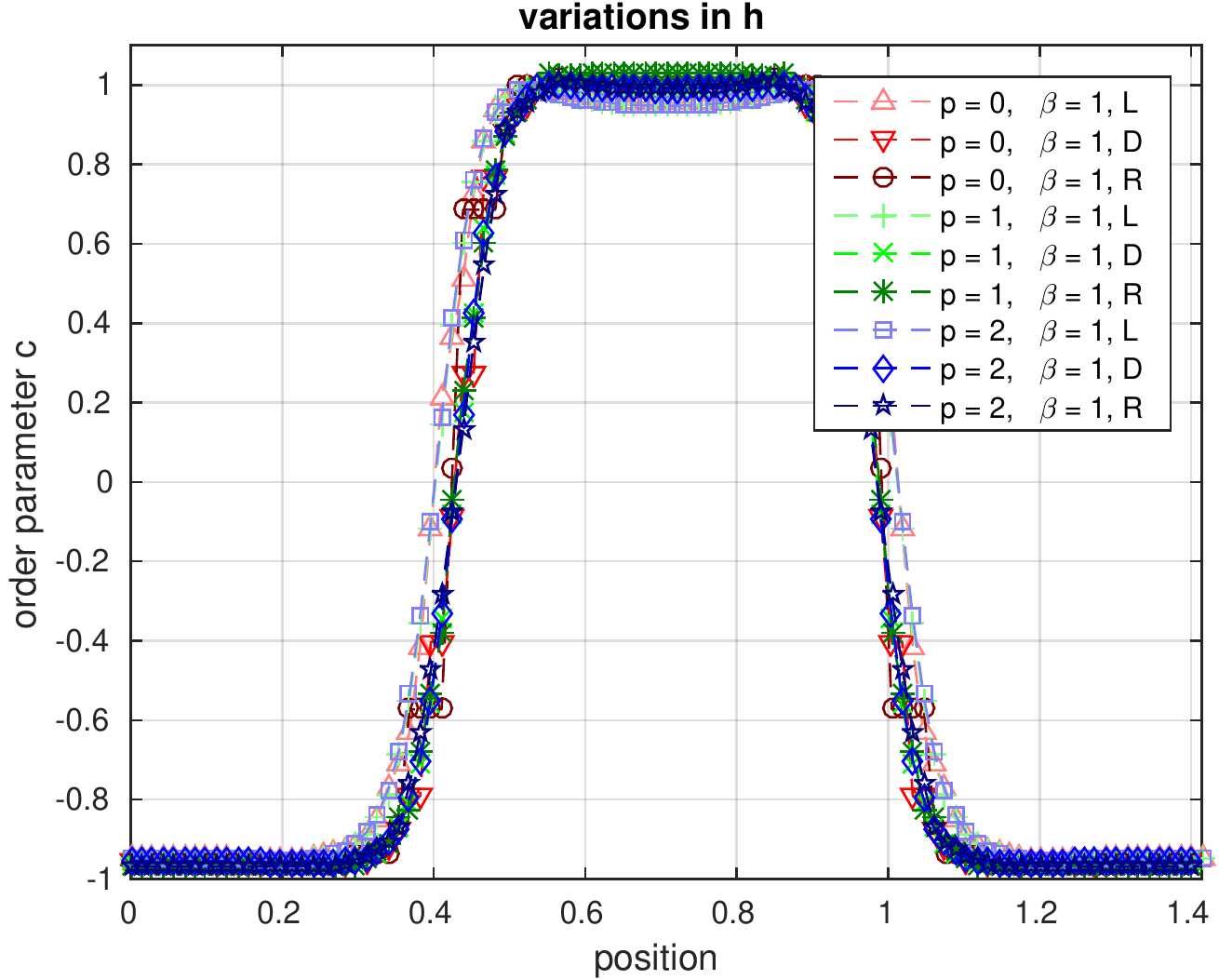} & 
\includegraphics[width=\linewidth]{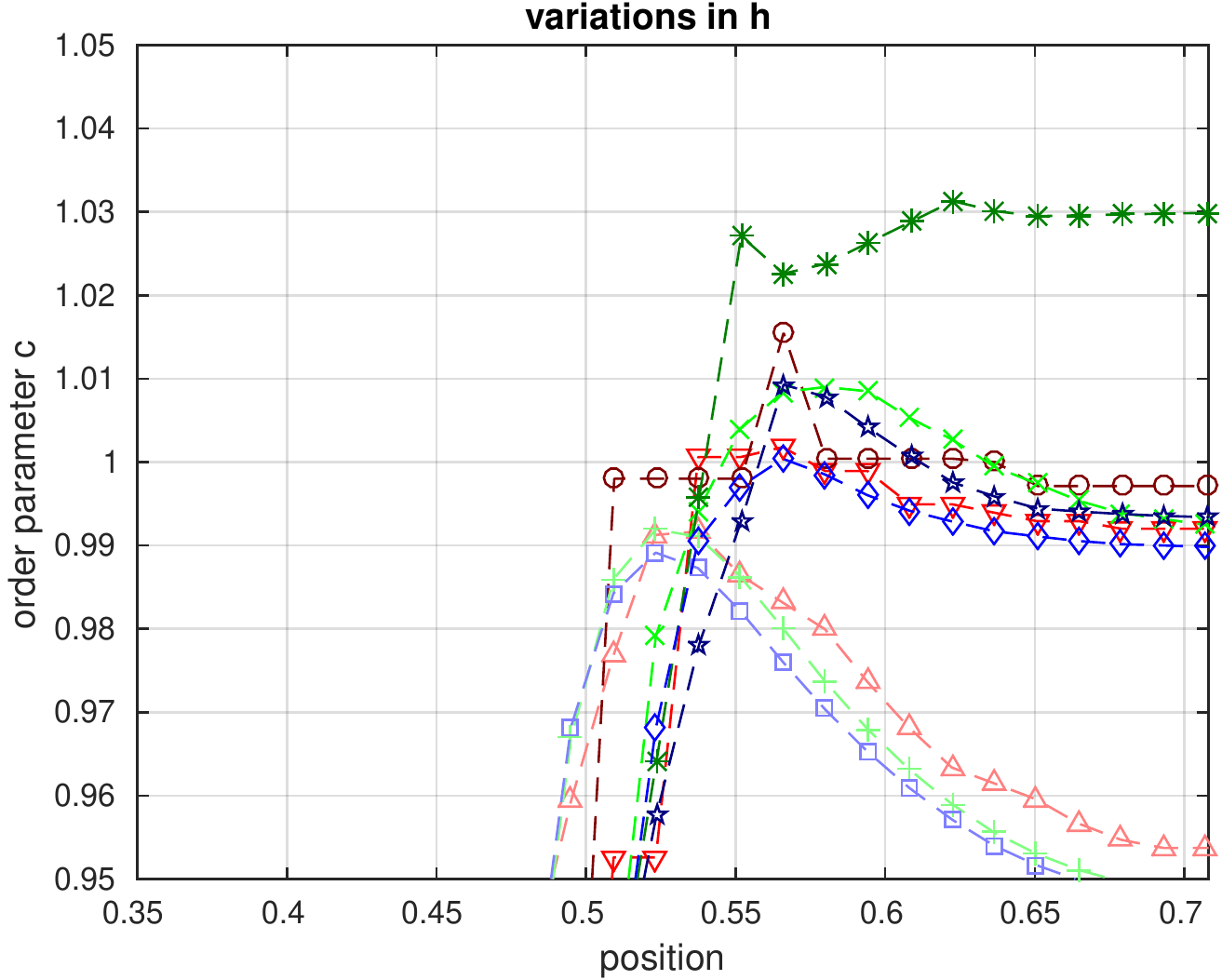} 
\\
\includegraphics[width=\linewidth]{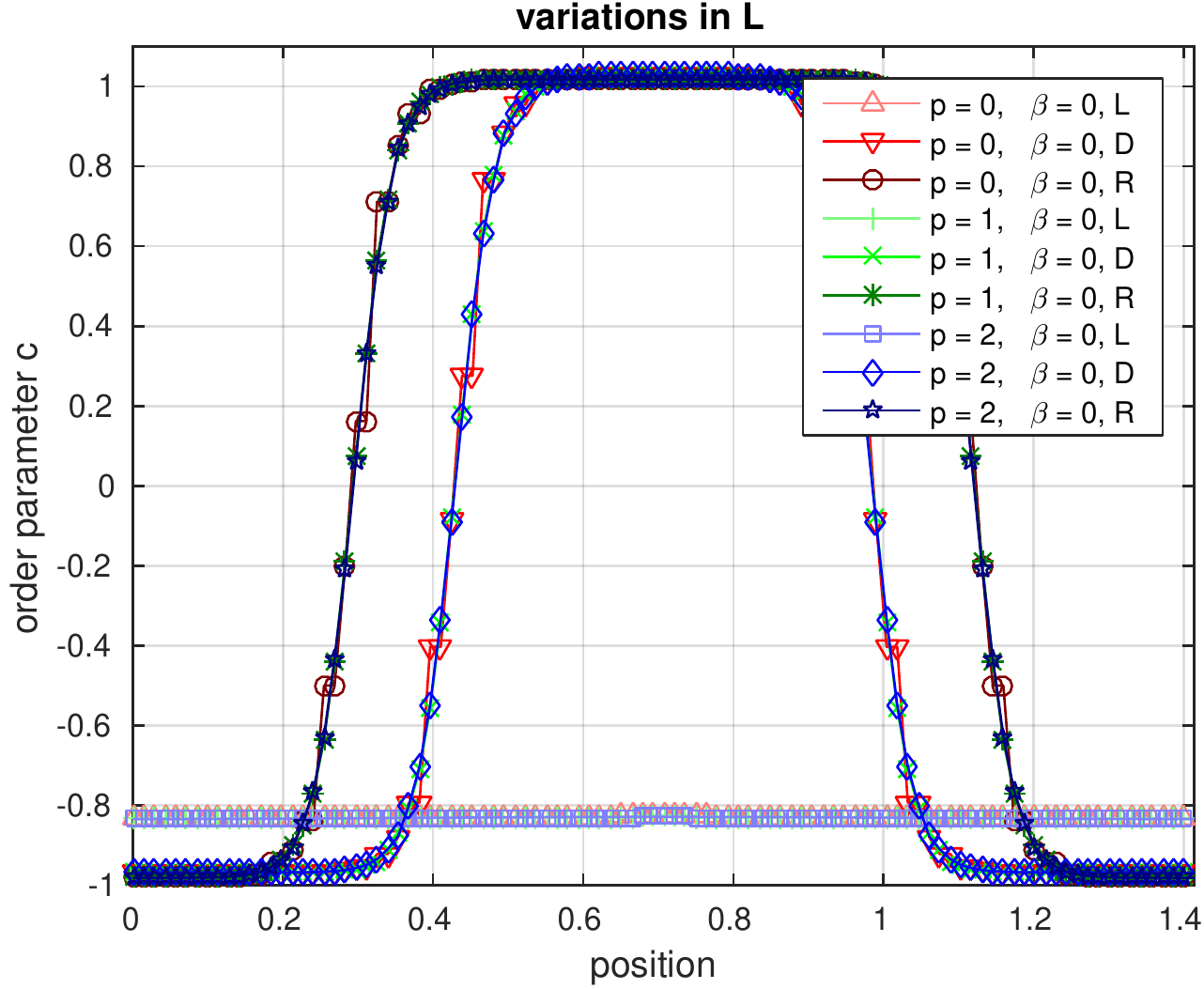} & 
\includegraphics[width=\linewidth]{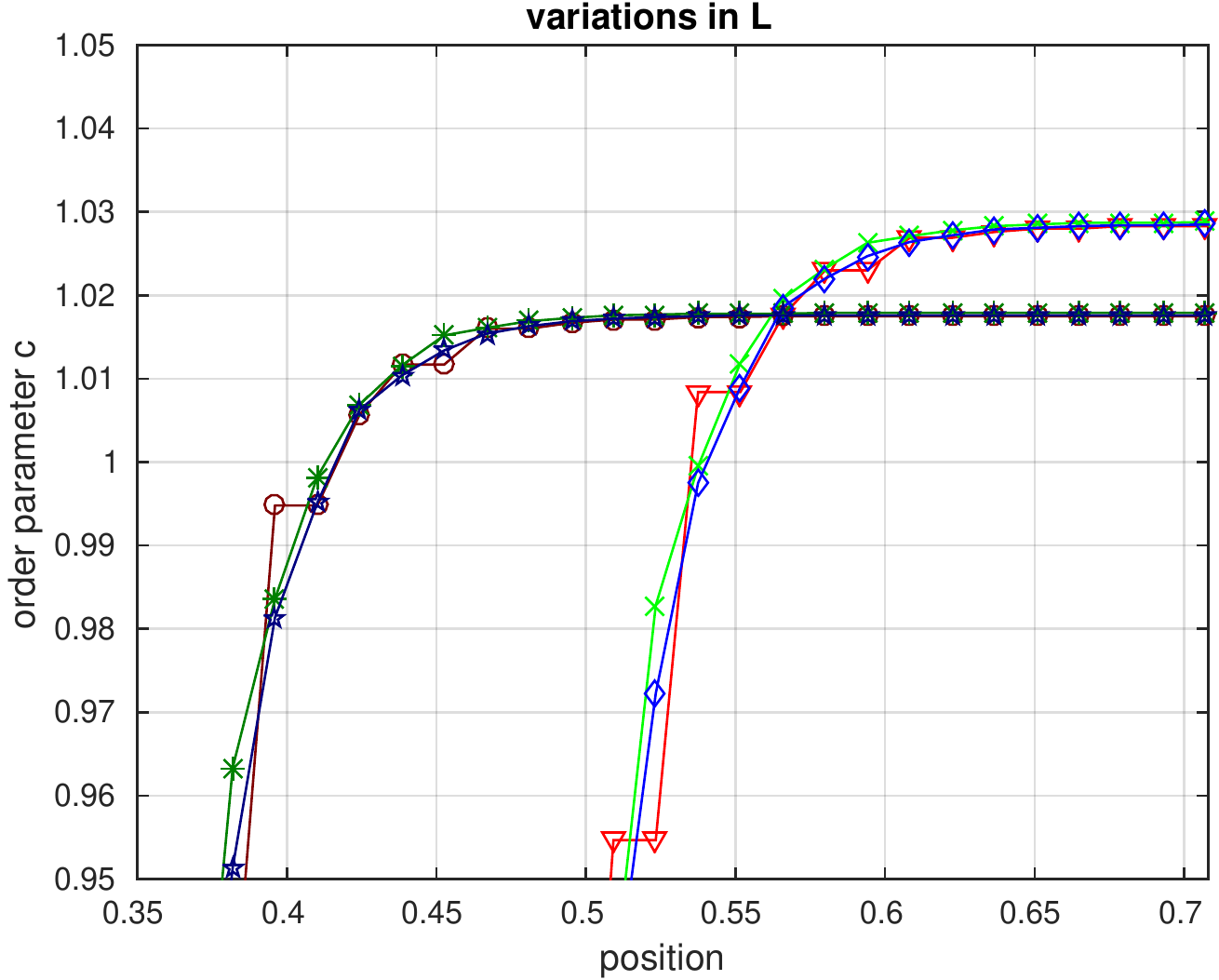} & 
\includegraphics[width=\linewidth]{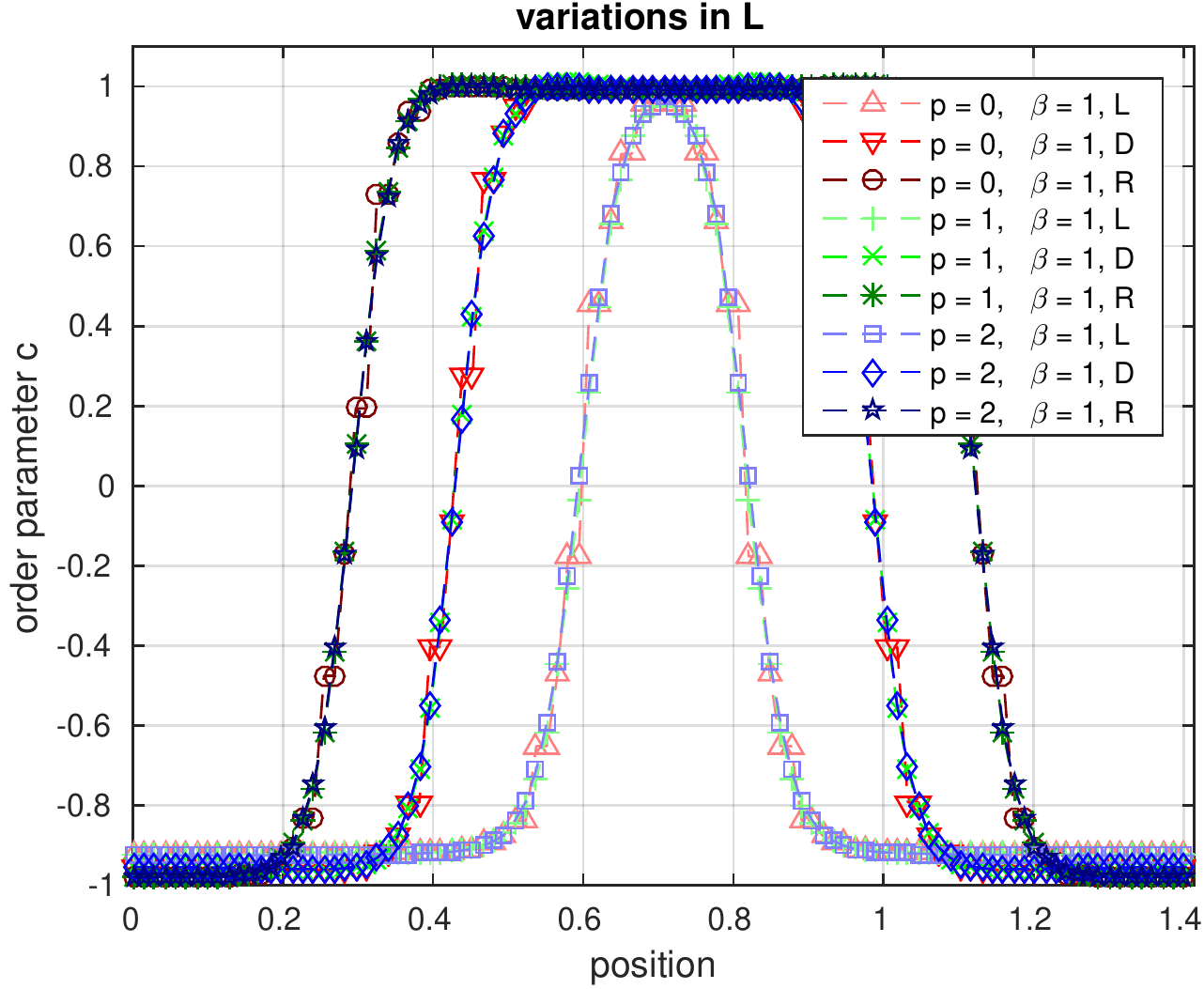} & 
\includegraphics[width=\linewidth]{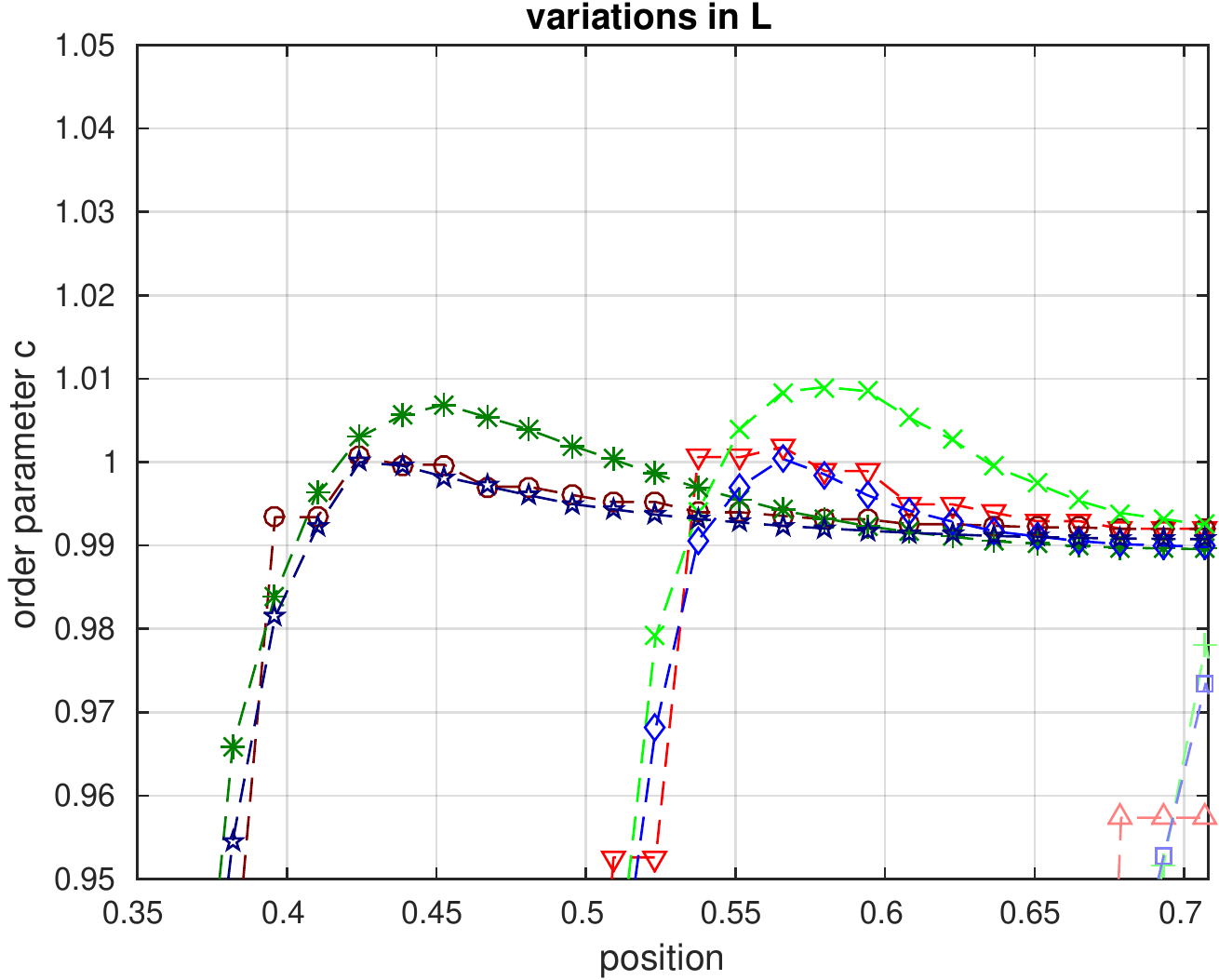} 
\\
\includegraphics[width=\linewidth]{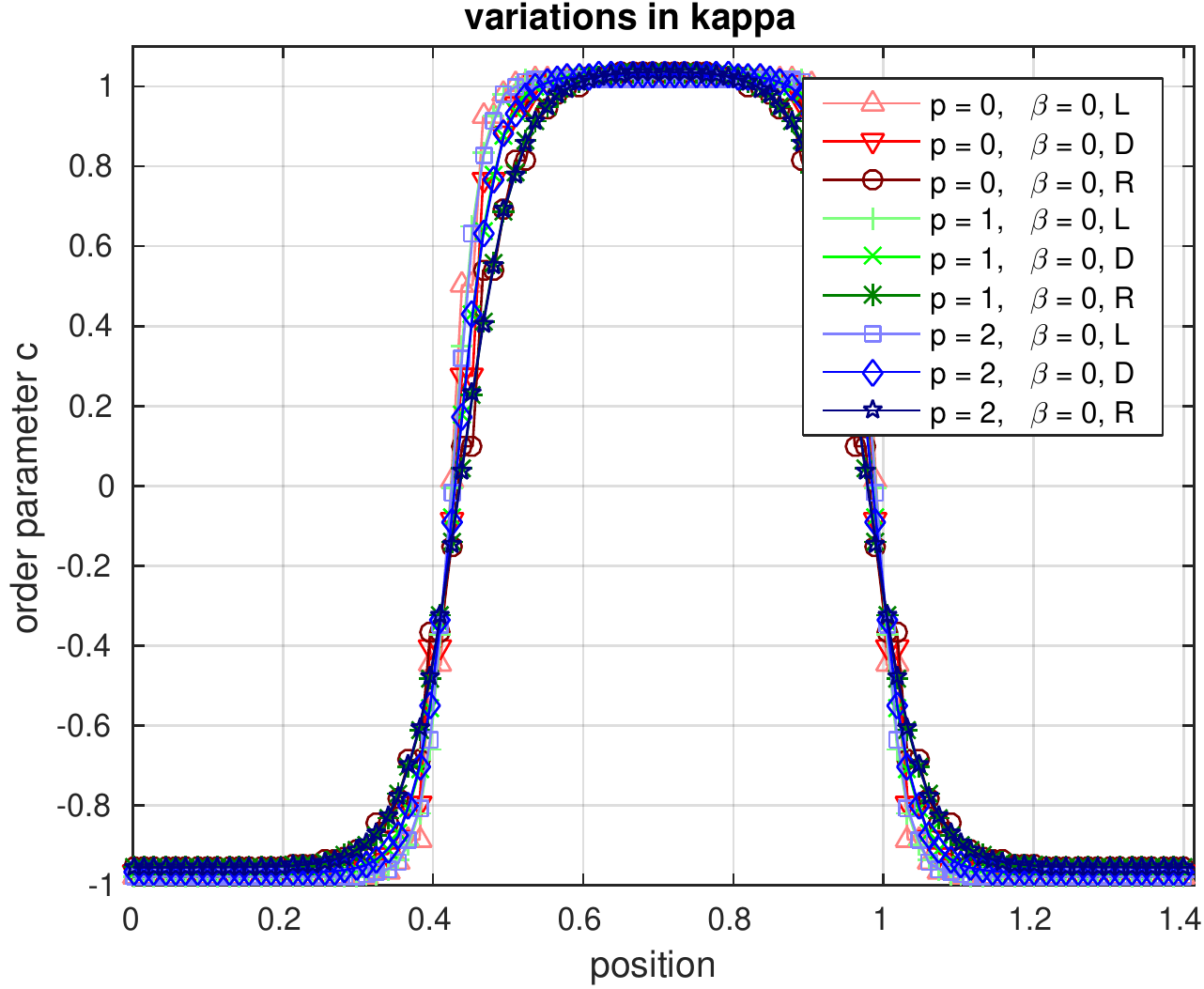} & 
\includegraphics[width=\linewidth]{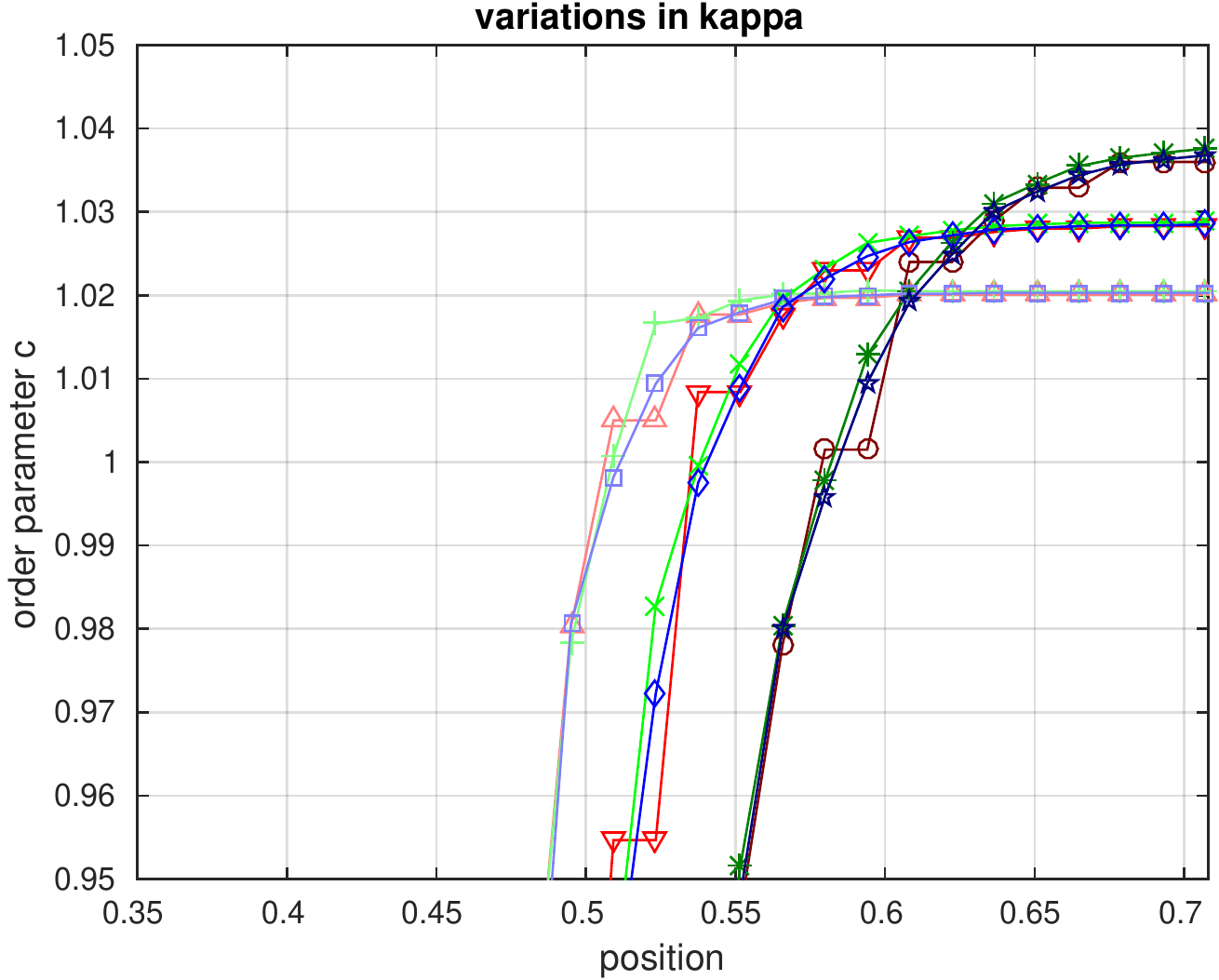} & 
\includegraphics[width=\linewidth]{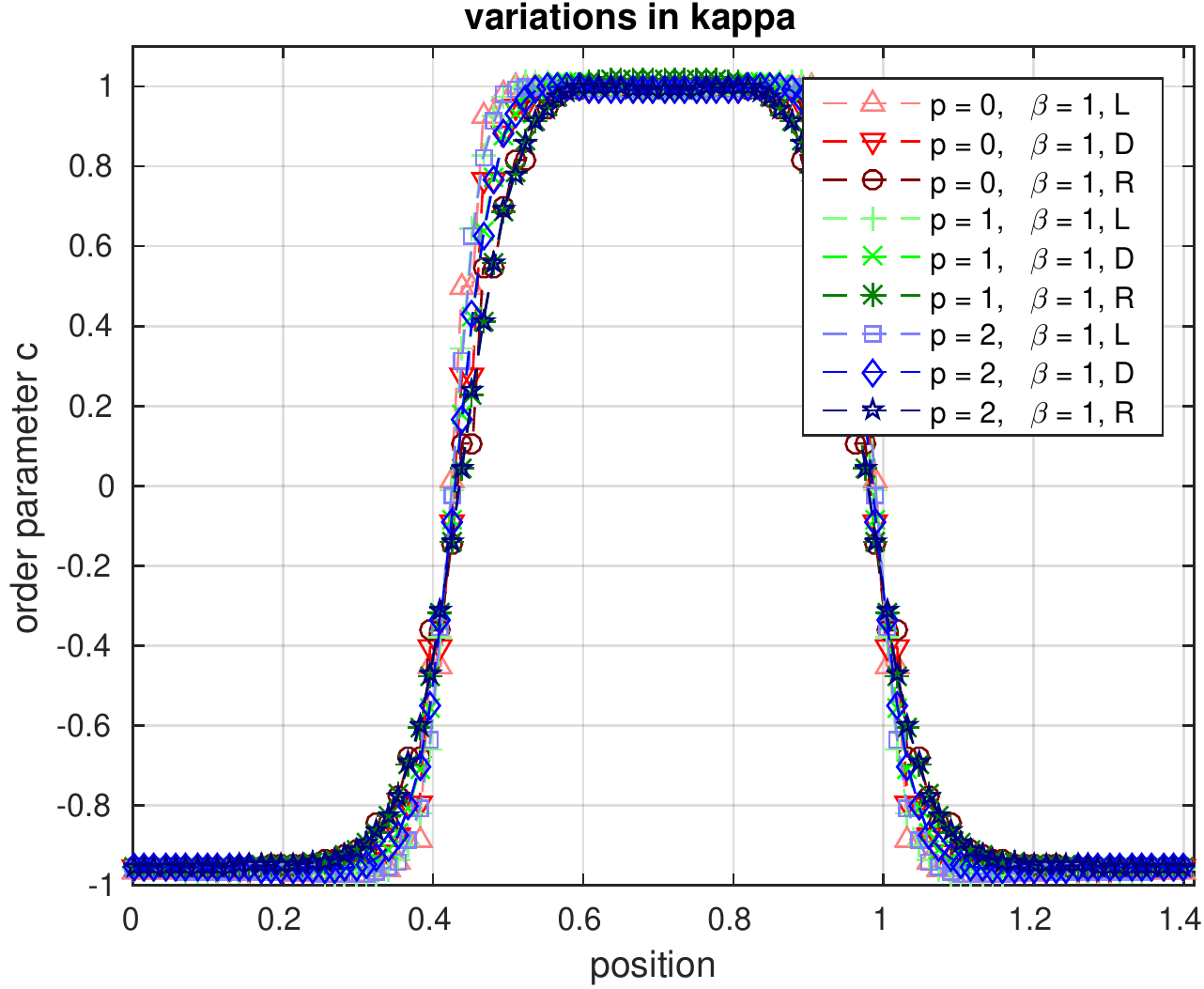} & 
\includegraphics[width=\linewidth]{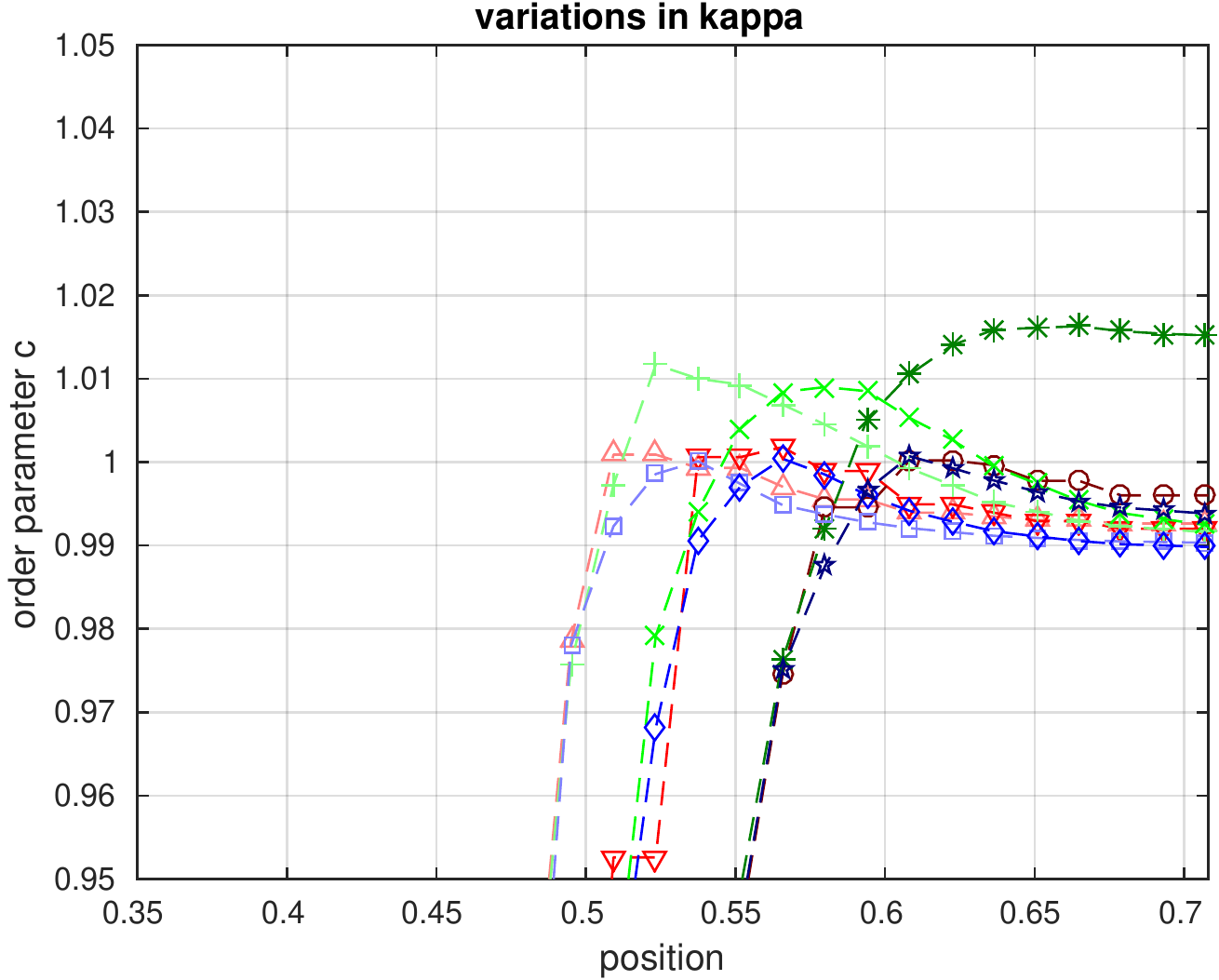} 
\end{tabularx}
\caption{Stationary states of the order parameter across a~diagonal line through the domain for the scenario of~\cref{sec:bulkShift} for $p=0$~\emph{(red)}, $p=1$~\emph{(green)}, $p=2$~\emph{(blue)}, and $\beta=0$~\emph{(solid)}, $\beta=1$~\emph{(dashed)}.  The default case \emph{(D)} is normal color, the smaller and larger variation is lighter and darker, respectively, cf.~\cref{tab:bulkShiftData} \emph{((L), (R) correspond to column on left and right of the default case, respectively)}.  The zoomed-in plots do not preserve the aspect-ratio.}
\label{fig:bulkShiftPlots}
\end{figure}
\begin{figure}[ht!]\centering\footnotesize
\begin{tabularx}{\linewidth}{@{}LL|LL@{}}
\includegraphics[width=\linewidth]{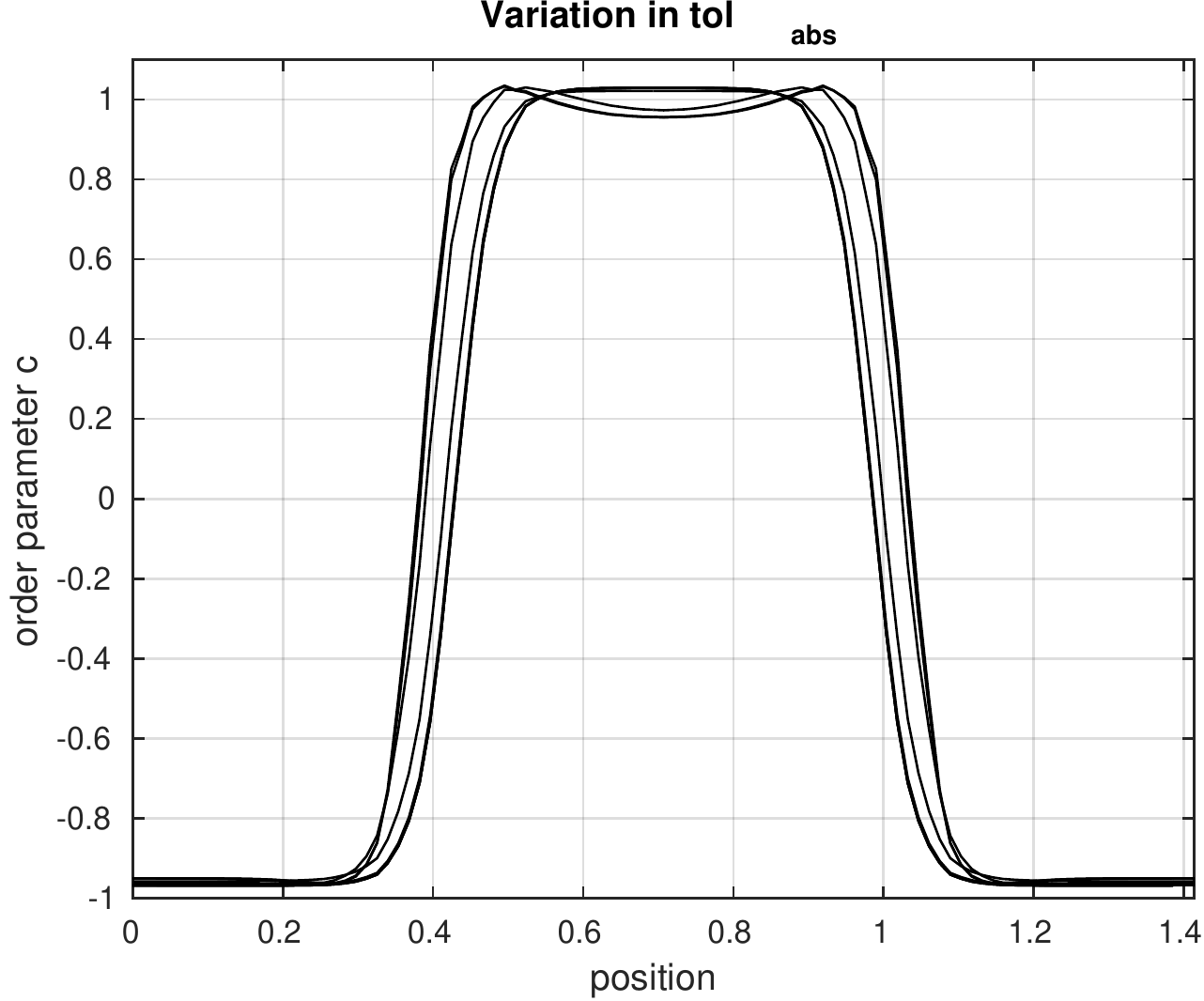} & 
\includegraphics[width=\linewidth]{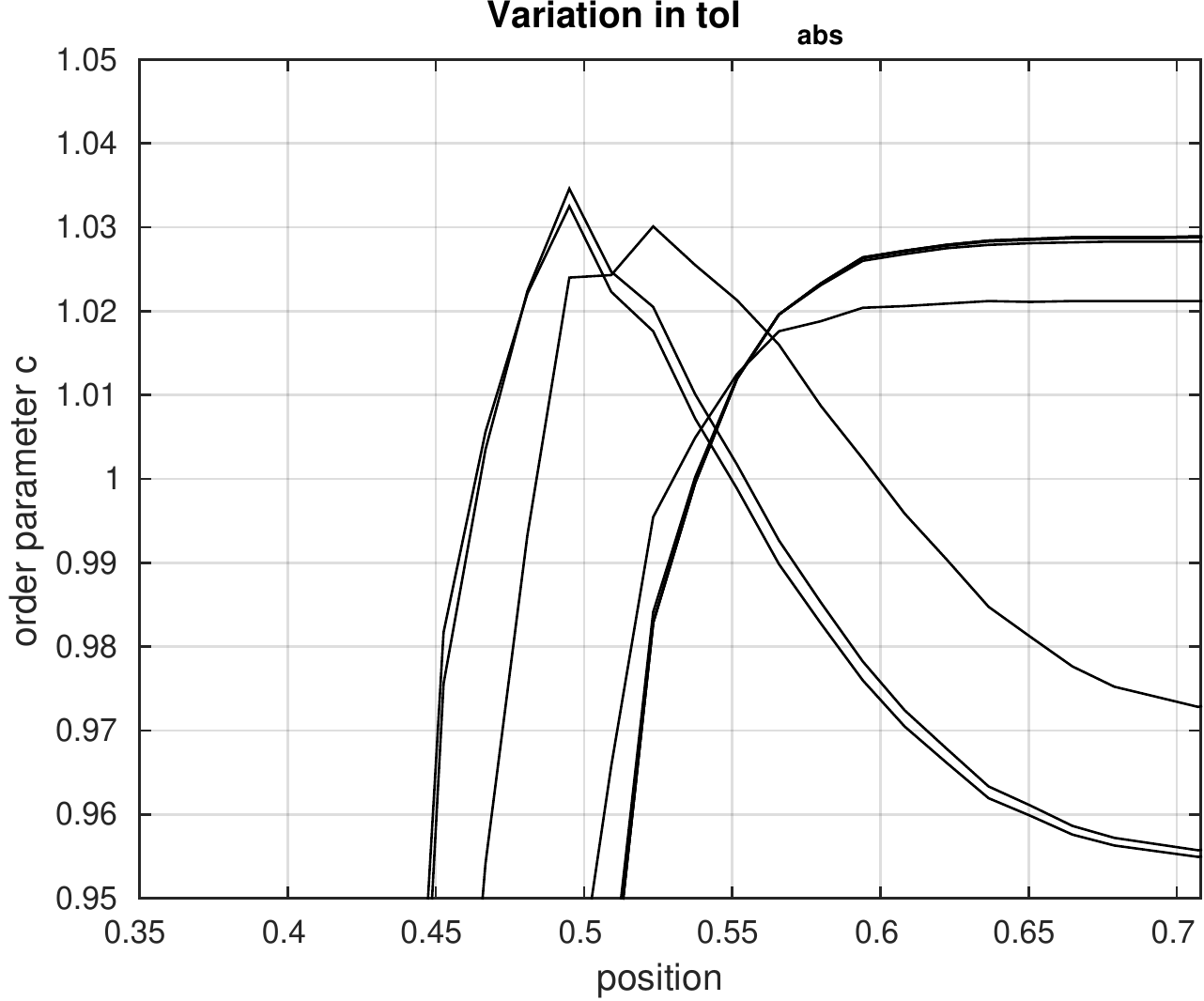} & 
\includegraphics[width=\linewidth]{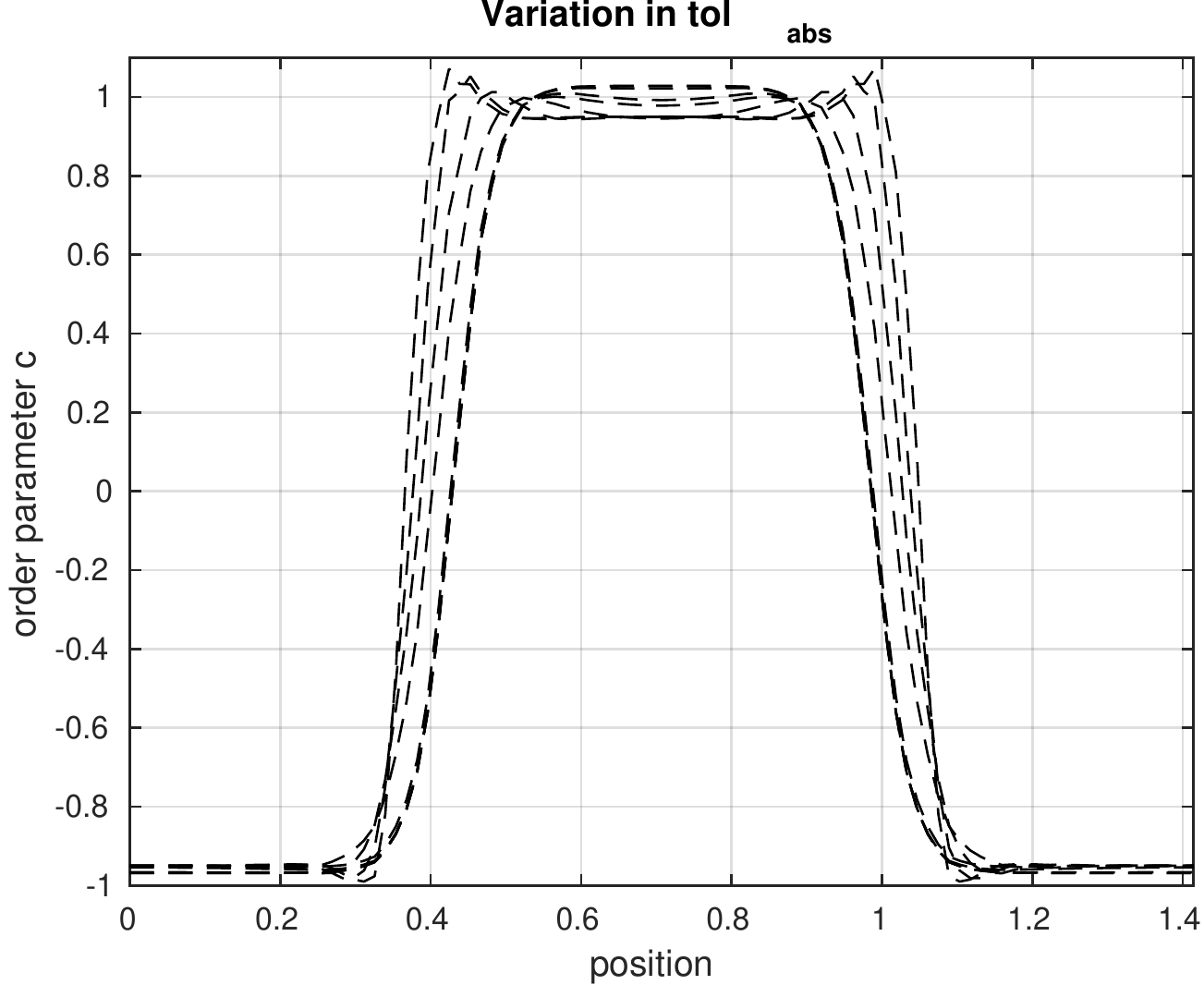} & 
\includegraphics[width=\linewidth]{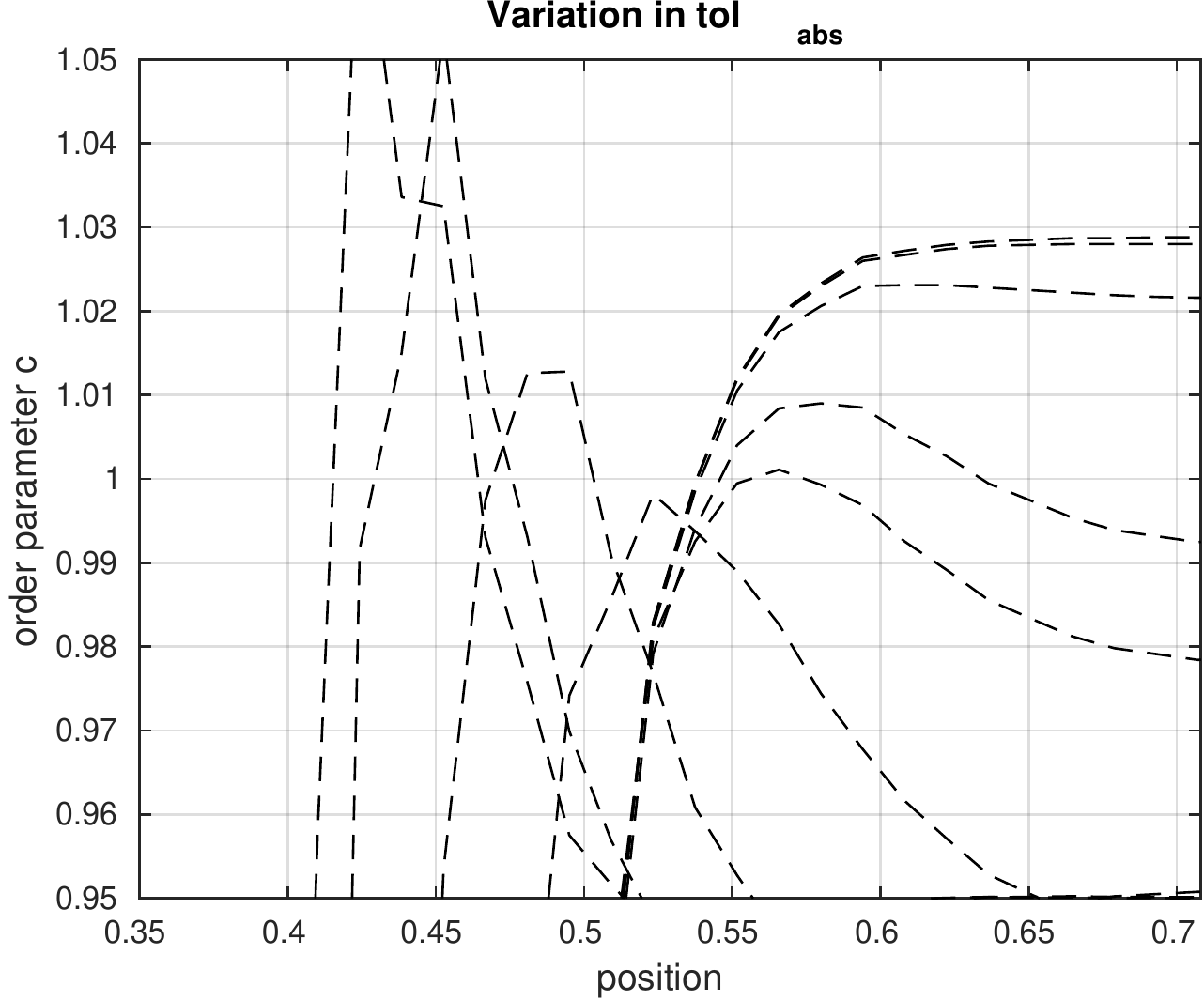} 
\end{tabularx}
\caption{Order parameter across a~diagonal line through the domain for the scenario of~\cref{sec:bulkShift} for $p=1$, $\beta=0$~\emph{(solid)}, $\beta=1$~\emph{(dashed)}.
The profiles represent the stationary state controlled by the absolute Newton tolerance~$\mathrm{tol_{abs}}$, cf.~\eqref{eq:Newton:stop}, with values of~$1\E{-16}$ \emph{(inner profile)}, $1\E{-15}$, $\ldots$, $1\E{-8}$ \emph{(outer profile)}.  The zoomed-in plots do not preserve the aspect-ratio.}
\label{fig:bulkShiftPlotsAbsTol}
\end{figure}
While in~\cref{sec:CHmixed:spinodalDecomposition}, we investigated the dynamics of an~evolving binary system for the constant and degenerate mobility case, in this section we perform a~sensitivity study toward a~stationary state solution.  Analytical solutions~$c$ of~\eqref{eq:model} in one-dimensional domains take the value of~$-1$ and $1$ in the bulk phases~\cite[cf.~e.\,g.][]{AristotelousThesis2011}.  Even though these two values minimize the chemical energy density~$\Phi$, it is observed that these bulk values are perturbed by numerical solutions especially in the case of multiple space dimensions, but even in the one-dimensional case~\cite[e.\,g][]{LeeEtAl2016}.  Even though this perturbation---that we call \emph{bulk shift} since all bulk values are either shifted up or down---is bounded by the mesh size (cf.~\cref{sec:convergenceStudy}), in practical simulations, a~false bulk value will imply a~false droplet size due to mass conservation reasons.
This section will give some indication on how the amount of bulk shift depends on certain parameters, namely time step size~$\tau$, mesh size~$h$, droplet size~$L$, and interface parameter~$\kappa$.  The baseline scenario is defined by the choices
$\tau\in \{ 1/1000, 1/2000, 1/4000\}$, $h=1/40$, $L=1/2$, and $\kappa = 1/800$.
\paragraph{Scenario definition}
We consider the following \enquote{droplet in a fluid} scenario in a~\mbox{pseudo-2D} setting for constant and degenerate mobility, and polynomial degrees~$p = 0,1,2$:
Let the initial data be given by $c^0 = 0.95$ for $\vec{x}\in [0.5 - L/2, 0.5+L/2]^2\times [0,h]$ and $c^0 = -0.95$ elsewhere.  We choose the parameters as listed in~\cref{tab:bulkShiftData} and penalty~$\sigma=2^p$.  Values in bold font correspond to the baseline scenario. When we decrease or increase a~particular parameter value with respect to its baseline, the other parameters are kept equal to their baseline values.
This cube of size~$L$ will eventually develop into a~droplet with the expected values of approximately~$1$ inside and~$-1$ outside the droplet.  We could have also used the initial values of $c^0\in\{-1,1\}$, however, the case $p=0$ with degenerate mobility would not develop as there is no element with non-zero mobility, cf.~\eqref{eq:mobility}.  Note that using $c^0\in\{-1,1\}$ instead of $c^0\in\{-0.95,0.95\}$ will yield the same stationary state provided that the size of the rectangle is chosen such that there is the same total amount of mass in the domain.  It is important to note that the absolute nonlinear solver tolerance controls the stationary state, cf.~\eqref{eq:Newton:stopInitial}. Thus, unless explicitly stated so, all simulations use same absolute nonlinear solver tolerance of $\mathrm{tol}_\mathrm{abs}=1\E{-16}$.
\paragraph{Results}
The order parameter~$c_h$ at stationary state was plotted against a~diagonal line across the domain, results are visualized in~\cref{fig:bulkShiftPlots}.  
We observe that the phenomenon of bulk shift occurs only if the interface has a~curvature (cf.~stationary state of~\cref{fig:sim:SD00p1beta0:energyDissipation} for the case of a~flat interface).  Loosely speaking, mass is preferred in droplet-types of phases, i.\,e., is pulled into droplets and thus causes a~shift over~$1$ if the droplet is associated with~$c=1$ and vice versa.
\par
\emph{Variations in $\tau$:}
Let us first investigate the impact of varying~$\tau$ as documented in~\cref{tab:bulkShiftData}. 
For $\beta = 0$, the solutions do not exhibit any notable sensitivity to~$\tau$.
For $\beta = 1$, the bulk shift is slightly reduced for smaller values of~$\tau$.
Clearly, with increasing order of approximation, the solution for~$c$ becomes smoother.  It was observed that the magnitude of the maximum bulk shift is smaller for the degenerate mobility cases in an overall comparison against the constant mobility cases.  Interestingly, among the degenerate mobility cases, the amount of bulk shift is largest for $p=1$ while the profiles for $p=0,2$ are close to each other.
\par
\emph{Variations in~$h$:}
The effect of $h$ is examined next by varying $h$ according to the settings in~\cref{tab:bulkShiftData}.  Overall, the degenerate mobility solutions exhibit relatively more sensitivity to the choice of~$h$.  On the other hand, the magnitude of bulk shift is smaller for the degenerate mobility cases.  For example, $p=0$ cases do not overshoot the value of $1$.  One exception is the $p=1, \beta=1$ case, for which the $c$ profiles appear to be almost identical for $\beta=0$ and $\beta=1$.  This is another registration of the fact that the stable configuration is the product of a complex interplay between mobility, absolute nonlinear solver tolerance~$\mathrm{tol}_\mathrm{abs}$, and the order of approximation. For decreasing~$h$, the bulk shift is always reduced. 
\par
\emph{Variations in~$L$:}
The droplet size, represented by the parameter~$L$, is one of the major parameters that affects the curvature of the profile as it evolves from a non-equilibrium to a stabilized equilibrium configuration.  
The droplet size has a~significant impact on the bulk shift, especially for~\mbox{$(\beta=0)$}. The trend can be summarized as follows: The larger the droplet, the less bulk shift.  The reason is likely to be the smaller curvature of larger droplets.  This behavior is less prominent for \mbox{$(\beta=1)$}.
Smaller droplets \mbox{($L=1/4$)} completely diffuse away for~\mbox{$(\beta=0)$} since the complete removal of gradient energy at the expense of chemical energy seems to be energetically beneficial.  Interestingly, for \mbox{$(\beta=1)$} droplets of that size stay stable.
\par
\emph{Variations in $\kappa$:}
The width of the interface and thus the slope is controlled by~$\kappa$, cf.~\cref{sec:ModelProblem}. 
Especially for~\mbox{$(\beta=0)$}, the bulk shift is significantly reduced for smaller values of~$\kappa$. 
Using degenerate~\mbox{$(\beta=1)$} instead, $\kappa$ has no significant impact on the profile at stationary state except for~$p=0$.  
In general, the value in the center of the droplet exceeds the value of~$c=1$ only for $p=1$ and not for $p=0,2$. 
\par
\emph{Impact of Newton tolerance:} 
The effect of absolute nonlinear solver tolerance~$\mathrm{tol}_\mathrm{abs}$ is investigated in~\cref{fig:bulkShiftPlotsAbsTol}.
The profiles represent the stationary state controlled by~$\mathrm{tol_{abs}}$, cf.~\eqref{eq:Newton:stop}, with values of~$1\E{-16}$ \emph{(inner profile)}, $1\E{-15}$, $\ldots$, $1\E{-8}$ \emph{(outer profile)}. The profiles documented in this figure indicate that the final configuration of the profile strongly depends on~$\mathrm{tol}_\mathrm{abs}$.  The dependence is mildly larger and more complex for the degenerate mobility cases.   
\begin{figure}[t!]
\centering
\begin{tabularx}{\linewidth}{@{}CCCC@{}}
\includegraphics[width=\linewidth]{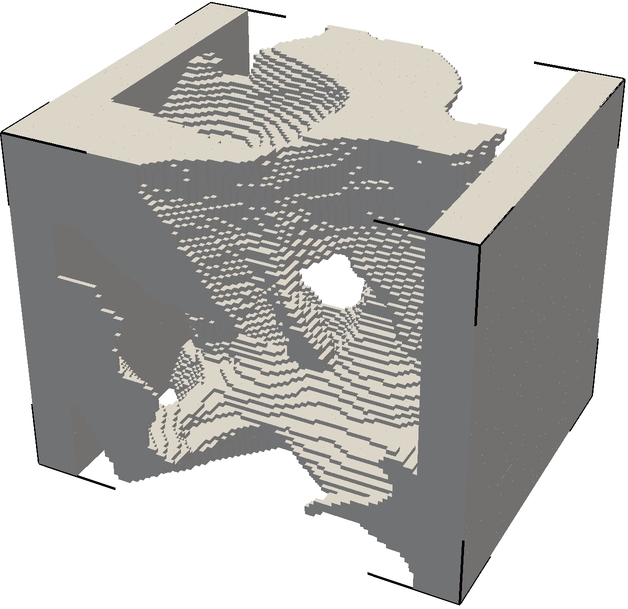} &
\includegraphics[width=\linewidth]{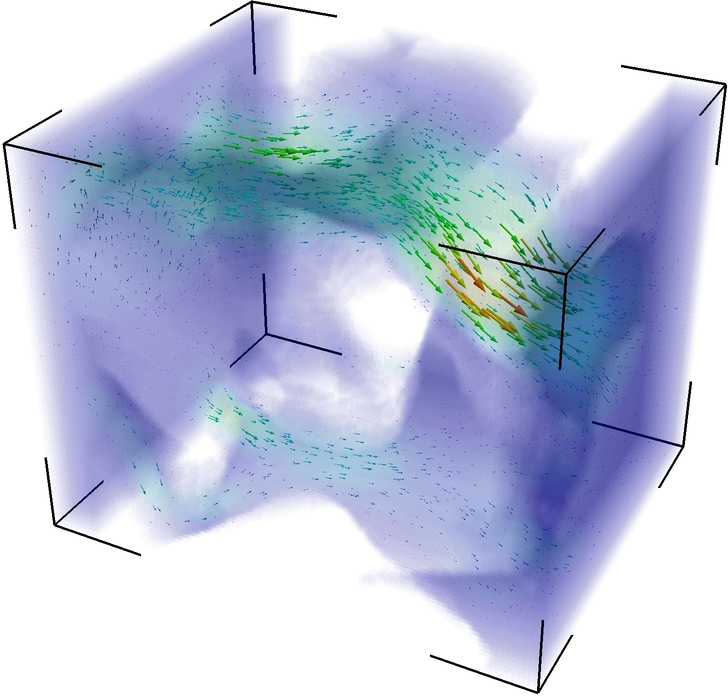} & 
\includegraphics[width=\linewidth]{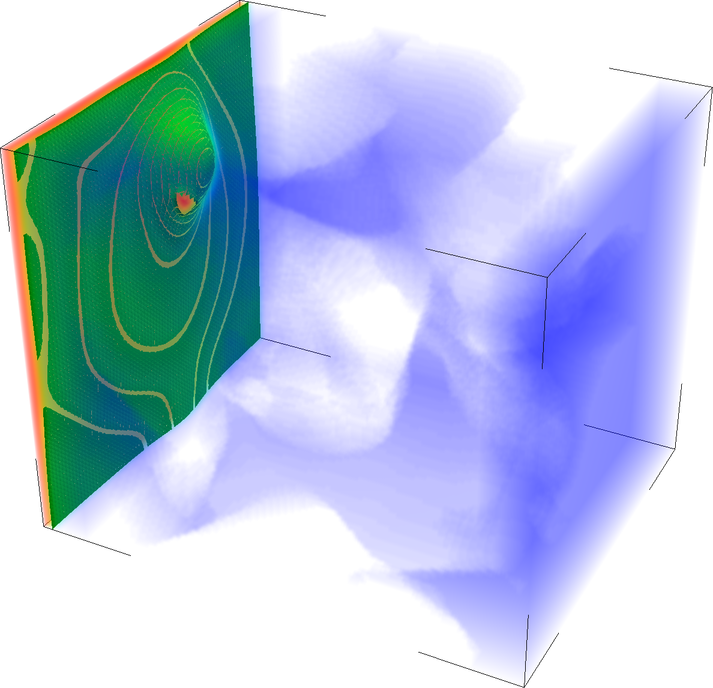} & 
\includegraphics[width=\linewidth]{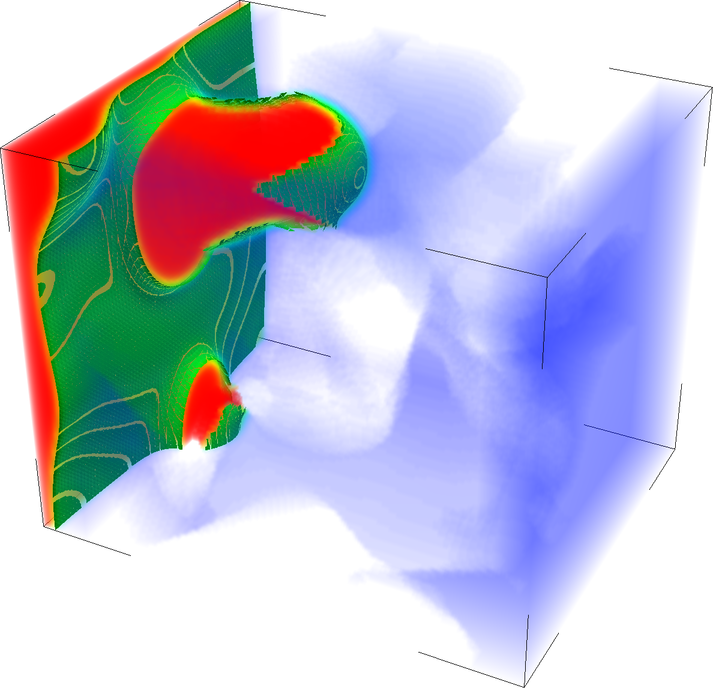} 
\\
\includegraphics[width=\linewidth]{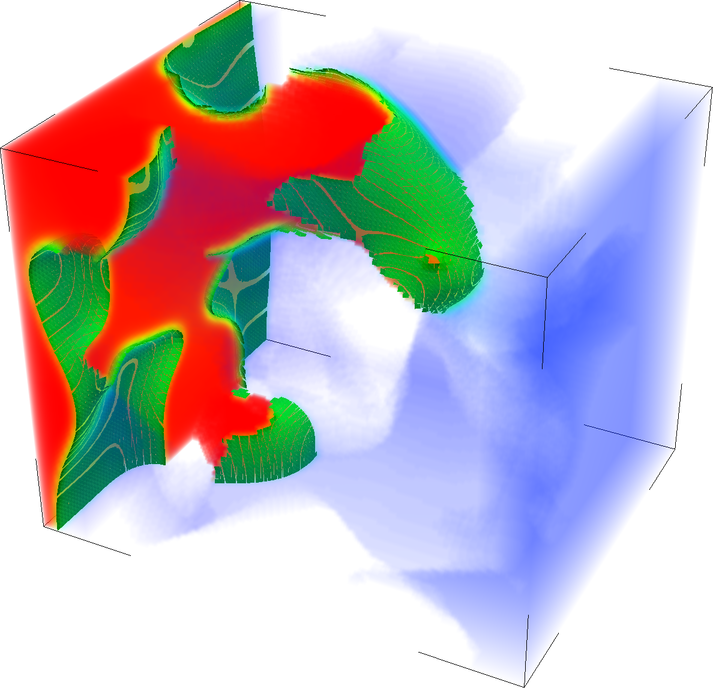} & 
\includegraphics[width=\linewidth]{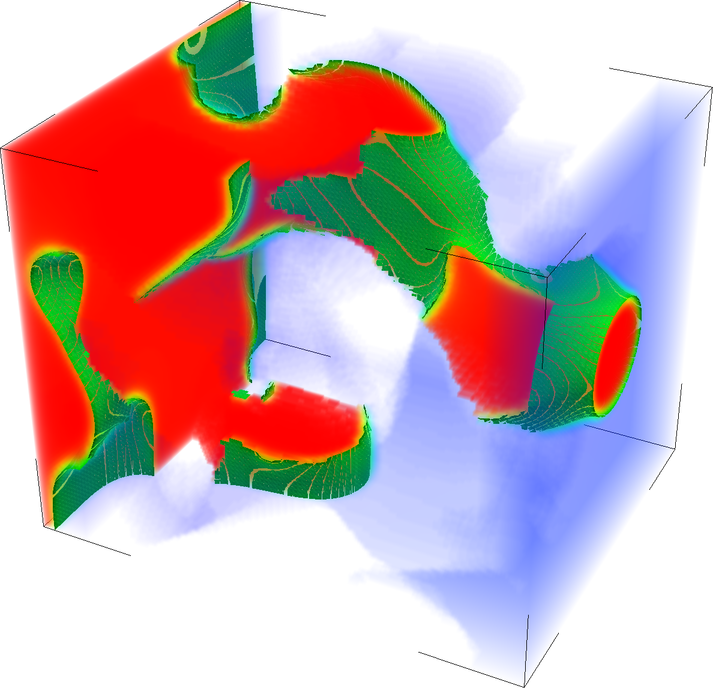} &
\includegraphics[width=\linewidth]{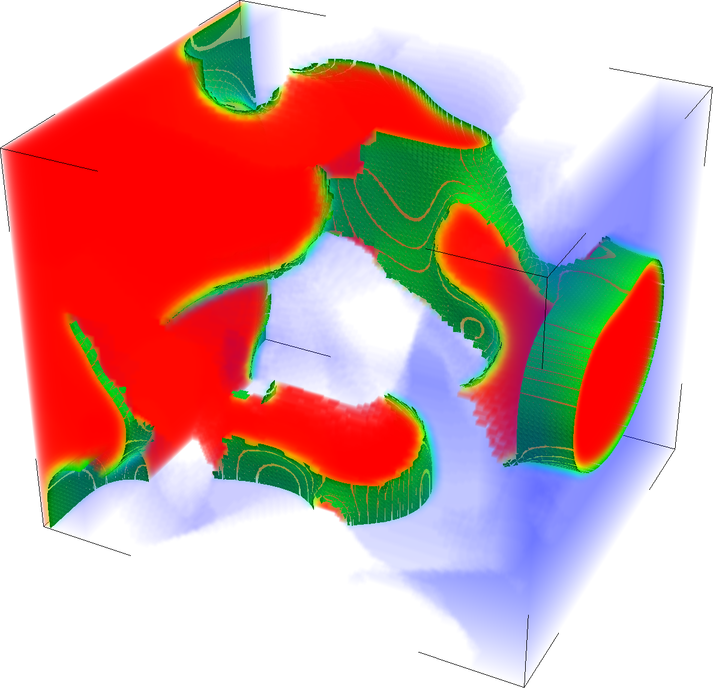} & 
\includegraphics[width=\linewidth]{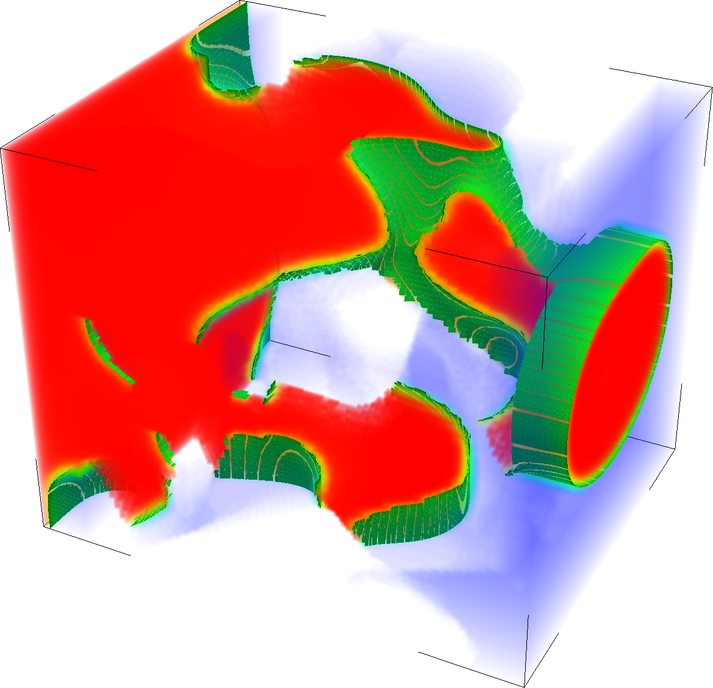} 
\\
\includegraphics[width=\linewidth]{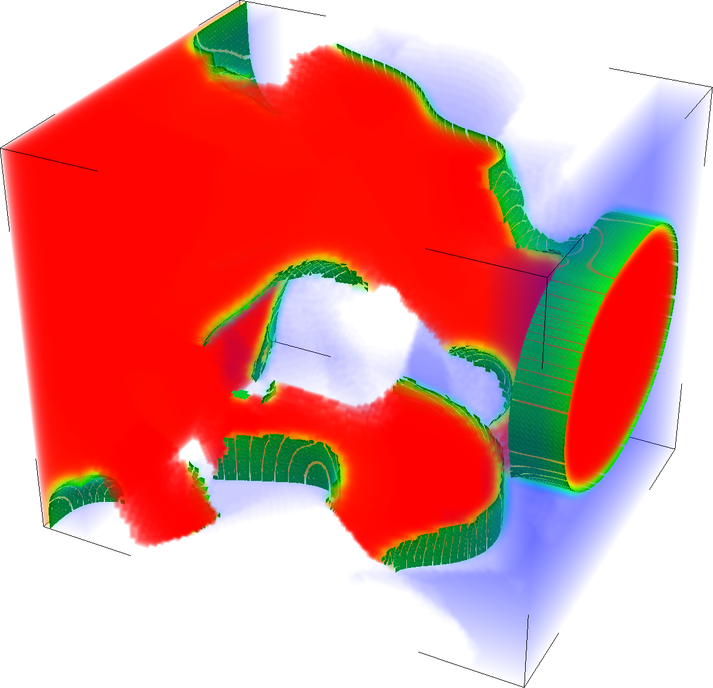} & 
\includegraphics[width=\linewidth]{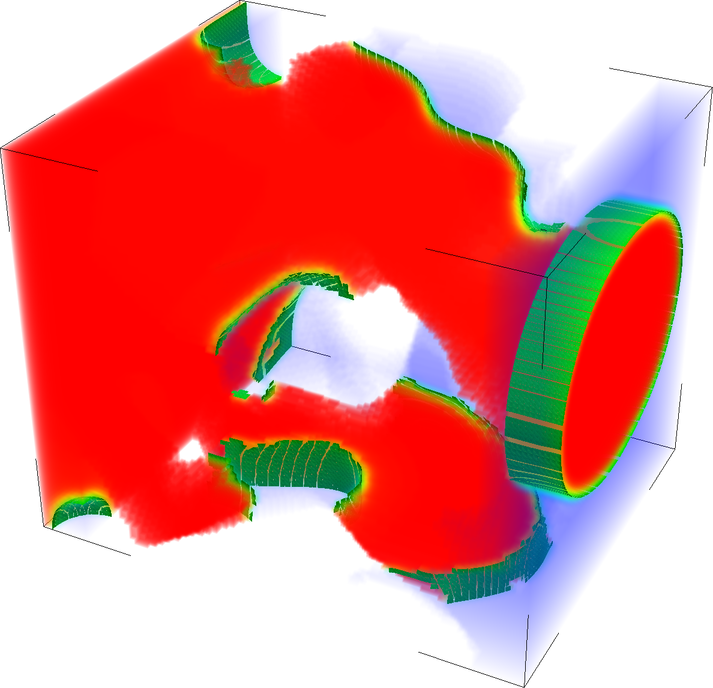} &
\includegraphics[width=\linewidth]{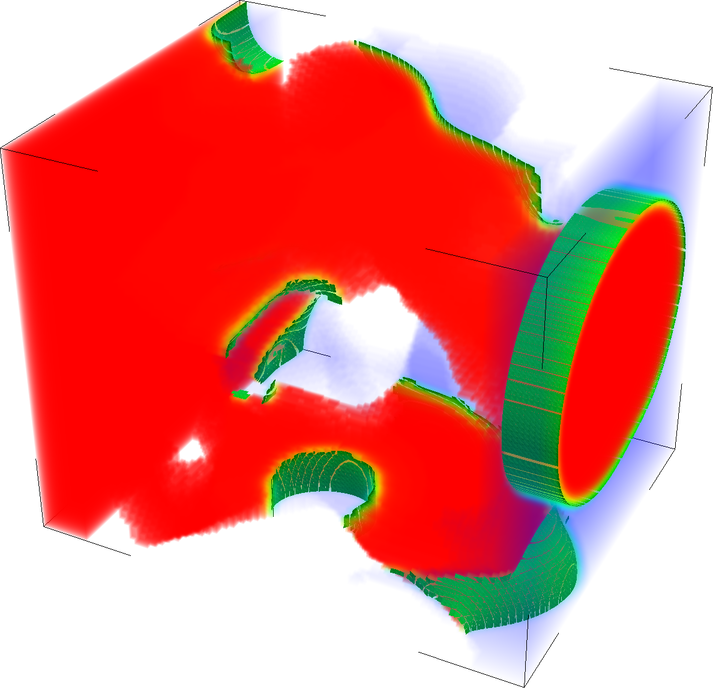} & 
\includegraphics[width=\linewidth]{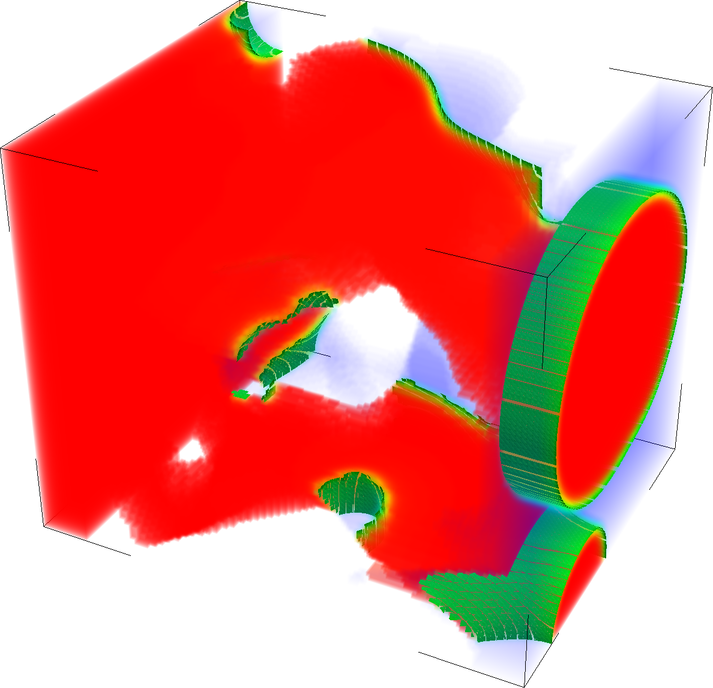} 
\end{tabularx}
\caption{Breakthrough simulation according to~\cref{sec:PorousFlow} using $p=1$, $\beta = 0$.  Illustrated are the computational domain~$\Omega_h$, the stationary flow field, and the order parameter after $100$, $200$, \ldots, $1000$ steps. The colors indicate the phase~\emph{(blue\,/\,transparent \mbox{\normalfont for $c{<}0$}, green \mbox{\normalfont for $c{=}0$}, red \mbox{\normalfont for $c{>}0$})}.}
\label{fig:sim:BreakThrough}
\end{figure}
\subsection{Flow through a porous domain}\label{sec:PorousFlow}
This numerical experiment features immiscible two-phase two-component flow through a~porous domain that is generated by $\muCT$-imaging of a~real sandstone rock sample.
The focus here is the solution of the advective CH~equation subjected to a~realistic flow field.
Thus, in this scenario, we compute a~solenoidal~flow field~$\vec{v}$ by solving the incompressible Navier--Stokes equations and use this velocity in the advective term of the CH~equations~\eqref{eq:model}.  
\par
We consider a~voxel set of $100^3$ voxels containing $248\,309$~elements, the union of which builds~$\Omega_h$.   
The interior part with a~size of~$80^3$ voxels stems from $\muCT$-imaging and represents the porous rock, cf.~\cref{fig:sim:BreakThrough}. 
In order to induce a~flow field~$\vec{v}$ into the pore space, we attach buffers (also called jackets) with a~size of $10\times80\times80$ voxels to two facing sides of~$\Omega_h$ aligned with the main flow direction.
Basically, the buffers allow the flow field to develop, and thereby, permit the imposition of physically realistic boundary conditions at the inlet and outlet~$\Gamma_h^\mathrm{ext}$. 
\par
At initial time, let the domain be saturated by one component, i.\,e.~$c^0 = -1$ in~$\Omega_h$.
The second component is then injected through the inlet buffer. 
\Cref{fig:sim:BreakThrough} shows the evolution of the order parameter through the pore space~$\Omega_h$ as the flow progresses.
The upper channel system is better connected than the lower one. Therefore, it leads to larger velocity magnitudes, which in turn gives rise to effective transport of the phase associated with~$c=1$.  This phase in the upper channel system breaks through earlier and spreads into the outlet buffer, while the lower channel system experiences a~delayed breakthrough.
From the physical viewpoint, the diffuse interfaces deform but do not degenerate or \enquote{diffuse away} while the order parameter~$c$ does not exhibit a~significant bulk shift (cf.~\cref{sec:bulkShift}). 

\section{Conclusion and outlook}\label{sec:Conclusion}
We developed a novel numerical method for the solution of the Cahn--Hilliard equation and its advective variant in~3D. The method is designed to treat both constant and degenerate mobility and specialized for pseudo-Cartesian voxel grids generated via \muCT-imaging of porous materials. A~locally conservative DG method is devised for the spatial discretization of the underlying partial differential equations. A~semi-implicit convex--concave splitting technique is used for the temporal discretization. The resulting nonlinear system of equations are solved by use of a Newton iterative solver. Application of the Schur complement method is designed for the efficient solution of the linear systems arising within the iterative nonlinear solution procedure. The protocol for generating the linear system that arises in the nonlinear solution procedure takes special advantage of the voxel-set nature of the pore-scale phase-separation/transport problems. Owing to the structure of the grid (voxel set), it is possible to reduce the polynomial order to zero, yielding a~cell-centered finite volume approach.
\par
In the absence of the advection term, the numerical method simulates the quasi-static (non-equilibrium to equilibrium state) evolution of two immiscible phases subject to a diffuse interface. When the advection term is included, the algorithm simulates the advective steady-state transport of two immiscible fluid phases (in the absence of inertial effects) with a~diffuse interface separating them. The continuum model for the Cahn--Hilliard equation is based on the minimization of Helmholtz free energy. We demonstrate that the proposed numerical algorithm minimizes the Helmholtz free energy in a~discrete sense. 
\par
The code implementation of the numerical method formulated for the Cahn--Hilliard equation is first validated using numerical convergence tests. Results of these tests signify optimal theoretical convergence rates. The code is then applied to a number of fundamental problems for validation and numerical experimentation purposes. In this context, simulations of the spinodal decomposition phenomenon performed with the degenerate mobility approach are compared to those that use the constant mobility approach.  
A~phenomenon which we call \enquote{bulk shift}---the shift of the numerical solution in either direction---is investigated in a~sensitivity study on the simulation of a~droplet in fluid: while being insensitive to the time step size, a~decrease of interface width or grid size, an~increase of droplet size, or the use of degenerate mobility reduces the amount of bulk shift.
Subject to the working assumptions of the advective Cahn--Hilliard equation, the physical applicability and robustness of the numerical model is demonstrated on a relatively complex computational domain obtained from a \muCT-scan of a real sandstone rock sample.
\par
The demonstrated flexibility of the numerical method in terms of changing the order of the spatial approximation could be used in a~multi-numerics scheme\,/\,\mbox{$p$-adaptive} framework such that the constant bulk areas are approximated with finite volume cells while DG-elements of higher order are selected for the interface area.  Specially designed preconditioners as presented in~\cite{Bansch2010Preconditioning} can be used to further improve performance.  Although not addressed in this paper, DG~methods are suitable for unstructured meshes, local mesh refinement, and parallelization.
\par
A~number of model and numerical parameters are analyzed to evaluate their impact on the droplet in a~fluid problem simulation. In the degenerate mobility case, the absolute nonlinear solver tolerance emerges as an~important factor not only for controlling the evolution but also for the equilibrium state of the droplet.  However, the computational cost is significantly higher than for the constant mobility case, which may render its practical usage a~considerable compromise between computational efficiency vs.~diffusion characteristics one would like to model.
\par
A~breakthrough simulation through a~porous medium demonstrates that a~real-life pore-scale flow problem can be effectively approximated by the numerical scheme proposed in this paper.
If one focuses on the first-order physics that can be captured with one-way coupling of the Cahn--Hilliard and Navier--Stokes equations, this numerical experiment also validates the successful implementation of the advective Cahn--Hilliard solver from the physical viewpoint.
We fully realize that there is significant additional work necessary to render the two-phase flow model more physically realistic through two-way coupling of the advective Cahn--Hilliard and Navier--Stokes equations~\cite[cf.][]{LowengrubTruskinovsky1998} including the capillary effects, which are neglected here. Future work will focus on capturing more realistic flow physics. Another branch of the future work will focus distributed parallel computing~\cite{Grossman2016Analysis,GrossmanEtAl2016Survey} to allow the solution of larger physical problems. 


\begin{table}[H]\centering\footnotesize
\begin{tabularx}{\linewidth}{@{}lL@{}}
\toprule
\textbf{Symbol}      & \textbf{Definition}\\
\midrule
$\avg{\,\cdot\,}$, $\jump{\,\cdot\,}$   & Average and jump of a~quantity on a~face, cf.~\eqref{eq:AverageAndJump}.\\
$\delta_c$            & Functional derivative with respect to~$c$.\\
$\beta$              & Switch between constant and degenerate mobility, cf.~\eqref{eq:mobility}.\\
$c$                   & Order parameter, physically meaningful in~$[-1,1]$.\\
$c_h^n$              & Discrete approximation of~$c$ at time level~$t_n$, $n\in\{0,\ldots,\Nst\}$.\\
$e_{km}$             & $m$th face of the $k$th element~$E_k$ (a square).\\
$E_k$                & $k$th element (rectangular cuboid), $k\in\{0,\ldots,\Nel-1\}$.\\
$\setE_h$            & Set of elements, triangulation of~$\Omega_h$.\\
$F$; $F_h$           & Helmholtz free energy as a~functional in~$c$, cf.~\eqref{eq:HelmholtzFreeEnergy}; discrete version, cf.~\eqref{eq:HelmholtzFreeEnergyDiscrete}.\\
$\Gamma_h^\mathrm{int}$; $\Gamma_h^\mathrm{ext}$; $\Gamma_h^\mathrm{wall}$;   & Set of interior faces; set of exterior faces; set of faces on impermeable wall, cf.~\cref{fig:compDomainsAndIndex}.\\
$h$                  & Mesh size, edge length of each element~$E\in\setE_h$.\\
$J$                  & $\coloneqq (0,T)$, open time interval.\\
$\kappa$             & Interface parameter, cf.~\eqref{eq:HelmholtzFreeEnergy}.\\
$\mu$                & Chemical potential, cf.~\eqref{eq:definitionChemicalPotential}.\\
$M$; $M_{+}$         & Mobility as a function of~$c$; $M_{+}(c)\coloneqq\max\{M(c),0\}$, cf.~\eqref{eq:mobility}.\\
$\normal$; $\normal_e$; $\normal_E$  & Unit normal on~$\partial \Omega$ outward of~$\Omega$; --- on face~$e\in\Gamma_h^\mathrm{int}$ with fixed orientation, cf.~\cref{sec:spaceDiscretization}; --- on~$\partial E$ outward of~$E\in\setE_h$.\\
$\Nel$; $\Nloc$; $\Nst$  & Number of elements; number of local degrees of freedom; number of time steps.\\
$\IN$; $\IN_0$       & Set of natural numbers; set of natural numbers including zero.\\
$\Omega$; $\Omega_h$ & Spatial domain, subset of~$\IR^3$; subset of voxels building spatial domain (pore space).\\
$\partial\Omega$     & Boundary of $\Omega$.\\
$\partial\Omega^\mathrm{in}(t)$  & $\coloneqq \{ \vec{x} \in \partial\Omega \,;\, \vec{v}(t, \vec{x}) \cdot \normal < 0 \}$, inflow boundary.\\
$\partial\Omega^\mathrm{out}(t)$ & $\coloneqq \partial\Omega \setminus \partial\Omega^\mathrm{in}(t)$, outflow and noflow boundary.\\
$p$                  & Polynomial degree, $p\in\IN_0$, cf.~$\IP_p(E)$.\\
$\varphi_{kj}$       & $=\hat{\varphi}_j\circ\vec{F}_k^{-1}$, $j$th hierarchical basis function on~$E_k\in\setE_h$.\\
$\hat{\varphi}_j$    & $=\varphi_{kj}\circ\vec{F}_k$, $j$th hierarchical basis function on~$\hat{E}$.\\
$\IP_p(E)$           & Space of polynomials of degree at most~$p$ on element~$E$, $p\in\IN_0$.\\
$\IP_p(\setE_h)$     & $\coloneqq~\prod\nolimits_{E\in\setE_h}\IP_p(E)$, $p\in\IN_0$.\\
$\Phi$               & Chemical energy density as a function of~$c$, cf.~\eqref{eq:PhiOfc}.\\
$\IR$; $\IR^+$       & Set of real numbers; set of strictly positive real numbers.\\
$\sigma$             & Penalty parameter ($\sigma=2^p$ used in this paper), cf.~\eqref{eq:bilinearform:DG}.\\
$t_n$                & $n$th time level, $n\in\{0,\ldots,\Nst\}$.\\
$\tau_n$             & Time step size of $n$th step, $n\in\{0,\ldots,\Nst\}$.\\
$T$                  & End time.\\
$\vec{v}$            & Advection velocity of the mixture, $\vec{v}: J \times\Omega \rightarrow\IR^3$.\\
$\vec{x}$            & $=\transpose{[x,y,z]}$, space variable in the domain~$\Omega$.\\
\bottomrule
\end{tabularx}
\caption{List of symbols.}
\label{tab:listofsymbols}
\end{table}


\bibliography{references} 
\bibliographystyle{elsarticle-num}


\end{document}
\endinput